%% file: Lavoro201112.tex
\documentclass[11pt,a4paper]{article}
\usepackage{amsmath}
\usepackage{amsfonts}
\usepackage{amssymb}
\usepackage{mathrsfs}
\newtheorem{theorem}{\bf Theorem}[section]
\newtheorem{lemma}[theorem]{\bf Lemma}
\newtheorem{corollary}[theorem]{\bf Corollary}

\newtheorem{remark}[theorem]{\bf Remark}
\newtheorem{proposition}[theorem]{\bf Proposition}
\newtheorem{example}[theorem]{\bf Example}

\def\K{{\mathcal {K}} }

\def\D{{\mathbb{D}} }
\newcommand \Ass[1]{{\bf [#1]}}
\renewcommand \bar[1]{{\overline{#1}}}
\def\Tr{{{\rm Tr}}}
\def \X1n{{{x_{1,n}  }}}
\def \xx1n{{\x_{1,n}}}

\def \o {{\omega}}

\def \d {{\delta}}
\def \l {{\lambda}}
\def \L {\mathscr{L}}

\def \A {\mathscr{A}}

\def \G {{\Gamma}}
\def \s {{\sigma}}

\def \R {{\mathbb {R}}}
\def \N {{\mathbb {N}}}

\def \GG {{\mathbb {G}}}

\def \x {{\xi}}
\def \e {{\varepsilon}}

\def \t {{\tau}}
\def \ll {{\lambda}}

\def \y {{\eta}}

\def \z {{\zeta}}
\def \g {{\gamma}}
\def \O {\mathcal{O}}

\def \tilde {\widetilde}

\def\p{\partial}
\def \P {{\cal{P}}}

\def \dd {{\tilde \delta}}

\def \rnn {{\mathbb {R}}^{N+1}}

\input{macro}

\usepackage[english]{babel}
\begin{document}
\title{Two-sided bounds for degenerate processes with densities supported in subsets of $\mathbb{R}^N$
}
\author{{\sc{Chiara Cinti\thanks{Dipartimento di Matematica, Universit\`{a} di Bologna,
Piazza di Porta S. Donato 5, 40126 Bologna (Italy).
E-mail: cinti@dm.unibo.it},
St\'{e}phane Menozzi \thanks{Universit\'{e} d'Evry Val d'Essonne,
Laboratoire d'Analyse et ProbabilitŽs, 23 Boulevard de France, 91037 Evry Cedex (France).
E-mail: stephane.menozzi@univ-evry.fr} and Sergio Polidoro\thanks{Dipartimento di Scienze Fisiche, Informatiche e Matematiche,
Universit\`{a} di Modena e Reggio Emilia, via Campi 213/b, 41125 Modena (Italy).
E-mail: sergio.polidoro@unimore.it}}}
}
\date{ }
\maketitle

\bigskip

\noindent{\it Abstract:}
We obtain two-sided bounds for the density of stochastic processes satisfying a weak H\"ormander condition. In particular  we consider the cases when the support of the density is not the whole space and when the density has various asymptotic regimes depending on the starting/final points considered (which are as well related to the number of brackets needed to span the space in H\"ormander's theorem).
The proofs of our lower bounds are based on Harnack inequalities for positive solutions of PDEs whereas the upper bounds are derived from the probabilistic representation of the density given by the Malliavin calculus.


\bigskip


\medskip


\normalsize

\section{Introduction}

\setcounter{equation}{0}
\setcounter{theorem}{0} 

We present a methodology to derive two-sided bounds for the density of some $\R^N $-valued  degenerate processes of the form
\begin{eqnarray}
\label{THE_EQUATION}
X_t=x+\sum_{i=1}^n\int_0^t Y_i(X_s)\circ dW_s^i+\int_0^tY_0(X_s)ds
\end{eqnarray}
where the $(Y_i)_{i\in \leftB 0,n\rightB} $ are smooth vector fields defined on $\R^{N}$, $((W_t^i)_{t\ge 0})_{i\in \leftB 1,n\rightB} $ stand for $n$-standard monodimensional independent Brownian motions defined on a filtered probability space $(\Omega,\F,(\F_t)_{t\ge 0},\P) $ satisfying the usual conditions. Also $\circ \, dW_t $ denotes the Stratonovitch integral.
The above stochastic differential equation is associated to the Kolmogorov operator
\begin{equation}
\label{OPERATOR}
  \L  = \tfrac12 \sum_{i=1}^n Y_i^2 +Z,\qquad Z=Y_0- \p_t.
 \end{equation}
We assume that the H\"ormander condition holds:
 \begin{description}
 \item[{\rm \Ass{H}}] $\qquad\qquad
 {\rm {Rank}}( {\rm {Lie}}\{Y_1,\cdots, Y_n ,Z \}(x))=N+1, \quad \forall \, x\in \R^N $.
\end{description}

We will particularly focus on processes satisfying a \emph{weak} H\"ormander condition, that is ${\rm {Rank}}( {\rm {Lie}}\{Y_1,\cdots, Y_n ,-\p_t \}(x))<N+1, \ \forall \, x\in \R^N$.
This means that the first order vector field $Y_0$ (or equivalently the drift term of the SDE) is needed to span all the directions.

As leading examples we have in mind processes of the form
\begin{eqnarray}
\label{PROC}
X_t^i=x_i+W_t^i,\quad \forall i\in \leftB 1,n\rightB, \quad  X_t^{n+1}=x_{n+1}+\int_{0}^t |X_s^{1,n}|^{k} ds,
\end{eqnarray}
where $X_s^{1,n}=(X_s^1,\cdots,X_s^n) $ (and correspondingly for every $x\in\R^{n+1},\ x_{1,n}:=(x_1,\cdots,x_n)$), $k$ is any \emph{even} positive integer and $|.| $ denotes the Euclidean norm of $\R^n $.
Note that we only consider even exponents in \eqref{PROC} in order to keep $Y_0$ smooth. Our approach also applies to
\begin{eqnarray}
\label{PROC_2}
X_t^i=x_i+W_t^i,\quad \forall i\in \leftB 1,n\rightB, \quad  X_t^{n+1}=x_{n+1}+\int_{0}^t \sum_{i=1}^n (X_s^{i})^{k} ds,
\end{eqnarray}
for any given positive integer $k$. 

It is easily seen that the above class of processes satisfies the weak H\"ormander condition.  Also for equation \eqref{PROC},  the density $p(t,x,.)$ of $X_t$ is  supported on $\R^n\times (x_{n+1},+\infty)$ for any $t>0$. Analogously, for equation \eqref{PROC_2}, the support of $p(t,x,.) $ is $\R^{n+1}$ when $k$ is odd and $\R^n\times (x_{n+1},+\infty) $ when  $k$ is  even.

\medskip
Let us now briefly recall some
known results concerning these two examples.
First of all, for $k=1$, equation \eqref{PROC_2} defines a Gaussian process. The explicit expression of the density goes back to Kolmogorov \cite{Kolmogorov} and writes for all $t>0, \ x,\x \in \R^{n+1} $:
\begin{eqnarray}
\label{DENS_K}
p_K(t,x,\x)=\frac{\sqrt 3}{(2\pi)^{\frac{n+1}2} t^{\frac{n+3}2}}\!\exp\!\left (\!-\left\{\! \frac 14\frac{|\x_{1,n}-x_{1,n}|^2}{t}+3\frac{|\x_{n+1}-x_{n+1}-\frac{\sum_{i=1}^n (x_i+\x_i)}{2}t|^2}{t^3}\right\} \right)\!.
\end{eqnarray}
We already observe the two time scales associated respectively to the Brownian motion (of order $t^{1/2}$) and to its integral (of order $t^{3/2}$) which give the global diagonal decay of order $t^{n/2+3/2}$.
The additional term $\frac{x_1+\x_1}{2}t $ in the above estimate is due to the transport of the initial condition by the unbounded drift. We also refer to the works of Cinti and Polidoro \cite{CintiPolidoro1} and Delarue and Menozzi \cite{dela:meno:10} for similar estimates in the more general framework of variable coefficients, including non linear drift terms with linear growth.

For equation \eqref{PROC} and $k=2,\ n=1$, a representation of the density of $X_t $ has been obtained from the seminal works of Kac on the Laplace transform of the integral of the square of the Brownian motion \cite{kac:49}. We can refer to the monograph of Borodin and Salminen \cite{boro:salm:02} for an explicit expression in terms of special functions. We can also mention the work of Tolmatz \cite{Tolmatz} concerning the distribution function of the square of the Brownian bridge already characterized in the early work of Smirnov \cite{Smirnov}. Anyhow, all these explicit representations are very much linked to Liouville type problems and this approach can hardly be extended to higher dimensions for the underlying Brownian motion. Also, it seems difficult from the expressions of \cite{boro:salm:02} to derive explicit quantitative bounds on the density.

Some related examples have been addressed by Ben Arous and L\'eandre \cite{bena:lean:91} who obtained asymptotic expansions for the density on the diagonal for the process $X_t^1=x_1+W_t^1,\ X_t^2=x_2+\int_{0}^t (X_s^1)^{m} dW_s^2+\int_0^t (X_s^1)^k ds $. Various asymptotic regimes are deduced depending on $m$ and $k$. Anyhow, the strong H\"ormander condition is really required in their approach, i.e. the stochastic integral is needed in $X^2$.

From the applicative point of view, equations with quadratic growth naturally appear in some turbulence models, see e.g. the chapter concerning the dyadic model in Flandoli \cite{flan:11}. This model is derived from the formulation of the Euler equations on the torus in Fourier series after a simplification consisting in considering a nearest neighbour interaction in the wave space. This operation leads to consider an infinite system of differential equations whose coefficients have quadratic growth. In order to obtain some uniqueness properties, a Brownian noise is usually added on each component. In the current work, we investigate from  a quantitative viewpoint what can be said for a drastic reduction of this simplified model, that is when considering 2 equations only, when the noise only acts on one component and is transmitted through the system thanks to the (weak) H\"ormander condition.

\medskip

Our approach to derive two-sided estimates for the above examples  is the following. The lower bounds are obtained using local Harnack estimates for positive solutions of $\L u=0 $ with $ \L$ defined in \eqref{OPERATOR}.
Once the Harnack inequality is established, the lower bound for $p(t,x,\x) $ is derived applying it recursively along a suitable path joining $x$ to $\x$ in time $t$. The set of points of the path to which the Harnack inequality is applied is commonly called a \emph{Harnack chain}. For $k=1$ in \eqref{PROC_2} the path can be chosen as the solution to the deterministic controllability problem associated to \eqref{PROC_2}, that is taking the points of the Harnack chain along the path $\gamma $ where
\begin{eqnarray*}
\gamma_i'(s)=\omega_i(s),\ \forall i\in \leftB 1,n \rightB,\ \gamma_{n+1}'(s)=\sum_{i=1}^n \gamma_i(s),\ \gamma(0)=x,\ \gamma(t)=\x.
\end{eqnarray*}
and $\omega: L^2([0,t])\rightarrow \R^n $ achieves the minimum of $\int_0^t|\omega(s)|^2ds $, see e.g. Boscain and Polidoro \cite{BoscainPolidoro}, Carciola \textit{et al.} \cite{CarciolaPascucciPolidoro} and Delarue and Menozzi \cite{dela:meno:10}.

In the more general case $k>1$ it is known that uniqueness fails for the associated control problem, i.e. when $\gamma_{n+1}'(s)=\sum_{i=1}^n (\gamma_i(s))^k $ in the above equation (see e.g. Tr\'elat \cite{Trelat}). Therefore, there is not a single natural choice for the path $\g$. Actually, we will consider suitable paths in order to derive homogeneous two-sided bounds. After the statement of our main results, we will see in Remark \ref{REM_TRELAT} that the paths we consider allow to obtain a cost similar to the one found in \cite{Trelat} for the abnormal extremals of the value function associated to the control problem.

Anyhow, the crucial point in this approach is to obtain a Harnack inequality \emph{invariant} w.r.t. scale and translation.
Introducing for all $(m,x)\in \N^*\times \R^{n+1}$ the space $V_m(x):=\{\bigl( (Y_{i_1})_{i_1\in \leftB 1,n\rightB}$, $([Y_{i_1},Y_{i_2}](x))_{(i_1,i_2)\in \leftB 0,n\rightB^2}, \cdots, ([Y_{i_1},[Y_{i_2},\cdots,[Y_{i_{m-1}},Y_{i_m}]]](x))_{(i_1,\cdots,i_m)\in \leftB 0,n\rightB^m} \bigr) \}$, the above invariance properties imply that ${\rm dim}({\rm Span} \{ V_{m}(x)\})$ does not depend on $x$ for any $m$.
This property fails for $k>1$ since we need exactly $k$ brackets to span the space at $x=(0_{1,n},x_{n+1})$ and exactly one bracket elsewhere.
Hence, we need to consider a \emph{lifting} procedure of $\L$ in \eqref{OPERATOR} introduced by Rotschild and Stein \cite{RothschildStein} (see also Bonfiglioli and Lanconelli \cite{BonfiglioliLanconelli1}).
Our strategy then consists in obtaining an invariant Harnack inequality for the lifted operator $\tilde \L $. We then conclude applying the previous Harnack inequality to $\L$-harmonic functions (which are also $\tilde \L $-harmonic).
A first attempt to achieve the whole procedure to derive a lower bound for \eqref{PROC_2} and odd $k $ can be found in Cinti and Polidoro \cite{CintiPolidoro3}.

\medskip


Concerning the upper bounds, we rely on the representation of the density of $p$ obtained by the Malliavin calculus. We refer to Nualart \cite{nual:98} for a comprehensive treatment of this subject. The main issues then consist in controlling the tails of the random variables at hand and the $L^p$ norm of the Malliavin covariance matrix for $p\ge 1$.  The tails can be controlled thanks to some fine properties of the Brownian motion or bridge and its local time. The behavior of the Malliavin covariance matrix has to be carefully analyzed introducing a dichotomy between the case for which the final and starting points of the Brownian motion in \eqref{PROC}-\eqref{PROC_2} are close to zero w.r.t. the characteristic time-scale, i.e. $|\X1n|\vee |\xx1n|\le Kt^{1/2} $ for a given $K>0$, which means that the non-degenerate component is in diagonal regime, and the complementary set.
In the first case, we will see that the characteristic time scales of the system \eqref{PROC}, \eqref{PROC_2} and the probabilistic approach to the proof of H\"ormander theorem, see e.g. Norris \cite{norr:86} will lead to the expected bound on the Malliavin covariance matrix whereas in the second case a more subtle analysis is required in order to derive a diagonal behavior of the density similar to the Gaussian case \eqref{DENS_K}.  Intuitively, when the magnitude of either the starting or the final point of the Brownian motion is above the characteristic time-scale, then only one bracket is needed to span the space and the Gaussian regime prevails in small time.

\medskip

Note that our procedure can be split in two steps. In the first one, purely PDEs methods provide us with lower bounds of the density $p$. In this part useful information about its asymptotic behavior in various regimes are obtained by elementary arguments. Once the lower bounds have been established, we rely on some \emph{ad hoc} tools of the Malliavin calculus to prove the analogous upper bounds. However, aiming at improving the readability of our work, we reverse our exposition: we first prove the upper bounds, as well as the diagonal ones, by using probabilistic methods, then we prove the lower bounds by PDEs arguments.

The article is organized as follows. We state our main results in Section \ref{Results}. We then recall some basic facts of Malliavin calculus in Section \ref{GLIMPSE_SECTION} and obtain the upper bounds as well as a diagonal lower bound in Gaussian regime in Section \ref{MALLIAVIN_SECTION}.In Section \ref{PDE_SECTION}, we recall some aspects of abstract potential theory needed to derive the invariant Harnack inequality.  We also give a geometric characterization of the set where the inequality holds.
Section \ref{CHAINS_SECTION} is devoted to the proof of the lower bounds. 

\section{Main Results}
\label{Results}
\setcounter{equation}{0}
\setcounter{theorem}{0}

Before giving the precise statement of our bounds for the the density $p$ of $X$ in \eqref{PROC} or \eqref{PROC_2}, we give some remarks. In the sequel $p(t,x,.)$ stands for the density of any stochastic process $X$ at time $t$ starting from $x$. It is well known that, if the vector fields $Y_1, \dots, Y_n$ (note that the drift term $Y_0$ does not appear) satisfy the H\"ormander condition, then the following two sided bound holds:
\begin{equation} \label{SUM_OF_SQUARES}
    p(t,x,\x)\asymp\frac{1}{\sqrt{\text{vol}\left( B_Y(x, t) \right)}} \exp \left (- \frac{d_Y(x, \x)^2}{t} \right).
\end{equation}
Here and in the sequel, for measurable functions $g:\R^{+*}\times \R^n \rightarrow \R,h:\R^{+*}\times \R^{2n}$, the above notation $p(t,x,\x)\asymp \frac1{g(t,x)} \exp(-h(t,x,\x))$ means that there exists a constant $C\ge 1$ s.t.
\begin{equation}
\label{DEF_ASYMP}
\frac{C^{-1}}{g(t,x)}\exp(-Ch(t,x,\x))\le  p(t,x,\x)\le \frac{C}{g(t,x)}\exp(-C^{-1}h(t,x,\x)).
\end{equation}
Moreover, in \eqref{SUM_OF_SQUARES}, $d_Y$ denotes the Carnot-Carath\'eodory distance associated to $Y_1, \dots, Y_n$, and $B_Y(x,r)$ is the relevant metric ball, with center at $x$ and radius $r$. On the other hand \eqref{DENS_K} shows that, when the drift term $Y_0$ is needed to check the H\"ormander condition, the density $p$ of the process $X$ doesn't satisfy \eqref{SUM_OF_SQUARES}. In this article we prove that, when considering processes \eqref{PROC} and \eqref{PROC_2} with $k>1$, different asymptotic behavior as $|x| \to + \infty$ appear.



To be more specific, we first remark that a behavior similar to \eqref{DENS_K} can also be observed for equations \eqref{PROC} and \eqref{PROC_2}.Conditioning w.r.t. to the non degenerate component we get
\begin{eqnarray*}
    p(t,x,\x) = p_{X^{1,n} }(t,x_{1,n},\x_{1,n}) \, p_{X^{n+1}}(t,x_{n+1},\x_{n+1}|X_0^{1,n}=x_{1,n},X_t^{1,n}=\x_{1,n}),
\end{eqnarray*}
where $p_{X^{1,n} }(t,x_{1,n},\x_{1,n})=\frac{1}{(2\pi t)^{n/2}} \exp(-\frac{|\x_{1,n}-x_{1,n}|^2}{2t})$ is the usual Gaussian density,
\begin{eqnarray*}
p_{X^{n+1}}(t,x_{n+1},\x_{n+1}|X_0^{1,n}=x_{1,n},X_t^{1,n}=\x_{1,n})=p_{Y_t}(\x_{n+1}-x_{n+1}),\\ \\
Y_t :=\begin{cases}\int_0^t |\frac{t-s}{t}x_{1,n}+\frac st\x_{1,n}+W_s^{0,t} |^kds\ \text{for}\ \eqref{PROC},\\
\int_0^t\sum_{i=1}^n (\frac{t-s}{t}x_{i}+\frac st\x_{i}+W_s^{0,t,i} )^kds \ \text{for}\ \eqref{PROC_2}, \end{cases}
\end{eqnarray*}
and $(W_s^{0,t})_{u\in [0,t]} $ stands for the standard $d$-dimensional Brownian bridge on $[0,t] $, i.e. starting and ending at 0.

For the sake of simplicity, we next focus on the case $n=1$ and $k=2$ so that \eqref{PROC} and \eqref{PROC_2} coincide. Moreover we assume $x_1 = \x_1$.
This leads to estimate the density of:
\begin{eqnarray}
\label{SVIL}
Y_t:=\int_0^t (x_1 + W_s^{0,t})^{2}ds=  t x_1^2 + 2 x_1 \int_0^t W_s^{0,t} ds + \int_0^t  (W_s^{0,t})^{2} ds.
\end{eqnarray}
Thus, when $|x_1|$ is sufficiently big w.r.t. the characteristic time scale $t^{1/2}$, the Gaussian  random variable
$$G:= 2 x_1\int_0^t W_s^{0,t} ds \overset{({\rm law})}{=}{\cal N}\left(0, \tfrac{t^3}3 |x_1|^{2}\right)$$ dominates in terms of fluctuation order w.r.t. the other random contribution whose variance behaves as $O(t^4)$  \footnote{The previous identity in law is derived from It\^o's formula and the differential dynamics of the Brownian bridge. Namely, $\int_0^t W_s^{0,t}ds=\left\{-(t-s)W_s^{0,t}|_{0}^t+\int_0 ^t (t-s)\left(-\frac{W_s^{0,t}}{t-s} ds+dW_s\right)  \right\}\iff \int_0^tW_s^{0,t}ds=\frac12\int_0^t (t-s) dW_s$.}.

If we additionally assume that $|\x_{2}-x_{2}-t x_{1}^2|\le \bar C |x_{1}|t^{3/2}  $, for some constant $\bar C:=\bar C(n=1,k=2)$ to be specified later on, that is the deviation from the deterministic system deriving from \eqref{PROC}, obtained dropping the Brownian contribution, has the same order as the standard deviation of $G$, we actually find:
\begin{eqnarray*}
 p(t,x,\x)\asymp\frac{1}{|x_{1}| t^{\frac{1+3}2}}\!\exp\!\left (\!-\left\{\! \frac{|\x_{1}-x_{1}|^2}{t}+\frac{|\x_{2}-x_{2}- x_1^{2}t|^2}{|x_{1}|^{2}t^3}\right\}\right)\!.
\end{eqnarray*}

When $|\x_{2}-x_{2}-t x_{1}^2| > \bar C |x_{1}|t^{3/2}$, that is when the deviation from the deterministic system exceeds a certain constant times the standard deviation, the term $\int_0^t(W_s^{0,t})^2 ds$  in \eqref{SVIL} is not negligeable any more and we obtain the following heavy-tailed estimate:
\begin{eqnarray*}
 p(t,x,\x)\asymp\frac{1}{ t^{\frac 12+2}}\!\exp\!\left (\!-\left\{\! \frac{|\x_{1}-x_{1}|^2}{t}+\frac{|\x_{2}-x_{2}- x_1^{2}t|}{t^{2}}\right\}\right)\!.
\end{eqnarray*}
The diagonal contribution of the degenerate component corresponds to the intrinsic scale of order $t^{2}$ of the term $\int_0^t (W_s^{0,t})^{2} ds $. In particular, if $x_{1,n}=0_{1,n}$ this is the only random variable involved. The off-diagonal bound can be explained by the fact that $\int_0^t (W_s^{0,t})^{2} ds$ belongs to the Wiener chaos of order $2$. The tails of the distribution function for such random variables can be characterized, see e.g. Janson \cite{jans:97}, and are homogeneous 
to the non Gaussian term in the above estimate.

Observe also that the density $p$ is supported on the half space $\{\xi \in \R^{2}: \xi_{2}>x_{2}\} $. We obtain as well an asymptotic behavior for the density close to the boundary. Precisely, for $0<\xi_{2}-x_{2}$ sufficiently small w.r.t. to the characteristic time-scale $t^{2}$ of $\int_0^t (W_s^{0,t})^2ds $, that is when the deviations of the degenerate component have the same magnitude as those of the highest order random contribution,  then
\begin{eqnarray*}
 p(t,x,\x)\asymp\frac{1}{ t^{1/2+2}}\!\exp\!\left (\!-\left\{\!\frac{|x_{1}|^{4}+|\x_{1}|^{4}}{\x_{2}-x_2} +\frac{t^{2}}{(\x_{2}-x_{2})}\right\}\right)\!.
\end{eqnarray*}

We summarize the above remarks with the assertion that processes of the form \eqref{PROC} or \eqref{PROC_2} do not have a single regime for $k>1$.The precise statements of the previous density bounds are formulated for general $n$ and $k$ in the following Theorem \ref{MTHM}.







\begin{theorem}
\label{MTHM}
Let $x=(x_{1,n},x_{n+1})\in \R^{n+1}$, and $\x=(\x_{1,n},\x_{n+1})\in \R^n\times (x_{n+1},+\infty) $ for $k$ even,
and $\x\in \R^{n+1} $ for $k$ odd, be given.
Define $$\Psi(x_{1,n},\x_{1,n}):=\begin{cases}
|x_{1,n}|^k+|\x_{1,n}|^k
,\  \text{for \ \eqref{PROC}},\\
\sum_{i=1}^n \{(x_i)^k+(\x_i)^k\},\  \text{for \ \eqref{PROC_2}}.
\end{cases}$$

\begin{itemize}
\item[{i)}] Assume $\frac{|\x_{n+1}-x_{n+1}-ct(|\X1n|^k+|\xx1n|^k)|}{t^{3/2}(|\X1n|^{k-1}+|\xx1n|^{k-1})}\ge \bar C $ where $c:=c(k)=2+\frac{2^{k-1}}{k+1}$ and $\bar C$ is fixed. Then there exists a constant $C_1:=C_1(n,k,\bar C)\ge 1$ s.t. for every $t>0$,
\begin{eqnarray}
\label{STIMA_0}
\frac{C_1^{-1}}{t^{\frac{n+k}2+1}}\exp \!\Big( \!\! -C_1 I(t,x,\x,\frac{1}{2^{k+4}})\! \Big) \! \!\!\! & \!\!\! \!\!\! \le \!\!\! & \!\!\!\!\!\! p(t,x, \x)  
\le \!\frac{C_1}{t^{\frac{n+k}2+1}} \exp\! \Big( \!\! -C_1^{-1} I(t,x,\x,\frac{2^{k-1}}{k+1})\Big),\!\!  \\
\forall c\in \R^+,\  I(t,x,\x,c)&:=& \frac{|\x_{1,n}-x_{1,n}|^2}t + \frac{|\x_{n+1}-x_{n+1}-c\Psi(x_{1,n},\x_{1,n})t|^{2/k}}{t^{1+2/k}}\nonumber.
\end{eqnarray}

	\item[{ii)}]  Assume $\frac{|\x_{n+1}-x_{n+1}-ct(|\X1n|^k+|\xx1n|^k)|}{t^{3/2}(|\X1n|^{k-1}+|\xx1n|^{k-1})}\le \bar C $ (with $c,\bar C$ as in point \textit{i)}) and $|x_{1,n}|\vee |\x_{1,n}|/ t^{1/2}\ge K$, with $K$ sufficiently large. Then, there exists $C_2:=C_2(n,k,K,\bar C)\ge 1 $ s.t. for every $t>0$:
\begin{eqnarray}
\label{STIMA_1}
\frac{C_2^{-1}\exp(-C_2I(t,x,\x))}{(|x_{1,n}|^{k-1}+|\x_{1,n}|^{k-1})t^{\frac{n+3}2}} \!&\le& p(t,x, \x)  
\le \! \frac{C_2 \exp(-C_2^{-1}I(t,x,\x))}{(|x_{1,n}|^{k-1}+|\x_{1,n}|^{k-1})t^{\frac{n+3}2}},\\
I(t,x,\x)&:=& \frac{|\x_{1,n}-x_{1,n}|^2}t + \frac{|\x_{n+1}-x_{n+1}- \Psi(x_{1,n},\x_{1,n})t |^2}{(|x_{1,n}|^{(k-1)}+|\x_{1,n}|^{(k-1)} )^2t^3}.\nonumber
\end{eqnarray}
	\item[{iii)}]  For $t>0$, assume $|\x_{n+1}-x_{n+1}|\le Kt^{1+k/2} $ for sufficiently small $K$. Then, there exists $C_3:=C_3(n,k,K)\ge 1$ s.t. we have:
\begin{eqnarray}
\label{STIMA_2}
\frac{C_3^{-1}}{t^{\frac{n+k}2+1}}\exp(-C_3 I(t,x,\x) ) &\le& p(t,x, \x)  
\le \!\frac{C_3}{t^{\frac{n+k}2+1}} \exp(-C_3^{-1}I(t,x,\x)),\!\!\nonumber\\
I(t,x,\x)&:=& \frac{|x_{1,n}|^{2+k}+|\x_{1,n}|^{2+k}}{|\x_{n+1}-x_{n+1}|} +\frac{t^{1+2/k}}{|\x_{n+1}-x_{n+1}|^{2/k}} .
\end{eqnarray}
\end{itemize}
\end{theorem}

\begin{remark} \label{REM-PDE}
It is known that the densities of the stochastic systems \eqref{PROC} and \eqref{PROC_2} are Fundamental Solutions of the Kolmogorov operators
\begin{equation}\label{K-PROC}
	\L = \frac12 \Delta_{x_{1,n}} + |x_{1,n}|^k \p_{x_{n+1}} - \p_t,
\end{equation}
and
\begin{equation}\label{K-PROC_2}
	\L = \frac12 \Delta_{x_{1,n}} + \sum_{j=1}^n x_{j}^k \p_{x_{n+1}} - \p_t,
\end{equation}
respectively. Then, from the PDEs point of view, Theorem \ref{MTHM} provides us with estimates analogous to those due to Nash, Aronson and Serrin for uniformly parabolic operators.
\end{remark}

We next give some comments about our main result. As already pointed out, processes of the form \eqref{PROC} or \eqref{PROC_2} do not have a single regime anymore for $k>1$. Let us anyhow specify that when $C^{-1}\sqrt t\le |x_{i}|\le C\sqrt t,\ \forall i\in \leftB 1,n \rightB, \ C\ge 1 $, then expanding $Y_t$ as in \eqref{SVIL}, we find that all the terms have the same order and thus a global estimate of type \eqref{STIMA_0} (resp. of type \eqref{STIMA_1}) holds for the upper bound (resp. lower bound) in both cases \eqref{PROC} and \eqref{PROC_2}. Observe also that in this case \eqref{STIMA_0} and \eqref{STIMA_1} give the same global diagonal decay of order $t^{(k+n)/2+1} $.

\begin{remark}
\label{REM_TRELAT}
As already mentioned in the introduction, for $k=2, n=1$, we observe from \eqref{STIMA_2} that the off-diagonal bound is homogeneous to the asymptotic expansion of the value function associated to the control problem at its abnormal extremals, see Example 4.2 in \cite{Trelat}. The optimal cost is asymptotically equivalent to $\frac 14 \frac{\x_1^4}{\x_2} $ when $x=(0,0)$ as $\x$ is close to $(0,0)$.
\end{remark}

\begin{remark}
Fix $|\xi_{n+1}-x_{n+1}| $ small, $t\in [K^{-1}|\xi_{n+1}-x_{n+1}|^{2-\varepsilon},K|\xi_{n+1}-x_{n+1}|^{2-\varepsilon}]$ for given $K\ge 1,\varepsilon>0$. We then get from \eqref{STIMA_2} that there exist $\tilde c:=\tilde c(n,k),\tilde C:=\tilde C(n,k,T) $ s.t. $p(t,x,\x)\le \tilde C\exp(-\tilde c/|\xi_{n+1}-x_{n+1}|^\varepsilon)$. This estimate can be compared to the exponential decay on the diagonal proved by Ben Arous and L\'eandre in \cite[Theorem 1.1]{bena:lean:91}.
\end{remark}

 \setcounter{equation}{0}
\setcounter{theorem}{0}
\label{MALLIAVIN}

\section{A Glimpse of Malliavin Calculus}
\label{GLIMPSE_SECTION}
\subsection {Introduction}
Introduced at the end of the 70s by Malliavin, \cite{malli:76}, \cite{malli:78}, the stochastic calculus of variations, now known as Malliavin  calculus, turned out to be a very fruitful tool. It allows to give probabilistic proofs of the celebrated H\"ormander theorem, see e.g. Stroock \cite{stro:83} or Norris \cite{norr:86}. It also provides a quite natural way to derive density estimates for degenerate diffusion processes. The most striking achievement in this direction is the series of papers by Kusuoka and Stroock, \cite{kusu:stro:84}, \cite{kusu:stro:85}, \cite{KusuokaStroock}. Anyhow, in those works the authors always considered ``strong" H\"ormander conditions, that is the underlying space is assumed to be spanned by brackets involving only the vector fields of the diffusive part. For the examples \eqref{PROC}, \eqref{PROC_2} we consider, this condition is not fulfilled. Anyhow a careful analysis of the Malliavin covariance matrix will naturally lead to the upper bounds of Theorem \ref{MTHM} and also to a Gaussian lower bound, when the initial or final point of the non-degenerate component is ``far" from zero w.r.t. the characteristic time scale on the compact sets of the underlying metric, see point \textit{ii)} of Theorem \ref{MTHM}.

We also point out that because of the non uniqueness associated to the deterministic control problem, the strategy of \cite{dela:meno:10} relying on a stochastic control representation of the density breaks down. For the systems handled in \cite{dela:meno:10}, we refer to Bally and Kohatsu-Higa for a Malliavin calculus approach \cite{ball:koha:10}. The Malliavin calculus remains the most robust probabilistic approach to density estimate in the degenerate setting.

We now briefly state some facts and notations concerning the Malliavin calculus that are needed to prove our results.
We refer to the monograph of Nualart \cite{nual:95}, from which we borrow the notations, or Chapter 5 in Ikeda and Watanabe \cite{iked:wata:89}, for further details.

\subsection{Operators of the Malliavin Calculus}
\label{OP_M}

Let us consider an $ n$-dimensional Brownian motion $W$ on the filtered probability space $(\Omega,\F,(\F_t)_{t\ge 0},\P) $ and a given $T>0$.
Define for $h\in L^2(\R^+,\R^n), \ W(h)=\int_0^T \langle h(s),dW_s\rangle $. We denote by ${\cal S} $ the space of simple functionals of the Brownian motion $W$, that is the subspace of $L^2(\Omega,\F,\P) $ consisting of real valued random variables $F$ having the form
$$F=f\bigl(W(h_1),\cdots, W(h_m)\bigr),$$
for some $m\in \N, h_i\in L^2(\R^+,\R^n)$, and where $f:\R^m\rightarrow \R$ stands for a smooth function with polynomial growth.\\

\textit{Malliavin Derivative.}\\

For $F\in{\cal S}$, we define the Malliavin derivative $(D_t F)_{t\in [0,T]}$ as the $\R^n $-dimensional (non adapted) process
$$D_t F=\bsum{i=1}^m \partial_{x_i} f\bigl(W(h_1),\cdots, W(h_m)\bigr) h_i(t).
$$
For any $q\ge 1$, the operator $D: {\cal S}\rightarrow L^q(\Omega, L^2(0,T))$ is closable. We denote its domain by $\D^{1,q} $ which is actually the completion of ${\cal S} $ w.r.t. the norm
$$\|F \|_{1,q}:=\left\{\E[|F|^q]+\E[|DF|_{L^2(0,T)}^q] \right\}^{1/q}.$$
Writing $D_t^j F$ for the $j^{{\rm{th}}} $ component of $D_tF$, we define the $k^{{\rm th}}$ order derivative as the random vector on $[0,T]^k\times \Omega $ with coordinates:
$$D_{t_1,\cdots, t_k}^{j_1,\cdots, j_k}F:=D_{t_k}^{j_k}\cdots D_{t_1}^{j_1} F.$$
We then denote by $\D^{N,q} $ the completion of ${\cal S} $ w.r.t. the norm
$$\|F \|_{N,q}:=\left\{\E[|F|^q]+\sum_{k=1}^{N}\E[|D^kF|_{L^2\bigl((0,T)^k\bigr)}^q] \right\}^{1/q}.   $$
Also, $\D^\infty:=\cap_{q\ge 1}\cap_{j\ge 1} \D^{j,q}$. In the sequel we agree to denote for all $q\ge 1,\ \|F\|_{q}:=\E[|F|^q]^{1/q} $.\\

\textit{Skorohod Integral.}\\

We denote by ${\cal P}$ the space of simple processes, that is the subspace of $L^2([0,T]\times \Omega,\F\times {\cal B}([0,T]),dt\otimes d\P) $ consisting of $\R^n $ valued processes processes  $(u_t)_{t\in [0,T]}$  that can be written
$$u_t=\sum_{i=1}^m  F_i(W(h_1),\cdots ,W(h_m))h_i(t),$$
for some $m\in \N$, where the $(F_i)_{i\in \leftB 1,m\rightB}$ are smooth real valued functions with polynomial growth, $\forall i\in \leftB 1,m\rightB,\ h_i \in L^2([0,T],\R^n)$ so that in particular $F_i(W(h_1),\cdots ,W(h_m))\in {\cal S}$.

Observe also that with previous definition of the Malliavin derivative for $F\in {\cal S}$ we have $(D_sF)_{s\in [0,T]}\in  {\cal P} $.
For $u\in {\cal P}$ we define the Skorohod integral
$$\delta (u):=\sum_{i=1}^m \biggl\{ F_i(W(h_1),\cdots,W(h_m))W(h_i)-\sum_{j=1}^{m}\partial_j F_i(W(h_1),\cdots,W(h_m)) \langle h_i,h_j\rangle_{L^2([0,T]}\biggr\}, $$
so that in particular $\delta (u)\in {\cal S}$. The Skorohod integral is also closable. Its domain writes ${\rm Dom }(\delta):=\{u\in L^2([0,T]\times \Omega): \exists (u_n)_{n\in {\cal P}}, \ u_n\overset{L^2([0,T]\times \Omega)}{\underset{n}{\longrightarrow}} u,\ \delta(u_n)\overset{L^2(\Omega)}{\underset{n}{\longrightarrow}} F:=\delta(u) \}$.\\

\textit{Ornstein Uhlenbeck operator.}\\

To state the main tool used in our proofs, i.e. the integration by parts formula in its whole generality, we need to introduce a last operator. Namely, the Ornstein-Uhlenbeck  operator $L$ which for $F\in {\cal S} $ writes:
\begin{eqnarray*}
LF&:=&\delta(DF) =\langle \nabla f\bigl( W(h)\bigr), W(h) \rangle-\Tr\bigl(D^2 f(W(h))  \langle h,h^*\rangle_{L^2(0,T)} \bigr),\\
 W(h)&=&\bigl(W(h_1),\cdots, W(h_m) \bigr)  .
\end{eqnarray*}
This operator is also closable and $\D^\infty$ is included in its domain ${\rm Dom}(L) $.\\

\textit{Integration by parts.}\\

\begin{PROP}[Integration by parts: first version]
\label{IBP_V0}
Let $F\in\D^{1,2}$, $u\in  {\rm Dom} (\delta)$, then the following indentity holds:
\begin{eqnarray*}
\E[\langle DF,u\rangle_{L^2([0,T])}]=\E[F\delta(u)],
\end{eqnarray*}
that is the Skorohod integral $\delta $ is the adjoint of the Malliavin derivative $D$.
As a consequence, for $F,G\in {\rm Dom}(L)$ we have
\begin{eqnarray*}
\E[FLG]=\E[F\delta(DG)]=\E[\langle DF,DG\rangle_{L^2([0,T])}]=\E[LFG],
\end{eqnarray*}
i.e. $L$ is self-adjoint.
\end{PROP}
These relations can be easily checked for $F,G \in {\cal S},\ u\in{\cal P}$, and extended to the indicated domains thanks to the closability.

\subsection{Chaos Decomposition}
\label{CAOS_DEC}
$$I_m(f_m):=m! \int_0^T \int_0^{t_1}\cdots \int_0^{t_{m}-1}f_m(t_1,\cdots,t_m) \otimes dW_{t_m} \otimes \cdots \otimes dW_{t_1}.$$
In the above equation $\otimes $ denotes the tensor product and $\left(dW_{t_m}\otimes \cdots\otimes dW_{t_1}\right) \in \left((\R^n)^{\otimes m}\right)^*$.

We now state a theorem that provides a decomposition of real-valued square-integrable random variables in terms of series of multiple integrals.
\begin{lemma}
\label{CHAOS_DEC}
Let $F$ be a real-valued random variable in $L^2(\Omega,\F,\P )$. There exists a sequence $(f_m)_{m\in \N} $ s.t.
\begin{equation}
\label{CD}
F=\sum_{m\in \N} I_m(f_m),
\end{equation}
where for all $m\in \N,\ f_m$ is a symmetric function in $L^2([0,T]^m,(\R^n)^{\otimes m}) $ and
$$\E[F^2]=\sum_{m\ge 0}^{} m!\|f_m\|_{L^2([0,T]^m,(\R^n)^{\otimes m})}^2<+\infty.$$
\end{lemma}
We refer to Theorem  1.1.2 in Nualart \cite{nual:95} for a proof. 
\begin{remark}
We use the term chaos decomposition for the previous expansion because the multiple integral $I_m$ maps  $L^2([0,T]^m,(\R^n)^{\otimes m}) $ onto the Wiener chaos ${\cal H}_m:= \{H_m(W(h)),\ h\in L^2([0,T],\R^n),\ \|h\|_{L^2([0,T],\R^n)}=1 \},$ where $H_m$ stands for the Hermite polynomial of degree $m$ (see again Theorem 1.1.2 in \cite{nual:95}). The orthogonality of the Hermite polynomials yields the orthogonality of the Wiener chaos, i.e. $\E[XY]=0$, for $(X,Y)\in ({\cal H}_n,{\cal H}_m),\ n\neq m$.
\end{remark}

The computation of Malliavin derivatives is quite simple for multiple integrals. Indeed,
\begin{eqnarray*}
D_t (I_m(f_m))&=&mI_{m-1}(f_m(t,.))\in \R^n.
\end{eqnarray*}
As a consequence, for a random variable $F$ having a decomposition as in \eqref{CD}, we have that it belongs to $\D^{1,2}$ if and only if
$\bsum{m\ge 1}^{} m m!\|f_m\|_{L^2([0,T]^m, (\R^n)^{\otimes m})}^2<+\infty$ in which case $D_t F=\sum_{m\ge 1}^{} mI_{m-1}(f_m(t,.)) $ and $\E[\int_0^T|D_t F|^2 dt]=\bsum{m\ge 1}^{} m m!\|f_m\|_{L^2([0,T]^m,(\R^n)^{\otimes m})}^2 $. Iterating the procedure, one gets $F\in \D^{N,2}\iff  \sum_{m=N}^{+\infty}\frac{(m!)^2}{(m-N)!} \|f_m \|_{L^2([0,T]^m,(\R^n)^{\otimes m})}^2<+\infty$ and $D_{t_1,\cdots,t_N}F=\bsum{m=N}^{+\infty} m(m-1)\cdots (m-N+1) I_{m-N}(f_m(t_1,\cdots,t_n,.))\in (\R^n)^{\otimes N} $.

Therefore, when a random variable is smooth in the Malliavin sense, i.e. $\D^\infty$, the Stroock formula, see \cite{stro:87}, provides a representation for the functions $(f_m)_{m\in \N}$ in the chaotic expansion in terms of Malliavin derivatives.
\begin{PROP}[Stroock's formula]
\label{STR_FORM}
Let $F\in \D^\infty$, then the explicit expression of the functions $(f_m)_{m\ge 1} $ in the chaotic expansion \eqref{CD} of $F$ writes:
$$\forall m\in \N,\ f_m(t_1,\cdots,t_m)=\E[D_{t_1,\cdots,t_m}^m F]\in (\R^n)^{\otimes m}.$$
\end{PROP}

For square integrable process, a result analogous to Lemma \ref{CHAOS_DEC} also holds.
\begin{lemma}
Let $(u_t)_{t\in [0,T]}$ be an $\R^n$-valued process in $L^2([0,T]\times \Omega,\F\times {\cal B}([0,T]),dt\otimes d\P)$. There exists a sequence of deterministic  functions $(g_m)_{m\in \N^*} $ s.t.
\begin{equation}
\label{CDP}
u_t=\bsum{m\ge 0}^{}I_m(g_{m+1}(t,.)),
\end{equation}
where the square integrable kernels $g_{m+1}$ are defined on $[0,T]^{m+1} $ with values in $(\R^n)^{\otimes (m+1)} $, are symmetric in the last $m$ variables and s.t. $\sum_{m\ge 0} m!\|g_{m+1}\|_{L^2([0,T]^{m+1}, (\R^n)^{\otimes (m+1)})}^2 <+\infty$.
\end{lemma}
We refer to Lemma 1.3.1 in \cite{nual:95} for a proof when $n=1$.

Also, the Skorohod integral of $u\in {\rm Dom (\delta)} $ is quite direct to compute from its chaotic decomposition \eqref{CDP}. Namely,
$$\delta(u):=  \sum_{m\ge 0}I_{m+1}(\tilde g_m),$$
where $\tilde g_m(t,t_1,\cdots,t_m):=\frac{1}{m+1}\bigl [g_m(t_1,\cdots,t_m,t)+\sum_{i=1}^m g_m(t_1,\cdots,t_{i-1},t,t_{i+1},\cdots,t_m,t_i) \bigr]$
is the symmetrization of $g_m$ in $[0,T]^{m+1}$.

\subsection{Representation of densities through Malliavin calculus}
For $F=(F_1,\cdots, F_N) \in (\D^\infty)^N$, we define the Malliavin covariance matrix $\gamma_F $ by
$$\gamma_F^{i,j}:=\langle DF^i, DF^j\rangle_{L^2(0,T)}, \forall (i,j)\in\leftB 1,N\rightB^2.$$

Let us now introduce the non-degeneracy condition
\begin{trivlist}
\item[\Ass{ND}] We say that the random vector $F=(F_1,\cdots, F_N)$ satisfies the non degeneracy condition if $\gamma_F $ is a.s. invertible and $\det(\gamma_F)^{-1}\in \cap_{q\ge 1}L^q(\Omega) $. In the sequel, we denote the inverse of the Malliavin matrix by
$$\Gamma_F:=\gamma_F^{-1}. $$
\end{trivlist}
This non degeneracy condition guarantees the existence of a smooth density, i.e. $C^\infty $, for the random variable $F$, see e.g. Corollary 2.1.2 in \cite{nual:95} or Theorem 9.3 in \cite{iked:wata:89}.

The following Proposition will be crucial in the derivation of an explicit representation of the density.
\begin{PROP}[Second integration by parts]
\label{IBP}
Let $F=(F_1,\cdots, F_N)\in (\D^\infty)^N$ satisfy the nondegeneracy condition \Ass{ND}. 
Then, for all smooth function $\varphi$ with polynomial growth, $G\in \D^\infty $ and all multi-index $\alpha$,
\begin{eqnarray*}
\E[\partial_{\alpha}\varphi(F)G]\!\!\!\!&=&\!\!\!\!\E[\varphi(F)H_{\alpha}(F,G)],\\
H_i(F,G)\!\!\!\!&=&\!\!\!\!-\bsum{j=1}^{N}\{ G\langle D \Gamma_F^{ij}, DF^j\rangle_{L^2(0,T)}+\Gamma_F^{ij}\langle DG, DF^j\rangle_{L^2(0,T)}-\Gamma_F^{ij}GLF^j \}, \ \forall i\in \leftB 1 ,N\rightB,\\
H_\alpha(F,G)\!\!\!\!&=&\!\!\!\!H_{(\alpha_1,\cdots,\alpha_m)}(F,G)=H_{\alpha_m}(F,H_{(\alpha_1,\cdots,\alpha_{m-1})}(F,G)).
\end{eqnarray*}
Also, for all $q>1$, and all multi-index $\alpha$, there exists $(C,q_0,q_1,q_2,r_1,r_2)$ only depending on $(q,\alpha) $ s.t.
\begin{equation}
\label{MEYER}
\|H_\alpha(F,G)\|_q\le C \| \Gamma_F \|_{q_0} \|G \|_{q_1,r_1}\|F\|_{q_2,r_2}.
\end{equation}
\end{PROP}
For the first part of the proposition
 we refer to Section V-9 of \cite{iked:wata:89}. Concerning equation \eqref{MEYER}, it can be directly derived from the Meyer inequalities on $\|LF\|_q $ and the explicit definition of $H$, see also Proposition 2.4 in Bally and Talay \cite{ball:tala:96:1}.

A crucial consequence of the integration by parts formula is the following representation for the density.
  \begin{corollary}[Expression of the density and 
  upper bound]  \label{REP_DENS}
  Let $F$ $=(F_1,\cdots,F_N)\in (\D^\infty)^N$ satisfy the nondegeneracy condition \Ass{ND}. The random vector $F$ admits a density on $\R^N$.
  Fix $y\in \R^N$. Introduce $\forall (u,v)\in \R^2, \varphi_0^u(v)=\I_{v>u} ,\ \varphi_1^u(v)=\I_{v\le u} $. For all multi-index $\beta=(\beta_1,\cdots,\beta_N)\in \{0,1 \}^N $ the density writes:
 \begin{equation}
 \label{EQ_REP_DENS}
 p_F(y)=\E[\prod_{i=1}^N\varphi_{\beta_i}^{y_i}(F_i) H_{\alpha}(F,1)](-1)^{|\beta|}, \ \alpha=(1,\cdots,N),\ |\beta|:=\sum_{i=1}^N \beta_i.
 \end{equation}
 As a consequence of \eqref{EQ_REP_DENS} and \eqref{MEYER} we get for all multi-index $\beta\in \{0,1 \}^N $:
\begin{equation}
\label{UPPER_BOUND}
\exists C>0,\ p_F(y)\le C \prod_{i=1}^N  \E[\varphi_{\beta_i}^{y_i}(F_i) ]^{\gamma(i)} \|H_{\alpha}(F,1)\|_{2},\quad \gamma(i)=2^{-(i+1)}.
 \end{equation}
 \end{corollary}
  \noindent \textit{Proof.}  Let $B:=\prod_{i=1}^N [a_i,b_i], \forall i\in \leftB 1,N\rightB, \ a_i<b_i$. Denote for all $u\in \R,\ I_0(u):=(-\infty,u) ,\ I_1(u):=[u,\infty)$.
  Set finally, for all multi-index $\beta \in \{0,1\}^N $, $\forall y\in \R^N,\ \Psi_B^\beta(y)=\int_{\prod_{i=1}^{N} I_{\beta_i}(y_i)}^{} \I_B(x)dx $.
  Proposition \ref{IBP} applied with $\alpha=(1,\cdots,N) $ and  $\Psi_B^\beta $ yields
  \begin{eqnarray}
\label{PREAL_1}
\E[\partial_\alpha \Psi_B^\beta(F)]=\E[\Psi_B^\beta(F)H_\alpha(F,1)].
 \end{eqnarray}
 Now, the r.h.s. of equation \eqref{PREAL_1} writes
\begin{eqnarray}
\label{RHS}
\E[\Psi_B^\beta(F)H_\alpha(F,1)]&=&\E[\int_{\prod_{i=1}^{N}I_{\beta_i}(F_i)}^{} \I_{B}(y)  dy H_\alpha(F,1)]=\bint{B}^{}\E[\prod_{i=1}^{N}\I_{y_i\in I_{\beta_i}(F_i)  }H_\alpha(F,1)] dy\nonumber\\
&=& \bint{B}^{}\E[\prod_{i=1}^N\varphi_{\beta_i}^{y_i}(F_i) H_{\alpha}(F,1)] dy.
\end{eqnarray}
The application of Fubini's theorem for the last but one equality is justified thanks to the integrability condition \eqref{MEYER} of Proposition \ref{IBP}.
 On the other hand, the l.h.s. in \eqref{PREAL_1} writes
 \begin{equation}
 \label{LHS}
\E[\partial_\alpha \Psi_B^\beta(F)]=\E[\prod_{i=1}^N \I_{F_i\in [a_i,b_i]}(-1)^{\beta_i}]=(-1)^{|\beta|}\bint{B}^{}p_F(y)dy.
 \end{equation}
Equation \eqref{EQ_REP_DENS} is now a direct consequence of \eqref{PREAL_1}, \eqref{RHS}, \eqref{LHS}.
Equation \eqref{UPPER_BOUND} is then simply derived  applying iteratively the Cauchy-Schwarz inequality
.\hfill $\square$

 \setcounter{equation}{0}
\setcounter{theorem}{0}

\section{Malliavin Calculus to Derive Upper and Diagonal Bounds in our Examples} \label{MALLIAVIN_SECTION}

\subsection{Strategy and usual Brownian controls}
We here concentrate on the particular case of the process \eqref{PROC} (indeed the estimates concerning \eqref{PROC_2} can be derived in a similar way). Since condition \Ass{H} is satisfied, assumption \Ass{ND} is fullfilled. It then follows from Theorem 2.3.2 in \cite{nual:95} that the  process $(X_s)_{s\ge 0}$ admits a smooth density $p(t,x,.) $ at time $t>0$. Our goal is to derive quantitative estimates on this density, emphasizing as well that we have different regimes in function of the starting/final points.

To do that, we condition w.r.t. to the non-degenerate Brownian component for which we explicitly know the density. For all $(t,x,\x)\in \R^{+*}\times (\R^{n+1})^2 $ we have:
\begin{eqnarray*}
&& p(t,x,\x)= p_{X^{1,n}}(t,x_{1,n},\x_{1,n})p_{X^{n+1}}(t,x_{n+1},\x_{n+1}|X_0^{1,n}=x_{1,n},X_t^{1,n}=\x_{1,n}  ),\nonumber \\
&& p_{X^{1,n}}(t,x_{1,n},\x_{1,n})=\frac{1}{(2\pi t)^{n/2}}\exp\left(-\frac{|\x_{1,n}-x_{1,n}|^2}{2t} \right)\nonumber.
\end{eqnarray*}
We then focus on the conditional density which agrees with the one of a smooth functional, in the Malliavin sense, of the Brownian bridge.
Precisely:
\begin{eqnarray}
&&p_{X^{n+1}}(t,x_{n+1},\x_{n+1}|X_0^{1,n}=x_{1,n},X_t^{1,n}=\x_{1,n}  ):=p_{Y_t}(\x_{n+1}-x_{n+1}),\nonumber\\
&&Y_t:=\int_0^t \left |x_{1,n}\frac {t-u}t+\x_{1,n}\frac ut +W_u^{0,t}\right|^kdu,\label{COND}
\end{eqnarray}
where $(W_u^{0,t})_{u\in [0,t]} $ is the standard $n$-dimensional Brownian bridge on the interval $[0,t] $.
The estimation of $p_{Y_t}$ is the core of the probabilistic part of the current work.


We recall, see e.g. \cite{revu:yor:99}, two ways to realize the standard $n$-dimensional Brownian bridge from a standard Brownian motion of $\R^n$. Namely, if $(W_t)_{t\ge 0}$ denotes a standard $n$-dimensional Brownian motion then
\begin{eqnarray}
\label{ID_LAW_BB}
(W_u-\frac ut W_t)_{u\in [0,t]} \overset{({\rm law})}{=}(W_u^{0,t})_{u\in [0, t]}, \label{BB_FINAL_V} \\
\left((t-u)\int_0^u\frac{dW_s}{t-s}\right)_{u\in [0,t]}\overset{({\rm law})}{=}(W_u^{0,t})_{u\in [0, t]} \label{BB_VIA_INT}.
\end{eqnarray}
To recover the framework of Section \ref{OP_M}, in order to deal
 with functionals of the Brownian increments, it is easier to consider the realization of the Brownian bridge given by \eqref{BB_VIA_INT}.
 \begin{remark}
\label{REM_TEMPO_RITROSO}
  The process $(\bar W_u)_{u\in [0,t]}:=(W_{t-u}-W_t)_{u\in [0,t]}$ is a Brownian
motion. Moreover, the processes $(\bar W_u-\frac ut \bar W_t)_{u\in
[0,t]} $ and  $((t-u)\int_0^u\frac{d\bar W_s}{t-s})_{u\in [0,t]} $ are
standard $n$-dimensional Brownian bridges on $[0,t]$ , as well.
  \end{remark}
For the sake of completeness, we recall some well known results
concerning the Brownian motion and Brownian bridge.
 \begin{PROP}
 \label{BROWNIAN_PROP}
 Let $q\ge 1$, and $(W_t)_{t\ge 0} $ be a standard $n$-dimensional Brownian motion. Then, there exists $C:=C(q,n)>0$  s.t.  for all $t\ge 0 $,
  \begin{eqnarray*}
  \E[|W_t|^q]\le C(q,n)t^{q/2},\quad \E[\sup_{s\in[0,t]}|W_s|^q]\le C(q,n)t^{q/2},\\
  \quad\E[\sup_{s\in [\tau ,t]}|W_s^{0,t}|^q]\le C(q,n)(t-\tau)^{q/2},\ 0\le \tau\le t.
  \end{eqnarray*}
  Moreover, there exists $\bar c:=\bar c(n)\ge 1$, s.t. for all $\zeta\ge 0$, and $0\le \tau\le t $,
\begin{eqnarray*}
\P[\sup_{s\in [\tau,t]} |W_s^{0,t}|\ge \zeta	]\le 2\exp\left(-\frac{|\zeta|^2}{\bar c(n)(t-\tau)}\right).
  \end{eqnarray*}
 \end{PROP}
 \noindent \textit{Proof.} The first inequality is a simple consequence of the Brownian scaling. The second one can be derived from convexity inequalities and L\'evy's identity that we now recall (see e.g. Chapter 6 in \cite{revu:yor:99}). Let $(B_t)_{t\ge 0}$ be a standard scalar Brownian motion.
Then:
\begin{equation}
\label{ID_LEVY}
\sup_{u\in [0,s]} B_u\overset{({\rm law})}{=} |B_s|, \forall s>0.
\end{equation}
The third inequality follows from the first two and the representation \eqref{BB_FINAL_V}.
Eventually, the deviation estimates follow from \eqref{ID_LEVY} as well. These deviations estimates can also be seen as special cases of Bernstein's inequality, see e.g. \cite{revu:yor:99} p. 153.
\hfill $\square $

\subsection{Some preliminary estimates on the Malliavin derivative and covariance matrix}
We now give the expressions of the Malliavin derivative and \textit{covariance} matrix of the \textit{scalar} random variable $Y_t$ defined in \eqref{COND} and some associated controls.

\begin{lemma}[Malliavin Derivative and some associated bounds]
\label{LEMMA_MAL}
Let us set\\
$m(u,t,x_{1,n},\x_{1,n}):=x_{1,n}\frac {t-u}t+\x_{1,n}\frac ut $. Rewrite
\begin{eqnarray}
Y_t&=&\int_0^t du |m(u,t,x_{1,n},\x_{1,n})+W_u^{0,t}|^{k}\nonumber\\
&=&\int_0^t du \{|m(u,t,x_{1,n},\x_{1,n})|^2+|W_u^{0,t}|^2+2\langle m(u,t,x_{1,n},\x_{1,n}),W_u^{0,t}\rangle \}^{k/2} \nonumber\\
&=&\bsum{i=0}^{k/2}C_{k/2}^i\int_0^t du |m(u,t,x_{1,n},\x_{1,n})|^{k-2i}\{|W_u^{0,t}|^2+2\langle m(u,t,x_{1,n},\x_{1,n}),W_u^{0,t}\rangle \}^{i}.\nonumber\\\label{DEC_Y}
\end{eqnarray}Considering the realization \eqref{BB_VIA_INT} of the Brownian bridge, the Malliavin derivative
of $Y_t$ (seen as a column vector) and the ``covariance" matrix (that is in our case a scalar) write for all $s\in [0,t] $:
\begin{eqnarray}
D_s Y_t&=&\bsum{i=1}^{k/2}C_{k/2}^i \int_s^t du  |m(u,t,x_{1,n},\x_{1,n})|^{k-2i} i \{|W_u^{0,t}|^2+2\langle m(u,t,x_{1,n},\x_{1,n}),W_u^{0,t}\rangle \}^{i-1}\nonumber\\
&&\times 2 \frac{t-u}{t-s} \bigl(W_u^{0,t}+m(u,t,x_{1,n},\x_{1,n}) \bigr):=\bsum{i=1}^{k/2}M_i(s,t,x_{1,n}, \xx1n),\nonumber \\
\gamma_{Y_t}&=&\int_0^t ds |D_s Y_t|^2. \label{DER_MALL}
\end{eqnarray}
Introduce now
\begin{eqnarray}
M_1(s,t,\X1n,\xx1n)&:=& k \int_s^t du  |m(u,t,x_{1,n},\x_{1,n})|^{k-2}\frac{t-u}{t-s}m(u,t,x_{1,n},\x_{1,n})\nonumber\\
&&+ M_1^R(s,t,\X1n,\xx1n):=(M_1^D+M_1^R)(s,t,\X1n,\xx1n), \label{MDECOMP-1} \\
R(s,t,\X1n,\xx1n)&:=&M_1^R(s,t,\X1n,\xx1n)+\bsum{i=2}^{k/2}M_i(s,t,\X1n,\xx1n),\nonumber\\
 \gamma_{Y_t}&=&\int_0^t ds |(M_1^D+R)(s,t,\X1n,\xx1n)|^2.\label{MDECOMP}
\end{eqnarray}
Set for all $\tau \in [0,t]$,
\begin{eqnarray}
\label{DEFINIZIONI_MR}
{\cal M}_{\tau,t}&:=&\int_\tau^tds |M_1^D(s,t,\x_{1,n},\xi_{1,n})|^2,\quad {\cal M}_t:={\cal M}_{0,t},\nonumber\\
{\cal R }_{\tau,t}&:=&\int_\tau^t ds |R(s,t,\x_{1,n},\xi_{1,n})|^2,\quad {\cal R}_t:={\cal R}_{0,t}.
\end{eqnarray}
There exists $C:=C(k,n) \ge 1$ s.t. for all $\tau \in [0,t]$:
\begin{eqnarray}
\label{equivMT}
C^{-1 }(t-\tau)^3(|\X1n|^{2(k-1)}+|\xx1n|^{2(k-1)}) \le {\cal M}_{\tau, t}\le C (t-\tau)^3(|\X1n|^{2(k-1)}+|\xx1n|^{2(k-1)}).
\end{eqnarray}
Also, for all $q\ge 1$, there exists $C(k,n,q)$ s.t.
\begin{eqnarray}
\E[|{\cal R}_{\tau,t}|^{q}]^{1/q}\!\!&\le&\!\! C(k,n,q) (t-\tau)^{3}( |\X1n|\vee |\xx1n|  )^{2(k-1)}\nonumber\\
&&\times \frac{(t-\tau)}{(|\X1n|\vee |\xx1n|)^2} \biggl ( 1+ \frac{(t-\tau)^{1/2}}{|\X1n|\vee |\xx1n|}  \biggr)^{2(k-2)},\label{CTR_R_MEDIA} \\
\label{CTR_R_LM}
 \forall \kappa\ge 0,\quad \P[{\cal R}_{\tau,t}\ge  \kappa{\cal M}_{\tau,t}] \!\!&\le &\!\! \bar c(n,k)\exp\left(-\kappa^2\frac{ (|\xx1n|\vee |\X1n|)^2}{\bar c(n,k)(t-\tau)} \right),
 \end{eqnarray}
for some constant $\bar c(n,k)\ge 1$.
\end{lemma}

\begin{remark}
\label{REM_REGIMI}
From \eqref{equivMT} and \eqref{CTR_R_MEDIA}, it follows that
\begin{eqnarray*}
\E[{\cal R}_t^q]^{1/q}\le \frac{C(k,n,q)C}{K^2}\left(1+\frac 1K\right)^{2(k-
2)}{\cal M}_t, \
\end{eqnarray*}
 when  $|x_{1,n}|\vee |\x_{1,n}|\ge Kt^{1/2} $. For $K:=K(k,n,q)$ large enough, then the term ${\cal M}_t$ (corresponding to the Malliavin covariance matrix of a Gaussian contribution) dominates the remainder. This intuitively explains the Gaussian regime appearing in \textit{ii)} of Theorem \ref{MTHM}.
\end{remark}

\noindent\textit{Proof.}
Assertion \eqref{DER_MALL} directly follows from the chain rule (see e.g. Proposition 1.2.3 in \cite{nual:95}) and the identity $ D_s W_u^{0,t}=\I_{s\le u}\frac{t-u}{t-s},\ \forall  (u,s)\in [0,t]^2$ deriving from \eqref{BB_VIA_INT}.

Concerning \eqref{equivMT}, we only prove the claim for $\tau=0 $ for notational simplicity. Usual computations involving convexity inequalities yield that there exists $C:=C(k,n)\ge 1$ s.t.
\begin{equation}
\label{UPBDMT}
{\cal M}_t\le C t^3(|\X1n|^{2(k-1)}+|\xx1n|^{2(k-1)}).
\end{equation}
On the other hand to prove that a lower bound at the same ordre also holds for ${\cal M}_t$ one has to be a little more careful.

W.l.o.g. we can assume that $|\xx1n|\ge |\X1n|$.  Indeed, because of the symmetry of the Brownian Bridge and its reversibility in time (see Remark \ref{REM_TEMPO_RITROSO}), if $|\xx1n|< |\X1n|$ we can perform the computations w.r.t. to the Brownian bridge $(\bar W_u^{0,t})_{u\in [0,t]}:=(W_{t-u}^{0,t})_{u\in [0,t]} $ using the sensitivity w.r.t. to the Brownian motion $(\bar W_u)_{u\in [0,t]}:=(W_{t-u}-W_t)_{u\in [0,t]} $.
Note that $|\xx1n|\ge  |\X1n| \Rightarrow |\xx1n|_\infty\ge \frac 1{n^{1/2}} |\X1n|_\infty  $. Let $i_0\in \leftB 1,n\rightB$ be the index s.t. $|\xx1n|_\infty:=|\x_{i_0}| $, then $|\x_{i_0}|\ge \frac1{n^{1/2}}|x_{i_0}| $. Let us now write
\begin{eqnarray*}
{\cal M}_t\ge k^2\int_0^t ds \left( \int_s^t du \left|m(u,t,\X1n,\xx1n)\right|^{k-2}\left(\frac{t-u}{t} x_{i_0} +\frac ut \x_{i_0} \right)\frac{t-u}{t-s}\right)^2.
\end{eqnarray*}
Observe now that for $s\ge \frac{n^{1/2}}{n^{1/2}+1} t$ we have that $\forall u\in [s,t],  \frac{t-u}{t} x_{i_0} +\frac ut \x_{i_0} $ has the sign of $\x_{i_0} $. Hence,
\begin{eqnarray}
\label{CTR_PREL_Mt}
{\cal M}_t\ge k^2\int_{\frac{n^{1/2}}{n^{1/2}+1}t}^tds\left(\int_s^t du \left|\frac{t-u}{t} x_{i_0} +\frac ut \x_{i_0} \right|^{k-1}\frac{t-u}{t-s}\right)^2.
\end{eqnarray}
Now, for $s\ge  
t\left(\frac{n^{1/2}}{n^{1/2}+2^{-1}} \right)
$, we have for all $u\in [s,t] $:
\begin{eqnarray}
 \left|\frac{t-u}{t} x_{i_0} +\frac ut \x_{i_0} \right|^{k-1}&\ge &  \left( \frac{u}{t}\right)^{k-1}\frac{|\x_{i_0}|^{k-1}}{2^{k-2}}-\left( \frac{t-u}{t}\right)^{k-1}|x_{i_0}|^{k-1}\nonumber \\
 &\ge & |\x_{i_0}|^{k-1}\left(\frac{n^{1/2}}{2n^{1/2}+1}\right)^{k-1}. \label{CTR_CONV}
\end{eqnarray}
Equation \eqref{equivMT} thus follows from \eqref{CTR_PREL_Mt}, \eqref{CTR_CONV}  and \eqref{UPBDMT}.

Concerning the remainders we get that there exists $C_3:=C_3(n,k), C_4:=C_4(n,k)$ s.t.:
\begin{eqnarray}
|M_1^R(s,t,\X1n,\xx1n)|^2&\le& C_3 (t-s)^2|\xx1n|^{2(k-1)} \sup_{u\in [s,t]}|W_u^{0,t}|^2|\xx1n|^{-2},\nonumber \\
\forall i\in \leftB 2,k/2\rightB,\ |M_i(s,t,\X1n,\xx1n)|^2&\le& C_4(t-s)^2|\xx1n|^{2(k-1)} \biggl\{\sup_{u\in[s,t]}|W_u^{0,t}|^{2(2i-1)}|\xx1n|^{-2(2i-1)}\nonumber \\
&&+\sup_{u\in[s,t]}|W_u^{0,t}|^{2(i-1)}|\xx1n|^{-2(i-1)}\biggr\}.
\label{CTR_MI}
\end{eqnarray}
From \eqref{MDECOMP} and a convexity inequality, we derive $|R(s,t,\X1n,\xx1n)|^2\le \frac{k}2( |M_1^R(s,t,\X1n,\xx1n)|^2 +\sum_{i=2}^{k/2} |M_i(s,t,\X1n,\xx1n)|^2) $. Thus, from \eqref{CTR_MI}, \eqref{DEFINIZIONI_MR}, we obtain that there exists $C(k,n,q)$ s.t. for all $\tau\in [0,t] $:\begin{eqnarray*}
\E[|{\cal R}_{\tau,t}|^q]^{1/q}&\le& C(k,n,q)(t-\tau)^{3}|\xx1n|^{2(k-1)}\biggl\{\sum_{i=1}^{k/2} \E[|\sup_{u\in[\tau,t]}|W_u^{0,t}|^{2(2i-1)}|\xx1n|^{-2(2i-1)}|^q]^{1/q}\\
&&+\sum_{i=2}^{k/2} \E[|\sup_{u\in[\t,t]}|W_u^{0,t}|^{2(i-1)}|\xx1n|^{-2(i-1)}|^q]^{1/q}\biggr\},
\end{eqnarray*}
which, thanks to Proposition \ref{BROWNIAN_PROP}, gives \eqref{CTR_R_MEDIA}.

On the other hand, from \eqref{equivMT} and the previous convexity inequality for $R$ we get:
\begin{eqnarray*}
&&\P[ \kappa {\cal M}_{\t,t} \le {\cal R}_{\t,t}]\le \P\biggl[ C^{-1}(t-\tau)^3 |\xx1n|^{2(k-1)}\kappa \le \frac k2 \biggl[ \int_{\t}^t ds |M_1^R(s,t,\X1n,\xx1n)|^2   \nonumber\\
&& +\bsum{i=2}^{k/2}\int_{\t}^{t} ds |M_i(s,t,\X1n,\xx1n)|^2\biggr] \biggr]\nonumber\\
&\le& \P\left[ \left(\frac{2}k\right)^2 C^{-1}(t-\tau)^3 |\xx1n|^{2(k-1)} \kappa \le  \int_{\t}^t ds |M_1^R(s,t,\X1n,\xx1n)|^2\right]\\
&&+\bsum{i=2}^{k/2}\ \P\left[\left(\frac{2}k\right)^2 C^{-1}(t-\tau)^3 |\xx1n|^{2(k-1)}\kappa\le \int_{\t}^{t} ds |M_i(s,t,\X1n,\xx1n)|^2 \right] \\
&\overset{\eqref{CTR_MI}}{\le}& \P\left[ \left(\frac{2}k\right)^2 C^{-1}(t-\tau)^3 |\xx1n|^{2(k-1)} \kappa\le \frac {C_3}3(t-\tau)^3|\xx1n|^{2(k-1)}  \sup_{u\in [\t,t]}|W_u^{0,t}|^2|\xx1n|^{-2}\right]\label{THE_PM}\\
&&\!\!\!+\bsum{i=2}^{k/2}\ \P\left[\left(\frac{2}k\right)^2\!\! C^{-1}\kappa \le \! \frac{C_4}3\biggl\{   \biggl( \frac{\sup_{u\in[\tau,t]}|W_u^{0,t}|}{|\xx1n|}\biggr)^{2(2i-1)}\!\!+\biggl(\frac{\sup_{u\in[\t,t]}|W_u^{0,t}|}{|\xx1n|}\biggr)^{2(i-1)} \biggr\}  \right]. 
\end{eqnarray*}
Equation \eqref{CTR_R_LM} then follows from Proposition \ref{BROWNIAN_PROP}.\hfill $\square$


\subsection{Control of the weights}
Now to exploit Corollary \ref{REP_DENS} to give estimates on $p_{Y_t} $ we need to have bounds on the Malliavin weights.  Formula \eqref{EQ_REP_DENS} involves two kinds of terms: the inverse of the Malliavin Matrix and the Ornstein-Uhlenbeck operator. Lemma \ref{LEMMA_MAL} provides tools to analyze the Malliavin matrix. Concerning the Ornstein-Uhlenbeck operator we will rely on the chaos expansion techniques introduced in Section \ref{CAOS_DEC}.

\subsubsection{``Gaussian" regime}
In this section we assume that $|x_{1,n}|\vee |\x_{1,n}|\ge Kt^{1/2} $, for $K:=K(n,d)$ sufficiently large. That is we suppose that the starting or the final point of the non-degenerate component has greater norm than the characteristic time-scale $t^{1/2}$. In this case, we show below that the dominating term in the Malliavin derivative is the one associated to the non-random term $M_1^D$ in \eqref{MDECOMP-1}. This term corresponds to the Malliavin derivative of a Gaussian process. This justifies the terminology ``Gaussian" regime. 

In order to give precise asymptotics on the density of $Y_t$, the crucial step consists in controlling the norm of $\Gamma_{Y_t}:=\gamma_{Y_t}^{-1} $ in $L^q(\Omega),\ q\in [ 1,+\infty)$ spaces.
\begin{lemma}[Estimates on the Malliavin covariance]
\label{LEMMA_COV}
\vspace*{2pt}Assume that $|\X1n|\vee |\xx1n|\ge  K t^{1/2} $. Then, for all $q\in  [1,+\infty)$ there exists $C_{q,\ref{LEMMA_COV}}:=C_{q,\ref{LEMMA_COV}}(n,k,K)\ge 1$ s.t. 
$$ \frac{C_{q,\ref{LEMMA_COV}}^{-1}}{(|\X1n|^{2(k-1)}+|\xx1n|^{2(k-1)}  ) t^{3}}\le  \|\Gamma_{Y_t}\|_q\le \frac{C_{q,\ref{LEMMA_COV}}}{(|\X1n|^{2(k-1)}+|\xx1n|^{2(k-1)}  ) t^{3}}.$$\end{lemma}

\medskip

\noindent \textit{Proof.} As in Lemma \ref{LEMMA_MAL}, we assume, without loss of generality, that $|\xx1n|\ge |\X1n|$.
To give the $L^q$ estimates of the Malliavin derivative we recall the definition of ${\cal M}_t$ given in \eqref{DEFINIZIONI_MR}, and we use the following partition:
\begin{eqnarray}
\E[|\Gamma_{Y_t}|^q]= \bsum{m \in \N}^{}\E[|\Gamma_{Y_t}|^q\I_{  \Gamma_{Y_t}\in [ \frac{4m}{{\cal M}_t},\frac{4(m+1)}{{\cal M}_t}] }]\le \nonumber \\
\left(\frac{4}{{\cal M}_t}\right)^q+\bsum{m\ge 1}^{} \left(\frac{4(m+1)}{{\cal M}_t}\right)^q \P[\gamma_{Y_t}\le \frac{{\cal M}_t}{4m}].
\label{DECOMP_DIAD}
\end{eqnarray}
Equation \eqref{equivMT} in Lemma \ref{LEMMA_MAL} provides us with an useful bound for ${\cal M}_t$. We next give estimates of $\P\left[\gamma_{Y_t}\le \frac{{\cal M}_t}{4m}\right ],\ m\ge 1 $ in the spirit of Bally \cite{ball:90}.

Introduce $t_m:=\inf\{v\in [0,t]: {\cal M}_{v,t}
\le{\cal M}_t/m\} $.
We first show that there exists $m_0\in \N$ and $\bar C:=\bar C(n,k) $ such that $t_m\ge  t(1-\bar Cm^{-1/3})$ for all $m\ge m_0$.
From \eqref{equivMT} we obtain
\begin{eqnarray*}
t_m&\ge& \inf\{v\in [0,t]: {\cal M}_{v,t}
\le Ct^3(|\X1n|^{2(k-1)}+|\xx1n|^{2(k-1)})/m \}\\
&\ge&\inf\{v\in [0,t]: {\cal M}_{v,t}
\le 2C t^3|\xx1n|^{2(k-1)}/m \} =: \bar t_m,
\end{eqnarray*}
recalling we have assumed $|\xx1n|\ge |\X1n| $ for the last inequality.
Equations \eqref{CTR_PREL_Mt} and \eqref{CTR_CONV} also yield that there exists $C_2:=C_2(n,k)$ s.t. for all $v\ge \frac{n^{1/2}}{n^{1/2}+2^{-1}}t$,
$$
{\cal M}_{v,t}
\ge C_2(t-v)^3|\xx1n|^{2(k-1)} .
$$
Note that $\bar t_m \to t$ as $m \to + \infty$, then there exists $\bar m$ such that $\bar t_m \ge \frac{n^{1/2}}{n^{1/2}+2^{-1}}t$ for every $m \ge \bar m$, and the above inequality holds for every $v \in [\bar t_m, t]$. Set $\bar C:=(2C/C_2)^{1/3}$, and $m_0 = \lfloor \bar C^3 \rfloor \vee \bar m$. For every $m\ge m_0$ we have that:
\begin{eqnarray}
\P\left[\gamma_{Y_t}\le \frac{{\cal M}_t}{4m}\right ]& \le& \P\left[\int_{t_m}^t ds |(M_1^D+R)(s,t,\X1n,\xx1n)|^2\le \frac{{\cal M}_t}{4m}\right]\nonumber\\
&\le& \P\left[\frac12 \int_{t_m}^t ds |M_1^D(s,t,\X1n,\xx1n)|^2-\int_{t_m}^{t} ds|R(s,t,\X1n,\xx1n)|^2 \le \frac{{\cal M}_t}{4m}\right]\nonumber\\&\le & \P\left[\frac{{\cal M}_{t_m,t}}{4}\le  {\cal R}_{t_m,t}  \right]\le  \bar c(n,k) \exp\left(-\frac{|\xx1n|^2m^{1/3}}{16\bar C \bar c(n,k) t} \right),
\label{INTER_GAMMA}
\end{eqnarray}
using \eqref{CTR_R_LM} for the last inequality.
Plugging this control into \eqref{DECOMP_DIAD},  using once again \eqref{equivMT} 
we derive that there exists $C_3:=C_3(n,k),(C_4,C_5):=(C_4,C_5)(n,k,q) $ s.t.:
\begin{eqnarray*}
\E[|\Gamma_{Y_t}|^q]\le \left(\frac{4Cm_0^2}{t^3|\xx1n|^{2(k-1)}} \right)^q+C_3\bsum{m\ge m_0}^{}\left(\frac{4C(m+1)}{t^3|\xx1n|^{2(k-1)}} \right)^q\exp\left(-C_3^{-1} \frac{|\xx1n|^2m^{1/3}}{t}\right)\\
\le  \left(\frac{4Cm_0^2}{t^3|\xx1n|^{2(k-1)}} \right)^q+\left( \frac{C_3(8C)^q}{|\xx1n|^{(2(k-1)+6)q}}\right)\bsum{m\ge m_0}^{}\left(\frac{m^{1/3}|\xx1n|^2}{t} \right)^{3q} \exp\left( -C_3^{-1} \frac{|\xx1n|^2m^{1/3}}{t}\right)\\
\le \left(\frac{4Cm_0^2}{t^3|\xx1n|^{2(k-1)}} \right)^q+\frac{C_4}{|\xx1n|^{(2k+4)q}}\bsum{m\ge m_0}^{} \exp\left( -C_4^{-1} \frac{|\xx1n|^2m^{1/3}}{t}\right)\\
\le C_5\left[\frac{1}{t^{3q}|\xx1n|^{2q(k-1)}}+\frac{1}{|\xx1n|^{(2k+4)q}}\frac{t^3}{|\xx1n|^6}\right]\le\frac{C_5}{t^{3q}|\xx1n|^{2q(k-1)}}\left[1+\frac{t^{3(q+1)}}{|\xx1n|^{6(q+1)}} \right],
\end{eqnarray*}
which for $|\xx1n|\ge Kt^{1/2} $ gives the upper bound of the lemma.
\\

 Let us now turn to the lower bound for $\|\Gamma_{Y_t}\|_{L^p(\P)} $. 
 Write:
 \begin{eqnarray*}
 \E[\Gamma_{Y_t}^q]&\ge &\E[\Gamma_{Y_t}^q\I_{\gamma_{Y_t}\le 3{\cal M}_t}]\ge \frac{1}{(3{\cal M}_t)^q}\P[\gamma_{Y_t}\le 3{\cal M}_t]
 \ge  \frac{1}{(3{\cal M}_t)^q}(1-\P[\gamma_{Y_t}> 3{\cal M}_t]).
 \end{eqnarray*}
 From equations \eqref{DER_MALL}-\eqref{MDECOMP} one has $\P[\gamma_{Y_t}> 3{\cal M}_t]\le \P[2{\cal M}_t+2{\cal R}_t>3{\cal M}_t]=\P[{\cal R}_t>\frac12 {\cal M}_t] $. Now, from Lemma \ref{LEMMA_MAL} equation \eqref{CTR_R_LM}, one gets $\P[\gamma_{Y_t}>3{\cal M}_t]\le \bar c(n,k)\exp\left(-
 \frac{|\xx1n|^2}{4 \bar c(n,k)t}\right) $. Therefore, for $|\xx1n|\ge K t^{1/2}$ and $K$ large enough, we get
 $\E[\Gamma_{Y_t}^q]\ge \frac{1}{2(3{\cal M}_t)^q}$, which thanks to \eqref{equivMT} completes the proof.\hfill $\square $\\

\textbf{Controls of the weight for the integration by parts.} \\
\\
From Proposition \ref{IBP} and Corollary \ref{REP_DENS}, we derive
\begin{eqnarray}
p_{Y_t}(\x_{n+1}-x_{n+1})&=&\E[H_t\I_{Y_t>\x_{n+1}-x_{n+1}}],\nonumber \\
H_t&=&-\langle D\Gamma_{Y_t}, DY_t \rangle_{L^2(0,t)}+\Gamma_{Y_t}LY_t=\gamma_{Y_t}^{-2}\langle D\gamma_{Y_t}, DY_t\rangle_{L^2(0,t)}+\Gamma_{Y_t}LY_t\nonumber\\
&:=&H_t^1+H_t^2,
 \label{SCRIT_IPP}
\end{eqnarray}
using the chain rule 
 for the last but one identity.

We have the following $L^q(\P),\ q\ge 1 $, bounds for the random variable $H_t$.
\begin{PROP}[Estimates for the Malliavin weight]
\label{CTR_W}
\vspace*{2pt}Assume that $|\X1n|\vee |\xx1n|\ge  K t^{1/2} $ for $K$ large enough. Then, for all $q\in  [1,+\infty)$ there exists $C_{q,\ref{CTR_W}}:=C_{q,\ref{CTR_W}}(n,k,K)\ge 1$ s.t.
$$ 
\|H_t\|_q\le \frac{C_{q,\ref{CTR_W}}}{(|\X1n|^{(k-1)}+|\xx1n|^{(k-1)}  ) t^{3/2}}.$$
\end{PROP}
\textit{Proof: Control of $H_t^1 $.} From \eqref{SCRIT_IPP} we get for all given $q\ge 1$,
\begin{eqnarray}
\label{CTR_H1_PROV}
\|H_t^1\|_q:=\E[\gamma_{Y_t}^{-2q} |\langle D \gamma_{Y_t},DY_t\rangle_{L^2(0,t)} |^q]^{1/q}\le \E[\gamma_{Y_t}^{-4q}]^{1/2q}\E[|\langle D\gamma_{Y_t},DY_t\rangle_{L^2(0,t)}|^{2q}]^{1/2q}\nonumber\\
\le \frac{C_{q,\ref{LEMMA_COV}}}{t^6(|\xx1n|^{2(k-1)}+|\X1n|^{2(k-1)})^2} \E[|\langle D\gamma_{Y_t},DY_t\rangle_{L^2(0,t)}|^{2q}]^{1/2q}\nonumber\\
\le \frac{C_{q,\ref{LEMMA_COV}}}{t^6(|\xx1n|^{2(k-1)}+|\X1n|^{2(k-1)})^2}\E[| D\gamma_{Y_t}|_{L^2(0,t)}^{4q}]^{1/4q}\E[|DY_t|_{L^2(0,t)}^{4q}]^{1/4q},
\end{eqnarray}
 using Lemma \ref{LEMMA_COV} for the last but one inequality. Now, from equations \eqref{DER_MALL}, \eqref{MDECOMP}, using the notations of Lemma \ref{LEMMA_COV},
 \begin{eqnarray*}
 \E[|DY_t|_{L^2(0,t)}^{4q}]^{1/4q}&=&\E[(\int_0^t ds |D_s Y_t|^2 )^{2q}]^{1/4q}:=\E[\gamma_{Y_t}^{2q}]^{1/4q}\\
 &\le & \left(2^{4q-1}\left\{ {\cal M}_t^{2q}+\E[{\cal R}_t^{2q}]\right\}\right)^{1/4q}\le 2^{1-1/4q}\left\{ {\cal M}_t^{1/2}+\E[{\cal R}_t^{2q}]^{1/4q}\right\}.
 \end{eqnarray*}
 On the one hand equation \eqref{equivMT} in Lemma \ref{LEMMA_MAL} readily gives ${\cal M}_t^{1/2}\le Ct^{3/2}(|\X1n|^{k-1}+|\xx1n|^{k-1}) $.
 On the other hand, equation \eqref{CTR_R_LM} of the same Lemma yields
\begin{eqnarray*}
\E[|{\cal R}_t|^{2q}]^{1/4q} &\le & C(k,q) \left( t^{3/2}\big( |\xx1n|^{k-1} \vee |\X1n|^{k-1} \big) K^{-1} k/2\right).
\end{eqnarray*}
Hence, there exists $C_{1}:=C_{1}(n,k,q,K)$ s.t.
\begin{eqnarray}
\label{CTR_DYT_LP}
 \E[|DY_t|_{L^2(0,t)}^{4q}]^{1/4q}=\E[\gamma_{Y_t}^{2q}]^{1/4q}&\le& C_{1} t^{3/2}(|\xx1n|^{k-1}+|\X1n|^{k-1}).
\end{eqnarray}
In order to get a bound for $\|H_t^1\|_q $, it remains to control $\E[| D\gamma_{Y_t}|_{L^2(0,t)}^{4q}]^{1/4q}$.
Equation \eqref{DER_MALL} and the chain rule yield that for all $u_2\in [0,t] $, $D_{u_2}\gamma_{Y_t}=2\int_0^t du_1  D_{u_2}D_{u_1} Y_t \times D_{u_1} Y_t$. We get
\begin{eqnarray}
\E[| D\gamma_{Y_t}|_{L^2(0,t)}^{4q}]^{1/4q}&\le& 4 \E[\gamma_{Y_t}^{4q}]^{1/8q}\E[(\int_0^t du_1\int_0^t du_2 |D_{u_2,u_1}Y_t|^2)^{4q}]^{1/8q}\nonumber \\
&\le &C_{2} t^{3/2}(|\xx1n|^{k-1}+|\X1n|^{k-1})\E[|D^2 Y_t|_{L^2((0,t)^2)}^{8q}]^{1/8q},\label{CTR_DER_GYT}
\end{eqnarray}
$C_2:=C_2(n,k,q,K) $ using \eqref{CTR_DYT_LP} for the last inequality.

With the notations of equations \eqref{DER_MALL}, \eqref{MDECOMP} we set for all $u_1\in [0,t] $,
\begin{eqnarray*}
D_{u_1}Y_t&:=&
\bsum{i=1}^{k/2} M_i(u_1,t,\X1n,\xx1n):=\bsum{i=1}^{k/2} \bar M_i(u_1,t),
\end{eqnarray*}
for simplicity. 

Observe now that for all $i\in\leftB 2,k/2\rightB$, $u_2\in [0,t] $,
\begin{eqnarray*}
D_{u_2}\bar M_i(u_1,t)=C_{k/2}^i\int_{u_1\vee u_2}^t dv|m(v,t,\X1n,\xx1n)|^{k-2i}\left\{|W_v^{0,t}|^2+2\langle m(v,t,\X1n,\xx1n), W_v^{0,t}\rangle \right\}^{i-2}\\
\times 2i \frac{(t-v)^2}{(t-u_1)(t-u_2)}\left\{ 2(i-1) (W_v^{0,t}+m(v,t,\X1n,\xx1n))\otimes(W_v^{0,t}
+m(v,t,\X1n,\xx1n))\right.\\
 \left.+\left\{|W_v^{0,t}|^2+2\langle m(v,t,\X1n,\xx1n), W_v^{0,t}\rangle \right\}I_n\right\},\\
D_{u_2}\bar M_1(u_1,t)=k\int_{u_1\vee u_2}^t dv|m(v,t,\X1n,\xx1n)|^{k-2} \frac{(t-v)^2}{(t-u_1)(t-u_2)} I_n.
\end{eqnarray*}
From the above equations, assuming once again w.l.o.g. $|\xx1n|\ge |\X1n| $, the arguments used in Lemma \ref{LEMMA_MAL} yield:
\begin{eqnarray*}
\E[|D^2 Y_t|_{L^2((0,t)^2)}^{8q}]^{1/8q}\le C\E[|\int_{[0,t]^2} du_1du_2(t-u_1\vee u_2)^2|\xx1n|^{2(k-2)}\\
\times(1+\sum_{i=2}^{k/2}|\xx1n|^{4(1-i)} \sup_{u\in [0,t]}|W_u^{0,t}|^{4(i-1)}) |^{4q}]^{1/8q}\\
\le  Ct^2|\xx1n|^{k-2}(1+\sum_{i=2}^{k/2}|\xx1n|^{2(1-i)}\E[\sup_{u\in [0,t]}|W_u^{0,t}|^{16q(i-1)}]^{1/8q} )\\
\overset{{\rm Prop.}\ \ref{BROWNIAN_PROP}}{\le} Ct^2|\xx1n|^{k-2}(1+\sum_{i=2}^{k/2}\left(\frac{t^{1/2}}{|\xx1n|}\right)^{2(i-1)}),
\end{eqnarray*}
where $C:=C(n,k,q)$ may change from line to line. Recalling that $|\xx1n|\vee |\X1n|\ge Kt^{1/2} $ we obtain
\begin{eqnarray*}
\E[|D^2 Y_t|_{L^2((0,t)^2)}^{8q}]^{1/8q}\le Ct^2|\xx1n|^{k-2}, \ C:=C(n,k,q,K).
\end{eqnarray*}
Plugging the above equation into \eqref{CTR_DER_GYT} we derive that
\begin{eqnarray*}
\E[|D\gamma_{Y_t}^{4q}|_{L^2(0,t)}]^{1/4q}\le Ct^{7/2}|\xx1n|^{2k-3},
\end{eqnarray*}
which together with \eqref{CTR_DYT_LP} and \eqref{CTR_H1_PROV}, eventually yields
\begin{eqnarray}
\label{CTR_HT1}
\|H_t^1\|_q \le \frac{\bar C_1}{t|\xx1n|^k}\le \frac{CK^{-1}}{t^{3/2}|\xx1n|^{k-1}}, \bar C_1:=\bar C_1(n,k,q,K).
\end{eqnarray}

\medskip

\noindent \textit{Control of $H_t^2$}. From \eqref{SCRIT_IPP} and Lemma \ref{LEMMA_COV}, for all $q\ge 1$, we get
 \begin{eqnarray}
\label{CTR_H2_PROV}
 \|H_t^2\|_q\le \E[|\Gamma_{Y_t}|^{2q}]^{1/2q}\E[|LY_t|^{2q}]^{1/2q}\le \frac{C_q}{t^{3}(|\X1n|^{2(k-1)}+|\xx1n|^{2(k-1)})}\E[|LY_t|^{2q}]^{1/2q}. \end{eqnarray}
Now, since $LY_t=\delta(DY_t)$, the idea is to provide a chaotic representation of $DY_t$.
To do that, we use Proposition \ref{STR_FORM} (Stroock's formula see \cite{stro:87}).
For a given $u_1\in [0,t]$, recalling $D_{u_1}Y_t:=\bsum{i=1}^{k/2}\bar M_i(u_1,t) $ where $ \bar M_i(u_1,t)\in \R^n$ is a random contribution involving Wiener chaos up to order $2i-1 $, one has:
\begin{eqnarray*}
\bar M_i(u_1,t)&=&\E[\bar M_i(u_1,t)]+\bsum{l=1}^{2i-1}I_l(g_l^i(.,u_1,t)),\\
I_l(g_l^i(.,u_1,t))&:=&\int_0^t\int_{0}^{v_1}\cdots  \int_{0}^{v_{l-1}} g_l^i(v_1,\cdots,v_l,u_1,t)\otimes dW_{v_l}\otimes \cdots \otimes dW_{v_1},\\
g_l^i(v_1,\cdots,v_l,u_1,t)&:=&\E[D_{v_l,\cdots,v_1}\bar M_i(u_1,t)] \in (\R^n)^{\otimes (l+1)},\ (dW_{v_l}\otimes \cdots \otimes dW_{v_1})\in ((\R^n)^{\otimes l})^*.
 \end{eqnarray*}
Hence, $D_{u_1}Y_t:=g_0(u_1,t)+\bsum{l=1}^{k-1}I_l(g_l(.,u_1,t))$, where $g_0(u_1,t):=\bsum{i=1}^{k/2}\E[\bar M_i(u_1,t)] $ and for all $l\in \leftB1 ,k-1\rightB,\ \ g_l(v_1,\cdots,v_l,u_1,t):=\bsum{i=\lfloor l/2 \rfloor+1}^{k/2} g_l^i(v_1,\cdots,v_l,u_1,t) $, so that
\begin{eqnarray}
\label{INT_SKO_CAOS}
LY_t=\int_0^t g_0(u_1,t)\otimes dW_{u_1}+\bsum{l=2}^{k}I_{l}(g_{l-1}(.,t)):=\bsum{l=1}^{k}I_l(g_{l-1}(.,t)).
\end{eqnarray}
Similarly to the proof performed to control $\E[|{\cal R}_t|^{2q}]^{1/4q}$ of \eqref{CTR_R_LM} in Lemma \ref{LEMMA_MAL}, 
we obtain that there exists $C:=C(n,k)$ s.t. for all $l\in \leftB 0,k-1\rightB $ and for all $(v_1,\cdots,v_l ,u_1)\in [0,t]^{l+1}$:
\begin{equation}
\label{CTR_GL}
|g_l(v_1,\cdots,v_l,u_1,t)|\le Ct(|\xx1n|^{k-(l+1)}+|\X1n|^{k-(l+1)}).
\end{equation}
Therefore,
\begin{eqnarray*}
\E[|LY_t|^{2q}]^{1/2q}&\le &C \bsum{l=1}^k t^{1+l/2}(|\xx1n|^{k-l}+|\X1n|^{k-l} )\\
&\le & Ct^{3/2}(|\xx1n|^{k-1}+|\X1n|^{k-1})\left\{\sum_{l=1}^{k}t^{(l-1)/2}(|\xx1n|^{1-l}+|\X1n|^{1-l} )\right\},
\end{eqnarray*}
where $C:=C(n,k,q)$ may change from line to line. Recalling that $ |\xx1n|\vee |\X1n| \ge Kt^{1/2}$, we derive from \eqref{CTR_H2_PROV} that there exists $\bar C_2:=\bar C_2(n,k,q,K)$ s.t.
\begin{eqnarray*}
 \|H_t^2\|_q\le  \frac{\bar C_2}{t^{3/2}(|\xx1n|^{k-1}+|\X1n|^{k-1})},
\end{eqnarray*}
which together with \eqref{CTR_HT1} and \eqref{SCRIT_IPP} completes the proof.\hfill $ \square$

\subsubsection{Non Gaussian regime}

We now consider the case $|\X1n|\vee |\xx1n|\le Kt^{1/2} $, which corresponds to a diagonal regime of the non-degenerate component w.r.t. the characteristic time scale. It turns out that the characteristic time-scale of the density $p_{Y_t}(\xi_{n+1}-x_{n+1})$ is $t^{1+k/2} $. Indeed, we have the following result.
\begin{PROP}[Estimates for the Malliavin weight in Non Gaussian regime ] \label{EST_NG}
Let $K>0$ be given and assume that $ |\X1n|\vee |\xx1n|\le Kt^{1/2}$. For every $q\ge 1$ there exists $C_{q,\ref{EST_NG}}:=C_{q,\ref{EST_NG}}(n,k,K)$ s.t.
$$\|H_t\|_q\le \frac{C_{q,\ref{EST_NG}}}{t^{1+k/2}}.$$
\end{PROP}
\noindent \textit{Proof.}  For $t>0$ write:
\begin{eqnarray*}
Y_t=\int_0^t\left |\X1n \frac{t-u}{t}+\xx1n \frac ut+W_u^{0,t}\right|^k du &= &t^{1+k/2}\int_0^1 \left| \frac{\X1n}{t^{1/2}}(1-u)+ \frac{\xx1n}{t^{1/2}}u+ \frac{W_{ut}^{0,t}}{t^{1/2}}\right|^k du\\
&=: &t^{1+k/2}\bar Y_1^t.
\end{eqnarray*}
Thus: 
\begin{eqnarray*}
p_{Y_t}(\x_{n+1}-x_{n+1})&:=&-\partial_{\x_{n+1}}\P[Y_t>\x_{n+1}-x_{n+1}]=-\partial_{\x_{n+1}}\P[\bar Y_1^t>\frac{\x_{n+1}-x_{n+1}}{t^{1+k/2}}]\\
&=&\frac1{t^{1+k/2}}p_{\bar Y_1^t}(\frac{\x_{n+1}-x_{n+1}}{t^{1+k/2}}).\end{eqnarray*}
From Corollary \ref{REP_DENS} (Malliavin representation of the densities), we obtain:
\begin{eqnarray*}
p_{Y_t}(\x_{n+1}-x_{n+1})&=&\E[H(Y_t,1)\I_{Y_t>\x_{n+1}-x_{n+1}}]=\frac{1}{t^{1+k/2}}\E[ H(\bar Y_1^t,1)\I_{\bar Y_1^t>\frac{\x_{n+1}-x_{n+1}}{t^{1+k/2}}}]\\
&=&\frac{1}{t^{1+k/2}}\E[H(\bar Y_1^t,1)\I_{Y_t>\x_{n+1}-x_{n+1}}],
\end{eqnarray*}
so that $H_t:= H(Y_t,1)=t^{-(1+k/2)}H(\bar Y_1^t,1):=t^{-(1+k/2)}H_1^{\bar Y_1^t}$. Hence, for all $q\ge 1$,
\begin{equation}
\label{SCAL_MM}
\|H_t\|_q \le \frac1{t^{1+k/2}}\|H_1^{\bar Y_1^t}\|_q.
\end{equation}
Now, as a consequence of the Brownian scaling we get $(\frac{W_{ut}^{0,t}}{t^{1/2}})_{u\in [0,1]}\overset{({\rm law})}{=} (W_u^{0,1})_{u\in [0,1]} $ so that $
\bar Y_1^t \overset{(\rm{law})}{=} t^{1+k/2}\int_0^1 \left| \frac{\X1n}{t^{1/2}}(1-u)+ \frac{\xx1n}{t^{1/2}}u+ W_u^{0,1}\right|^k du$. Recalling that $|\frac{\X1n}{t^{1/2}}|\vee |\frac{\xx1n}{t^{1/2}}|\le K $ we derive that the usual techniques used to prove the non degeneracy of the Malliavin covariance matrix under H\"ormander's condition (see e.g. Norris \cite{norr:86} or Nualart \cite{nual:95}) yield that there exists $C_q:=C_q(n,k,K)\in \R^{+*}$ s.t. $\|H_1^{\bar Y_1^t}\|_q\le C_q$ which from \eqref{SCAL_MM} concludes the proof. The crucial tool here is the global scaling.\hfill $\square $

\subsection{Deviation estimates}

\subsubsection{Off-diagonal bounds}

From the Malliavin representation of the density given by \eqref{SCRIT_IPP}, to derive off-diagonal bounds on the density, it remains to give estimates on $\P[Y_t>\x_{n+1}-x_{n+1}]$.
\begin{lemma}[Off-diagonal bounds]
\label{OFF_DIAG} Let $U_t^k(x,\x):=\x_{n+1}-x_{n+1}-\frac{2^{k-1}}{k+1}(|\X1n|^k+|\xx1n|^k) t$, and assume that $U_t^k(x,\x)>0$.
Then, there exists $C_{\ref{OFF_DIAG}}:=C_{\ref{OFF_DIAG}}(n,k) $ s.t.
\begin{trivlist}
\item[\textit{(i)}] If $|x_{1,n}|\vee |\x_{1,n}|\ge Kt^{1/2}$ for a given $K>0$,
\begin{eqnarray*}
\P[Y_t>\x_{n+1}-x_{n+1}]\le C_{\ref{OFF_DIAG}}\biggl\{
 \exp\left(-C_{\ref{OFF_DIAG}}^{-1}\frac{U_t^k(x,\x)^2}{(|\X1n|^{k-1}+|\xx1n|^{k-1})^2  t^3}\right) \\
+ \exp\left(-C_{\ref{OFF_DIAG}}^{-1}\frac{|\X1n|^2+|\xx1n|^2}{t} \right) \bsum{i=1}^{k/2} \exp\left(-C_{\ref{OFF_DIAG}}^{-1}\frac{U_t^k(x,\xi)^{1/i}}{\{|\X1n|^{k-2i}+|\xx1n|^{k-2i} \}^{1/i}t^{1+1/i}} \right)\biggr\}.
\end{eqnarray*}
\item [\textit{(ii)}] If $|x_{1,n}|\vee |\x_{1,n}|\le Kt^{1/2}$ for the same previous $K$,
\begin{eqnarray*}
\P[Y_t>\x_{n+1}-x_{n+1}]\le C_{\ref{OFF_DIAG}}\biggl\{
 \exp\left(-C_{\ref{OFF_DIAG}}^{-1}\frac{U_t^k(x,\x)^2}{  t^{k+2}}\right) 
+  \bsum{i=1}^{k/2} \exp\left(-C_{\ref{OFF_DIAG}}^{-1}\frac{U_t^k(x,\xi)^{1/i}}{ t^{k/(2i)+1/i}} \right)\biggr\}.
\end{eqnarray*}

\end{trivlist}
\end{lemma}


\noindent \textit{Proof.}  We only prove point \textit{(i)}, the second point can be derived in a similar way. According with \eqref{DEC_Y}, we first decompose $Y_t$ as
\begin{equation*}
    Y_t =\int_0^t |m(u,t,\X1n,\xx1n)|^k du + M_t^k(\X1n,\xx1n) + R_t^k(\X1n,\xx1n)
\end{equation*}
where
\begin{eqnarray*}
& M_t^k(\X1n,\xx1n) := k\int _0^t |m(u,t,\X1n,\xx1n)|^{k-2}\langle m(u,t,\X1n,\xx1n),W_u^{0,t}\rangle du, \\
& R_t^k(\X1n,\xx1n):=\frac k2\int _0^t |m(u,t,\X1n,\xx1n)|^{k-2}|W_u^{0,t}|^2du \\
& +\bsum{i=2}^{k/2}C_{k/2}^i\int_0^t |m(u,t,\X1n,\xx1n)|^{k-2i}(2\langle m(u,t,\X1n,\xx1n),W_u^{0,t}\rangle +|W_u^{0,t}|^2)^idu,
\end{eqnarray*}
then we have
\begin{eqnarray*}
&\P[Y_t>\x_{n+1}-x_{n+1}]= \P[ M_t^k(\X1n,\xx1n) + R_t^k(\X1n,\xx1n)> \\
&\x_{n+1}-x_{n+1}-\int_0^t |m(u,t,\X1n,\xx1n)|^k du].
\end{eqnarray*}
Note that all the terms in $R_t^k(\X1n,\xx1n) $ have characteristic time scales that are in small time negligible with respect to the one of the Gaussian contribution $M_t^k(\X1n,\xx1n)$. Moreover
\begin{equation*}
M_t^k(\X1n,\xx1n) \le 2^{k-2}(|\X1n|^{k-1}+|\xx1n|^{k-1} )t \sup_{u\in [0,t]}|W_u^{0,t}| =: \tilde M_t^k(\X1n,\xx1n).
\end{equation*}
Since by assumption $U_t^k(x,\x)<\x_{n+1}-x_{n+1} -\int_0^t |m(u,t,\X1n,\xx1n)|^k du$,
one gets:
\begin{eqnarray}
\P[Y_t>\x_{n+1}-x_{n+1}]
&\le& \P[(\tilde M_t^k +R_t^k)(\X1n,\xx1n)>U_t^k(x,\x)]\le \nonumber \\
&& \P[2 \tilde M_t^k(\X1n,\xx1n) >U_t^k(x,\x)]\nonumber \\
&+&\P[ (\tilde M_t^k+R_t^k)(\X1n,\xx1n)>U_t^k(x,\x)]^{1/2} \nonumber \\
&&
\times
\P[R_t^k(\X1n,\xx1n)\ge \tilde M_t^k(\X1n,\xx1n)]^{1/2}.\label{CTR_Y}
\end{eqnarray}
Standard computations, similar to the ones performed to prove the deviation estimate in \eqref{CTR_R_LM} in Lemma \ref{LEMMA_MAL},  give that there exist $C_1:=C_1(k), \ C_2:=C_2(n,k)\ge 1$ s.t.
\begin{eqnarray}
\P[R_t^k(\X1n,\xx1n)\ge \tilde M_t^k(\X1n,\xx1n)]&\le& (k-1) \P[\sup_{u\in [0,t]}|W_u^{0,t}|\ge C_1\{|\X1n|+|\xx1n| \}]\nonumber\\
&\le &  C_2\exp\left(-C_2^{-1}\frac{|\X1n|^2+|\xx1n|^2}{t} \right),\nonumber\\
\P[\tilde M_t^k(\X1n,\xx1n) >U_t^k(x,\x)/2] &\le & C_2 \exp\left( -C_2^{-1}\frac{U_t^k(x,\x)^2}{(|\X1n|^{k-1}+|\xx1n|^{k-1})t^3}\right).\label{CTR_EXP_AUX}
\end{eqnarray}
On the other hand, we have:
\begin{eqnarray}
\P[(\tilde M_t^k+R_t^k)(\X1n,\xx1n)\ge U_t^k(x,\x)]\nonumber \\
\le \P[\tilde M_t^k(\X1n,\xx1n)\ge\frac12 U_t^k(x,\x) ]
+\P[R_t^k(\X1n,\xx1n)\ge \frac 12 U_t^k(x,\x)].\label{DECOMP_PROBA_DEV}
\end{eqnarray}
Now,
\begin{eqnarray}
\P[R_t^k(\X1n,\xx1n)\ge \frac 12 U_t^k(x,\x)]\nonumber\\\
\le \P[\frac{2^{(k-3)\vee 0}k}{2(k-1)}\{|\X1n|^{k-2}+|\xx1n|^{k-2}\} t \sup_{u\in [0,t]}|W_u^{0,t}|^2\ge \frac{1}{2(k-1)}U_t^k(x,\x)]\nonumber \\
+\bsum{i=2}^{k/2}\biggl\{\P[C_{k/2}^i \frac{2^{k-i-2}}{k-2i+1}\{|\X1n|^{k-2i}+|\xx1n|^{k-2i} \}t\sup_{u\in [0,t]}|W_u^{0,t}|^{2i}\ge \frac{1}{2(k-1)}U_t^k(x,\x) ]\nonumber\\
+ \P[C_{k/2}^i \frac{2^{k+i-2}}{k-i+1}\{|\X1n|^{k-i}+|\xx1n|^{k-i} \}t\sup_{u\in [0,t]}|W_u^{0,t}|^{i}\ge \frac{1}{2(k-1)}U_t^k(x,\x) ]\biggr\}\nonumber\\
:=P_1+\bsum{i=2}^{k/2}(P_2^i+P_3^i).\nonumber\\
\label{DECOMP_P}
\end{eqnarray}
From Proposition \ref{BROWNIAN_PROP} one gets that there exists $C_3:=C_3(k,n)\ge 1$  s.t.
\begin{eqnarray}
P_1&\le& C_3\exp\left(- C_3^{-1}\frac{U_t^k(x,\x)}{ \{|\X1n|^{k-2}+|\xx1n|^{k-2} \} t^2} \right),  \ \forall i\in \leftB 2,k/2\rightB,\nonumber\\
P_2^i&\le& C_3\exp\left (-C_3^{-1}\frac{U_t^k(x,\x)^{1/i}}{\{|\X1n|^{k-2i}+|\xx1n|^{k-2i} \}^{1/i} t^{1+1/i}  } \right), \nonumber\\
P_3^i&\le& C_3\exp\left (-C_3^{-1}\frac{U_t^k(x,\x)^{2/i}}{\{|\X1n|^{k-i}+|\xx1n|^{k-i} \}^{2/i} t^{1+2/i}  } \right). \label{CTR_PI}
\end{eqnarray}
Hence, plugging \eqref{CTR_PI} in \eqref{DECOMP_P}
we derive the claim  from \eqref{DECOMP_P}, \eqref{DECOMP_PROBA_DEV}, \eqref{CTR_EXP_AUX} and  \eqref{CTR_Y}.\hfill $\square$


\subsubsection{Auxiliary deviation estimates
}

Still from the Malliavin representation of the density given by \eqref{SCRIT_IPP}, when $\x_{n+1}-x_{n+1} $ is small, that is when for the degenerate component the starting and final points are close, we have to give estimates on $\P[Y_t\le\x_{n+1}-x_{n+1}]$ (small and moderate deviations).

\begin{PROP}
\label{ASYMP_PROB}
There exist constants  $(c_1,c_2):=(c_1,c_2)(n,k)$ s.t. for all $(x_{1,n},\xx1n)\in (\R^{n}\backslash\{0\})^2$, $\x_{n+1}>x_{n+1}$ and $t\ge 
2^{k+3}\frac {\x_{n+1}-x_{n+1}}{|x_{1,n}|^k+|\x_{1,n}|^k}$ :
\begin{equation}
\label{PICCOLO_1}
\P[Y_t\le \x_{n+1}-x_{n+1}]\le  c_1\exp\left(-c_2\frac{|x_{1,n}|^{2+k}+|\xx1n|^{2+k}}{\x_{n+1}-x_{n+1}} \right).
\end{equation}
For a given $K\ge 0$, if $t\ge \left[(\x_{n+1}-x_{n+1})\frac34 (64K)^k \right]^{2/(k+2)} $, and $|\X1n|\vee |\xx1n|\le Kt^{1/2}$, then there exist $(\bar c_1,\bar c_2):=(\bar c_1,\bar c_2)(n,k,K) $: 
\begin{equation}
\label{PICCOLO_2}
\P[Y_t\le \x_{n+1}-x_{n+1}]\le \bar c_1\exp\left(-\bar c_2\frac{t^{1+2/k}}{(\x_{n+1}-x_{n+1})^{2/k}} \right).
\end{equation}
\end{PROP}

\noindent \textit{Proof.} We first begin with the proof of \eqref{PICCOLO_1}. As in the previous sections, we can assume w.l.o.g. that $|\X1n|\ge |\xx1n| $.For $s\in [0,t]$, we define $\tilde X_s:=\X1n\frac{t-s}{t}+\xx1n \frac st+W_s^{0,t} $ (where $(W_s^{0,t})_{s\in [0,t]} $ is a standard $n$-dimensional Brownian Bridge on $[0,t] $), so that $Y_t=\int_0^t|\tilde X_s|^kds $. Let us also set $\tau_{|\X1n|/2}:=\inf\{s\ge 0: |\tilde X_s|\le |\X1n|/2 \} $. Consider now the event
$A:=\{\tau_{|x_{1,n}|/2}\le 2^k \frac{\x_{n+1}-x_{n+1}}{|x_{1,n}|^k}  \} $ and denote by $A^C$ its complementary. Observe that $\P[\int_0^t |\tilde X_s|^k ds\le \x_{n+1}-x_{n+1},A^C]=\P[\int_0^{2^k \frac{\x_{n+1}-x_{n+1}}{|x_{1,n}|^k}} \left(\frac{|x_{1,n}|}2 \right)^k ds<\int_0^t |\tilde X_s|^k ds \le \x_{n+1}-x_{n+1} ,A^C ]=0$. Thus, $\P[Y_t\le \x_{n+1}-x_{n+1}]=\P[Y_t\le \x_{n+1}-x_{n+1},A]\le \P[A]$. Now
\begin{eqnarray*}
\P[A]&\le & \P[\inf_{s\in[0,2^k \frac{\x_{n+1}-x_{n+1}}{|x_{1,n}|^k}]}|\tilde X_s|\le |\X1n|/2]  \\
&\le & \P[\inf_{s\in[0,2^k \frac{\x_{n+1}-x_{n+1}}{|x_{1,n}|^k}]}\left|\X1n\frac{t-s}{t}+\xx1n\frac st \right|  +\inf_{s\in [0,2^k \frac{\x_{n+1}-x_{n+1}}{|x_{1,n}|^k}]}(-|W_s^{0,t}|) \le |\X1n|/2]\\
&\le & \P[|\X1n|/2+\inf_{s\in [0,2^k \frac{\x_{n+1}-x_{n+1}}{|x_{1,n}|^k}] }(-\frac st)\{|\X1n|+|\xx1n|\}-\sup_{s\in [0,2^k \frac{\x_{n+1}-x_{n+1}}{|x_{1,n}|^k}]}|W_s^{0,t}| \le 0]\\
&\le & \P[|\X1n|(1/2- \frac{2^{k+1} (\x_{n+1}-x_{n+1})} {|x_{1,n}|^kt}) \le \sup_{s\in [0,2^k \frac{\x_{n+1}-x_{n+1}}{|x_{1,n}|^k}]}|W_s^{0,t}|]\\&\le & \P[|\X1n|/4\le \sup_{s\in [0,2^k \frac{\x_{n+1}-x_{n+1}}{|x_{1,n}|^k}]}|W_s^{0,t}|],
\end{eqnarray*}
recalling $|\X1n|\ge |\xx1n| $ and $t\ge  2^{k+3}\frac {\x_{n+1}-x_{n+1}}{|x_{1,n}|^k} $ for the last two inequalities.
From Proposition \ref{BROWNIAN_PROP} 
we obtain:
\begin{eqnarray*}
\P[Y_t \le \x_{n+1}-x_{n+1}]\le \P[A]\le c_1\exp\left(-c_2 \frac{|x_{1,n}|^{2+k}}{\x_{n+1}-x_{n+1}} \right),
\end{eqnarray*}
which from the assumption $|\X1n|\ge |\xx1n|$ gives \eqref{PICCOLO_1} up to a modification of $c_2$.

Let us now turn to \eqref{PICCOLO_2}. Introduce $I_\beta(t):=\int_0^t \I_{|\tilde X_{s}|^k\le \beta(\x_{n+1}-x_{n+1})} ds $ for a parameter $\beta>0$ to be fixed later on. Define the set $A_\beta:=\{ I_\beta(t)\ge t/4\} $.
Observe that
\begin{eqnarray*}
\P[\int_0^t |\tilde X_s|^kds\le \x_{n+1}-x_{n+1}, A_\beta^C]=\\
 \P[\int_0^t \I_{|\tilde X_s|^k> \beta(\x_{n+1}-x_{n+1})}|\tilde X_s|^kds \le \int_0^t |\tilde X_s|^kds \le \x_{n+1}-x_{n+1}, A_\beta^C  ]\\
\le\P[\beta(\x_{n+1}-x_{n+1})3t/4< \int_0^t |\tilde X_s|^kds\le \x_{n+1}-x_{n+1} ,A_\beta^C].
\end{eqnarray*}
Choosing $\beta=\frac 4{3 t}$ we get from the above inequality $\P[\int_0^t |\tilde X_s|^kds\le \x_{n+1}-x_{n+1}, A_\beta^C]=0 $. Hence,
\begin{eqnarray}
&&\P[\int_0^t |\tilde X_s|^kds\le \x_{n+1}-x_{n+1}]=\P[\int_0^t |\tilde X_s|^kds\le \x_{n+1}-x_{n+1},A_{\frac4{3t}}]\le \P[A_{\frac 4{3t}}]\nonumber\\
&&\le  \P[\int_0^t\I_{|\tilde X_s|^k\le \frac{3(\x_{n+1}-x_{n+1})}{4t}} ds> t/4]\le \P[ \int_{0}^t \I_{|\tilde X_s^1|^k\le \frac{3(\x_{n+1}-x_{n+1})}{4t}}ds>t/4]\nonumber\\
&&\le  \P[\int_0^t\I_{|x_1\frac{t-s}t+ \x_1 \frac st +B_s^{0,t}|\le c(x,\x,t,k)} ds>t/4],\quad c(x,\x,t,k):=\left( \frac{3(\x_{n+1}-x_{n+1})}{4t}\right)^{1/k},
\nonumber\\
&&\le \P[\int_0^{t/2}\I_{|x_1\frac{t-s}t+ \x_1 \frac st +B_s^{0,t}|\le c(x,\x,t,k)} ds>t/8]+\P[\int_{t/2}^t\I_{|x_1\frac{t-s}t+ \x_1 \frac st +B_s^{0,t}|\le c(x,\x,t,k)} ds>t/8]\nonumber\\
&&:=P_1+P_2,
\label{PREAL_OT}
\end{eqnarray}
where $(B_s^{0,t})_{s\in [0,t]} $ stands for a one-dimensional Brownian bridge on $[0,t] $. Observing that $(\bar B_s^{0,t}):=(B_{t-s}^{0,t})_{s\in [0,t]} $ is also a Brownian bridge, we get that
\begin{eqnarray*}
P_2&:=&\P[\int_0^{t/2}ds \I_{|x_1 \frac st+\x_1\frac{t-s}{t}+\bar B_{s}^{0,t}|\le c(x,\x,t,k)}ds>t/8]\\
&=&\P[\int_0^{t/2}ds \I_{|x_1 \frac st+\x_1\frac{t-s}{t}+B_{s}^{0,t}|\le c(x,\x,t,k)}ds>t/8].
\end{eqnarray*}
Since we assumed $|x_1|\vee |\x_1|\le Kt^{1/2}$, $|x_1|$ and $|\x_1| $ have at most the same magnitude so that $P_1$ and $P_2$ can be handled exactly in the same way.
Let us deal with $P_1$. The occupation time formula for semimartingales (see Chapter 6 in \cite{revu:yor:99}) yields
$$\int_0^{t/2}\I_{|x_1\frac{t-s}t+ \x_1 \frac st +B_s^{0,t}|\le c(x,\x,t,k)}  ds=\int_{-c(x,\x,t,k) }^{c(x,\x,t,k)} dz L_{t/2}^z,$$ where $L_{t/2}^z $ stands for the local time at level $z$ and time $t/2$ of the process $(x_1\frac{t-s}t+ \x_1 \frac st +B_s^{0,t})_{s\in [0,t]}$. From the definition of $P_1$ in \eqref{PREAL_OT}:
\begin{eqnarray}
P_1&\le& \P[\sup_{z\in \left[-c(x,\x,t,k),c(x,\x,t,k) \right]}L_{t/2}^z \times 2 c(x,\x,t,k)> \frac t{8}]\nonumber \\
&=&  \P[\sup_{z\in [-\frac{c(x,\x,t,k)}{t^{1/2}},\frac{c(x,\x,t,k)}{t^{1/2}} ]} \bar L_{1/2}^z> \frac{t^{1/2}}{16c(x,\x,t,k)} ], \label{PN}
\end{eqnarray}
where $\bar L_{1/2}^z $ stands for the local time at level $z$ and time $1/2$ for the scalar process
$$(\bar X_u)_{u\in [0,1]}:=\left(\frac{x_1}{t^{1/2}}(1-u)+\frac{\x_1}{t^{1/2}}u+\frac{B_{ut}^{0,t}}{t^{1/2}}\right)_{u\in [0,1]}\overset{({\rm law})}{=} \left(\frac{x_1}{t^{1/2}}(1-u)+\frac{\x_1}{t^{1/2}}u+B_u^{0,1}\right)_{u\in [0,1]}. $$
The last equality in \eqref{PN} is a consequence of the scaling properties of the local time. From Tanaka's formula for semimartingales $\bar L_{1/2}^z=|\bar X_{1/2}-z|-|\bar X_0-z|-\int_0^{1/2} {{\rm sgn}}(\bar X_s-z)d\bar X_s $. Denoting with a slight abuse of notation $(\frac{B_{ut}^{0,t}}{t^{1/2}})_{u\in [0,1]}=(B_u^{0,1})_{u\in [0,1]} $, we have the following differential dynamics for $\bar X_u	 $:
$$d\bar X_u=-\frac{x_1-\x_1}{t^{1/2}} du +dB_u^{0,1} =-\frac{\bar X_u-\x_1}{1-u} du +dB_u,$$
where $(B_u)_{u\in [0,1]}$ is a standard scalar Brownian motion.

Therefore, from equation \eqref{PN} and the usual differential dynamics for the Brownian bridge:
\begin{eqnarray*}
P_1\le \P[\frac{|\x_1-x_1|}{2t^{1/2}}+|B_{1/2}^{0,1}|\\
+\sup_{z\in [-\frac{c(x,\x,t,k)}{t^{1/2}},\frac{c(x,\x,t,k)}{t^{1/2}} ]}\left|\int_0^{1/2}{{\rm sgn}}(\bar X_s-z)(-\frac{x_1-\x_1}{t^{1/2}}ds+dB_{s}^{0,1})\right|  \ge \frac{t^{1/2}}{8c(x,\x,t,k)} ]\\
\le \P[\frac{|\x_1-x_1|}{t^{1/2}}+  |B_{1/2}^{0,1}|+\\
 \int_0^{1/2}ds\frac{|B_s^{0,1}|}{1-s}+ \sup_{z\in [-\frac{c(x,\x,t,k)}{t^{1/2}},\frac{c(x,\x,t,k)}{t^{1/2}} ]}|\int_0^{1/2}{{\rm sgn}}(\bar X_s-z) dB_s| \ge   \frac{t^{1/2}}{16c(x,\x,t,k)}]\\
 \le \P[2K+3\sup_{s\in [0,1/2]}|B_s^{0,1}|+\sup_{z\in [-\frac{c(x,\x,t,k)}{t^{1/2}},\frac{c(x,\x,t,k)}{t^{1/2}} ]}|\int_0^{1/2}{{\rm sgn}}(\bar X_s-z) dB_s| \ge \frac{t^{1/2}}{16c(x,\x,t,k)}].
\end{eqnarray*}
Now from the definition of $c(x,\x,t,k)$ in \eqref{PREAL_OT}, for $t\ge \left[(\x_{n+1}-x_{x+1})\frac34 (64K)^k \right]^{2/(k+2)} $ one has $ \frac{t^{1/2}}{16c(x,\x,t,k)}-2K\ge \frac{t^{1/2}}{32c(x,\x,t,k)}$. Thus
\begin{eqnarray*}
P_1\le \P[3\sup_{s\in [0,1/2]}|B_s^{0,1}|\ge \frac{t^{1/2}}{64 c(x,\x,t,k)}]\\
+\P[\sup_{z\in [-\frac{c(x,\x,t,k)}{t^{1/2}},\frac{c(x,\x,t,k)}{t^{1/2}} ]}|\int_0^{1/2}{{\rm sgn}}(\bar X_s-z) dB_s|\ge \frac{t^{1/2}}{64 c(x,\x,t,k)}].
\end{eqnarray*}
%
Setting for all $t\in [0,1/2] $, $M_t:=\int_0^t{{\rm sgn}}(\bar X_s-z)dB_s$, $M_t:=\tilde B_{\langle M\rangle_t} =\tilde B_t$ (i.e. $\tilde B$ is the Dambis-Dubbins-Schwarz Brownian motion associated to $M$).
Hence, from Proposition \ref{BROWNIAN_PROP} 
we derive the announced bound for $P_1$. Since $P_2$ can be handled in a similar way, the claim then follows from equation \eqref{PREAL_OT}.

\subsection{Final derivation of the upper-bounds in the various regimes} \label{upp-bd}
In this section we put together our previous estimates in order to derive the upper bounds of Theorem \ref{MTHM} in the various regimes.
\subsubsection{Derivation of the Gaussian upper bounds}

In this paragraph we assume $|x_{1,n}|\vee |\x_{1,n}|\ge Kt^{1/2} $ for $K$ large enough. We also suppose $\frac{|\x_{n+1}-x_{n+1}-ct(|\X1n|^k+|\xx1n|^k)|}{t^{3/2}(|\X1n|^{k-1}+|\xx1n|^{k-1})}\le \bar C $ where $c:=c(k)=2+\frac{2^{k-1}}{k+1}$ and $\bar C$ is fixed.
From Corollary \ref{REP_DENS} (representation of the density), Proposition \ref{CTR_W} (controls of the weight in the integration by part) and Lemma \ref{OFF_DIAG} (deviation bounds), we have that there exists $C:=C(n,k,K,\bar C)\ge 1$, s.t. setting $U_t^k(x,\x):=\x_{n+1}-x_{n+1}-\frac{2^{k-1}}{k+1}(|x_{1,n}|^k+|\x_{1,n}|^k)t $ as in Lemma \ref{OFF_DIAG} one has:
\begin{eqnarray}
\label{GAUSSIAN_BOUND}
p(t,x,\x)\le \frac{C\exp\left(-\frac{|\x_{1,n}-x_{1,n}|^2}{2t}-C^{-1}\frac{U_t^k(x,\x)^2}{(|x_{1,n}|^{k-1}+|\x_{1,n}|^{k-1})^2t^3} \right)}{t^{n/2+3/2}(|x_{1,n}|^{k-1}+|\x_{1,n}|^{k-1})}.
\end{eqnarray}


\begin{remark}
The above result means that the Gaussian regime holds if the final point $\xx1n $ of the degenerate component has the same order as the ``mean" transport term $m_t(x,\x):=x_{n+1}+\frac{2^{k-1}}{k+1}(|\X1n|^k+|\xx1n|^k) t $ (moderate deviations). A similar lower bound holds true, see  Lemma \ref{MIN_CPT_MET}.
\end{remark}

\subsubsection{Derivation of the heavy-tailed upper bounds}
We here assume $\frac{|\x_{n+1}-x_{n+1}-ct(|\X1n|^k+|\xx1n|^k)|}{t^{3/2}(|\X1n|^{k-1}+|\xx1n|^{k-1})}\ge \bar C $ where $c:=c(k)=2+\frac{2^{k-1}}{k+1}$ and $\bar C$ is as in the previous paragraph.

If $|x_{1,n}|\vee |\x_{1,n}|\le Kt^{1/2} $ ($K$ being as in the previous paragraph), then Corollary \ref{REP_DENS}, Proposition \ref{EST_NG} and Lemma \ref{OFF_DIAG} yield that there exists $C:=C(n,k)\ge 1$ s.t.
\begin{eqnarray}
\label{EST_OFF_DIAG_DENS}
p(t,x,\x)\le \frac{C}{t^{(n+k)/2+1}}\exp\left(-\frac{|\x_{1,n}-x_{1,n}|^2}{2t}-C^{-1}\frac{(U_t^k(x,\x))^{2/k}}{t^{1+2/k}}\right).
\end{eqnarray}
On the other hand if $|x_{1,n}|\vee |\x_{1,n}|\ge Kt^{1/2}  $, then
Corollary \ref{REP_DENS}, Proposition \ref{CTR_W} and Lemma \ref{OFF_DIAG} yield that there exists $\tilde C:=\tilde C(n,k)\ge 1$ s.t.
\begin{eqnarray*}
p(t,x,\x)&\le& \frac{\tilde C}{t^{n/2+3/2}(|x_{1,n}|^{k-1}+|\x_{1,n}|^{k-1})}\exp\left(-\frac{|\x_{1,n}-x_{1,n}|^2}{2t}-\tilde C^{-1}\frac{U_t^k(x,\x)^{2/k}}{t^{1+2/k}}\right)\\
&\le & \frac{\tilde C}{K^{k-1} t^{(n+k)/2+1}}\exp\left(-\frac{|\x_{1,n}-x_{1,n}|^2}{2t}-\tilde C^{-1}\frac{U_t^k(x,\x)^{2/k}}{t^{1+2/k}}\right).
\end{eqnarray*}
Hence, up to a modification of $C$, the control given by \eqref{EST_OFF_DIAG_DENS} holds for all off-diagonal cases.

\subsubsection{Moderate deviations of the degenerate component}
In this paragraph we suppose $0< \x_{n+1}-x_{n+1}\le K t^{1+k/2}$, for $K$ sufficiently small.  This means that the deviation of the degenerate component is small w.r.t. its characteristic time scale. From Corollary \ref{REP_DENS}, Propositions \ref{CTR_W} and \ref{EST_NG} and Proposition \ref{ASYMP_PROB} we derive similarly to the previous paragraph that there exists $C:=C(n,k,K)$ s.t.
\begin{eqnarray} \label{MOMODERATE_BOUND}
p(t,x,\x)&\le &\frac{C \exp\left(-\frac{|\x_{1,n}-x_{1,n}|^2}{2t}-C^{-1}\left\{ \frac{|x_{1,n}|^{2+k}+|\x_{1,n}|^{2+k}}{\x_{n+1}-x_{n+1}}+\frac{t^{1+2/k}}{(\x_{n+1}-x_{n+1})^{2/k}}\right\}\right)}{t^{(n+k)/2+1}}.
\end{eqnarray}

\subsection{Gaussian lower bound on the compact sets of the metric}
We conclude this section with a proof of  a lower bound for the density on the compact sets of the metric associated to the Gaussian regime in Theorem \ref{MTHM}. A similar feature already appears in the appendix of \cite{dela:meno:10}.
\begin{lemma}
\label{MIN_CPT_MET}
Assume that $|\X1n|\vee |\xx1n |\ge K t^{1/2}, \ K\ge K_0:=K_0(n,k)$ and that for a given $\bar C\ge 0$ we have  $ \frac{|\x_{n+1}-x_{n+1}-ct(|\X1n|^k+|\xx1n|^k)|}{t^{3/2}(|\X1n|^{k-1}+|\xx1n|^{k-1})}\le \bar C $ where $c:=c(k)$ is fixed. Then, there exists $C_{\ref{MIN_CPT_MET}}:=C_{\ref{MIN_CPT_MET}}(n,k,\bar C) $ s.t.
$$\frac{C_{\ref{MIN_CPT_MET}}}{(|\X1n|^{k-1}+|\xx1n|^{k-1})t^{3/2}}\le p_{Y_t}(\x_{n+1}-x_{n+1}).$$
\end{lemma}
\begin{remark}
The condition in the Lemma means that the deviation $\x_{n+1}-x_{n+1} $ has \textbf{exactly} the same order as the transport term $t(|x_{1,n}|^k+|\x_{1,n}|^k) $, up to a neglectable fluctuation corresponding to the variance of the Gaussian contribution in $Y_t$.
\end{remark}

\noindent \textit{Proof.} We assume w.l.o.g. that $\x_{n+1}-x_{n+1}-ct(|\X1n|^k+|\xx1n|^k)\ge 0$ and $|\xx1n|\ge |\X1n| $.
From \eqref{SCRIT_IPP} we recall:
\begin{eqnarray*}
p_{Y_t}(\x_{n+1}-x_{n+1})&=&\E[H_t\I_{Y_t\ge \x_{n+1}-x_{n+1}}],\\
 H_t&:=&H_t^1+H_t^2:=\gamma_{Y_t}^{-2}\langle D\gamma_{Y_t},DY_t\rangle_{L^2(0,t)} +\Gamma_{Y_t} LY_t.
\end{eqnarray*}
Recalling the chaos decomposition of $LY_t$ introduced in Proposition \ref{CTR_W}, see equation \eqref{INT_SKO_CAOS}, we get:
\begin{eqnarray*}
p_{Y_t}(\x_{n+1}-x_{n+1})&\ge &\E[H_t^2\I_{Y_t\ge \x_{n+1}-x_{n+1}}]-\E[|H_t^1|]
                                        \ge
                                        \E[\frac{I_1(g_0(.,t)}{\gamma_{Y_t}}\I_{Y_t\ge \x_{n+1}-x_{n+1}}]\\
                                        &&-\left\{\E\biggl[\frac{|\sum_{l=2}^{k} I_{l}(g_{l-1}(.,t)) |}{\gamma_{Y_t}}\biggr]+\frac{\bar C_1}{t|\xx1n|^k}\right\},
\end{eqnarray*}
using the bound for $\E[|H_t^1|] $ given by equation \eqref{CTR_HT1}, with $\bar C_1:=\bar C_1(n,k,1,K) $, in the last inequality. From equation \eqref{CTR_GL}, there exists $\bar C_2:=\bar C_2(n,k)$,
$\E[|\sum_{l=2}^{k} I_{l}(g_{l-1}(.,t)) |^2]^{1/2}\le \bar C_2 t^{3/2}|\xx1n|^{k-1}\sum_{l=2}^k \left( \frac{t^{1/2}}{|\xx1n|}\right)^{l-1}\le \frac{(k-1)\bar C_2}Kt^{3/2}|\xx1n|^{k-1}$, recalling $|\xx1n|\ge Kt^{1/2} $ for the last inequality.
Also $I_1(g_0(.,t))=\int_0^t \E[\bar M_1(u,t)] dW_u +R_0^t $ where $\E[|R_0^t|^2]^{1/2}\le \frac{\bar C_2}K t^{3/2}|\xx1n|^{k-1} $.
From \eqref{DEC_Y}, we write:
\begin{eqnarray*}
Y_t&=&\int_0^t du |m(u,t,\X1n,\xx1n)|^k+k\int_0^t du|m(u,t,\X1n,\xx1n)|^{k-2}\langle m(u,t,\X1n,\xx1n),W_u^{0,t}\rangle\\
&&+R_t^k(\X1n,\xx1n)
=:(m_t^k+G_t^k+R_t^k)(\X1n,\xx1n)=:m_t^k+G_t^k+R_t^k,
\end{eqnarray*}
for simplicity. Proposition \ref{CTR_W} then yields:
\begin{eqnarray*}
p_{Y_t}(\x_{n+1}-x_{n+1})&\ge & \E\biggl[\frac{\int_0^t \E[\bar M_1(u,t)] dW_u}{\gamma_{Y_t}}\I_{m_t^k+G_t^k+R_t^k\ge \x_{n+1}-x_{n+1} }\biggr]\\
&&-\left\{\frac{C_2 \bar C_2k t^{3/2}|\xx1n|^{k-1}}{K|\xx1n|^{2(k-1)}t^{3}}+\frac{\bar C_1}{K t^{3/2}|\xx1n|^{k-1}}\right\}\\
&:=&p_{Y_t,1}(\x_{n+1}-x_{n+1})-r_1(t,x,\x).
\end{eqnarray*}
From the martingale representation theorem and the above computations we identify $G_t^k=\int_{0}^t \E[\bar M_1(u,t)] dW_u$. Still from Proposition \ref{CTR_W} we get:
\begin{eqnarray*}
p_{Y_t}(\x_{n+1}-x_{n+1})&\ge & \E\biggl[\frac{G_t^k}{\gamma_{Y_t}}\I_{G_t^k+R_t^k\ge \x_{n+1}-x_{n+1} -m_t^k}\I_{|R_{t}^k|\le |G_t^k|/2}\biggr]\\&&-\biggl[\frac{\bar C_3\P[|R_t^k|>|G_t^k|/2]^{1/2}}{(|\X1n|^{k-1}+|\xx1n|^{k-1})t^{3/2}}+r_1(t,x,\x)\biggr]\\
&=&p_{Y_t,2}(\x_{n+1}-x_{n+1})-r_2(t,x,\x),
\end{eqnarray*}
where $\bar C_3:=\bar C_3(n,k) $. One easily gets that there exists $c:=c(k)>0,\ m_t^k:=m_t^k(\X1n,\xx1n)\ge c t(|\X1n|^k+|\xx1n|^k) $.Thus, setting $U_t^k(x,\x):=\xi_{n+1}-x_{n+1}-ct(|\X1n|^k+|\xx1n|^k) $ and recalling as well that $U_t^k(x,\x)\ge 0$, one obtains that on the event $\{G_t^k+R_t^k\ge U_t^k(x,\x), |R_{t}^k|\le |G_t^k|/2\}$, $G_t^k\ge 0 $. Hence:
\begin{eqnarray}
p_{Y_t}(\x_{n+1}-x_{n+1})&\ge &\E[\frac{G_t^k}{\gamma_{Y_t}}\I_{G_t^k-|R_t^k|\ge U_t^k(x,\x)\ge 0}\I_{|R_t^k|\le |G_t^k|/2}]-r_2(t,x,\x)\nonumber\\
&\ge& \E[\frac{G_t^k}{\gamma_{Y_t}}\I_{G_t^k\ge 2 U_t^k(x,\x)\ge 0}\I_{|R_t^k|\le G_t^k/2}]-r_2(t,x,\x)\nonumber\\
&\ge & \E[\frac{G_t^k}{3{\cal M}_t}\I_{G_t^k\ge 2 U_t^k(x,\x)\ge 0}\I_{|R_t^k|\le G_t^k/2}\I_{\gamma_{Y_t}\le 3{\cal M}_t}]-\frac{\bar C_3\P[\gamma_{Y_t}>3{\cal M}_t]^{1/2}}{|\xx1n|^{k-1}t^{3/2}}\nonumber\\
&&-r_2(t,x,\x)\nonumber\\
&\ge & \E[\frac{G_t^k}{3{\cal M}_t}\I_{G_t^k\ge 2 \bar Ct^{3/2}(|\X1n|^{k-1}+|\xx1n|^{k-1})}\I_{|R_t^k|\le G_t^k/2}\I_{\gamma_{Y_t}\le 3{\cal M}_t}]-r_3(t,x,\x) \nonumber\\
&\ge & \frac{C^{-1}2\bar C\P[G_t^k\ge 2 \bar Ct^{3/2}(|\X1n|^{k-1}+|\xx1n|^{k-1}),|R_t^k|\le G_t^k/2]}{3t^{3/2}(|\X1n|^{k-1}+|\xx1n|^{k-1})}\nonumber \\
&&-
\frac{\bar C_4\P[\gamma_{Y_t}>3{\cal M}_t]^{1/2}}{t^{3/2}|\xx1n|^{k-1}}-r_3(t,x,\x),
\label{PRE_MINO}
\end{eqnarray}
where we used that $U_t^k(x,\x)\le \bar Ct^{3/2}(|\X1n|^{k-1}+|\xx1n|^{k-1})$ for the last but one inequality (compact sets of the metric). The constant $C$ is the one appearing in \eqref{equivMT}. To conclude it suffices to prove that
\begin{eqnarray}
P&:=&\P[G_t^k\ge 2 Ct^{3/2}(|\X1n|^{k-1}+|\xx1n|^{k-1}),|R_t^k|\le |G_t^k|/2]]\ge \tilde C,\label{CTR_P}\\
|r_4(t,x,\x)|&:=& \frac{\bar C_4\P[\gamma_{Y_t}>3{\cal M}_t]^{1/2}}{t^{3/2}|\xx1n|^{k-1}}+r_3(t,x,\x)\le \frac{C^{-1} \bar C\tilde C}{3t^{3/2}(|\X1n|^{k-1}+|\xx1n|^{k-1})} \label{CTR_R}.
\end{eqnarray}
Indeed, plugging \eqref{CTR_P} and \eqref{CTR_R} into \eqref{PRE_MINO} gives the statement.
Let us first prove \eqref{CTR_P}. Write:
\begin{eqnarray*}
P&\ge &\P[G_t^k\ge 2 Ct^{3/2}(|\X1n|^{k-1}+|\xx1n|^{k-1})]\\
&&-\P[G_t^k\ge 2 Ct^{3/2}(|\X1n|^{k-1}+|\xx1n|^{k-1}),|R_t^k|> |G_t^k|/2]\\
   &\ge &\P[{\cal N}(0,1)\ge 2\check C]-\P[|R_t^k|\ge Ct^{3/2}(|\X1n|^{k-1}+|\xx1n|^{k-1})], \check C:=\check C(n,k).
\end{eqnarray*}
Thus, similarly to the proof of \eqref{CTR_R_LM} in Lemma \ref{LEMMA_MAL} we can show that there exists $\bar C_5:=\bar C_5(n,k)\ge 1 $ s.t. $\P[|R_t^k|\ge Ct^{3/2}(|\X1n|^{k-1}+|\xx1n|^{k-1})]\le \bar C_5\exp\left(-\bar C_5^{-1}\frac{|\X1n|^2+|\xx1n|^2}{t} \right) $. Under the current assumptions, using standard controls on the Gaussian distribution function, this gives \eqref{CTR_P} for  $\tilde C:=\tilde C(n,k) $ for $K$ large enough.

Recall now that $|r_4(t,x,\x)|\le \frac{1}{t^{3/2}|\xx1n|^{k-1}} \bigl(\frac{\bar C_1+C_2\bar C_2k}{K}+(\bar C_3+\bar C_4)\P[\gamma_{Y_t}>3{\cal M}_t]^{1/2}+\bar C_3\P[|R_t^k|>|G_t^k|/2]^{1/2}\bigr):=\sum_{i=1}^3r_{4i}(t,x,\x)$. Under the current assumptions, we derive that for $K$ large enough,  $r_{41}(t,x,\x)\le  \frac{C^{-1} \bar C \tilde C}{9t^{3/2}(|\X1n|^{k-1}+|\xx1n|^{k-1})}$. On the other hand, writing $\P[|R_t^k|>|G_t^k|/2]^{1/2}\le (\P[|R_t^k|>\frac{\hat C}2 (|\X1n|^{k-1}$ $+|\xx1n|^{k-1}) t^{3/2}]+\P[|G_t^k|\le \hat C (|\X1n|^{k-1}+|\xx1n|^{k-1}) t^{3/2}])^{1/2} $ we derive similarly to \eqref{CTR_R_LM} (see also the proof of the lower bound in Lemma \ref{LEMMA_COV}) that $r_{43}(t,x,\x)\le \frac{C^{-1} \bar C\tilde C}{9t^{3/2}(|\X1n|^{k-1}+|\xx1n|^{k-1})}$ taking $\hat C$ small enough. Eventually, the same control holds true for $r_{42}(t,x,\x) $, still from arguments similar to those used to derive \eqref{CTR_R_LM}. This concludes the proof. \hfill $\square $

\label{CHAINS-2}

\section{Potential Theory and PDEs}
\setcounter{equation}{0}
\setcounter{theorem}{0}
\label{PDE_SECTION}

In this section we are interested in proving Harnack inequalities for non-negative solutions to
\begin{equation}
\label{e1}
\L u(z) = 0,\qquad z=(x,t) \in \rnn,
\end{equation}
with $\L$ defined in \eqref{OPERATOR}. Specifically, we consider any open set $\O\subseteq\rnn$, and any $z\in\O$, and we aim to show that there exists a compact ${\K} \subset \O$ and a positive constant $C_{\K}$ such that
\begin{equation}\label{e-tt-1}
   \sup_{\K} u \le C_{\K} \, u(z),
\end{equation}
for every positive solution $u$ to $\L u = 0$. We say that a set $\big\{z_0, z_1, \dots, z_k\big\} \subset \O$ is a \emph{Harnack chain of lenght} $k$ if
\begin{equation*}
   u(z_{j}) \le C_j \, u(z_{j-1}), \qquad \text{for} \ j= 1, \dots, k,
\end{equation*}
for every positive solution $u$ of $\L u = 0$, so that we get
\begin{equation}\label{eq-harn-rep}
   u(z_{k}) \le C_{1}  C_{2} \dots C_{k} \, u(z_{0}).
\end{equation}
In order to construct Harnack chains, and to have an explicit lower bound for the densities considered in this article, we will prove \emph{invariant} Harnack inequalities w.r.t. a suitable Lie group structure. By exploiting the properties of homogeneity and translation invariance of the Lie group, we will find Harnack chains with the property that every $C_j$ in \eqref{eq-harn-rep} agrees with the constant $C_{\K}$ in \eqref{e-tt-1}. As a consequence we find $u(z_{k}) \le C_{{\K}}^{k} \, u(z_{0})$, and the bound will depend only on the lenght of the Harnack chain connecting $z_0$ to $z_k$.

\medskip

Let us now recall some basic notations concerning homogeneous Lie groups (we refer to the monograph \cite{LibroBLU} by Bonfiglioli, Lanconelli and Uguzzoni for an  exhaustive treatment). Let $\circ$ be a given group law on $\rnn$ and suppose that the map $(z, \z) \mapsto \z^{-1}\circ z$
is smooth. Then $\GG = (\rnn, \circ)$ is called a \emph{Lie group}. Moreover, $\GG$ is said {\it homogeneous} if there exists a family of dilations $\left(\d_{\ll}\right)_{\ll>0}$ which defines an automorphism of the
group, \emph{i.e.},
\begin{equation*}
\d_\ll(z \circ \z) = \left(\d_\ll z\right) \circ \left( \d_\ll \z\right), \quad \text{for all} \ z, \z
\in \rnn \ \text{and} \ \ll >0.
\end{equation*}

We also make the following assumption.

\begin{description}
  \item[{\rm \Ass{L}}] $\L$ is {\it Lie-invariant} with respect to the Lie group $\GG = \big( \rnn, \circ,$ $ (\d_\l )_{\l >0} \big)$, {\it i.e.}
\begin{description}
   \item[{\it i)}] $ Y_1, \dots,  Y_{n}$ and $Z$ are left-invariant with respect to the composition law of $\GG$, {\it i.e.}
\begin{equation*}
\begin{split}
      Y_j \left(u \left(\z\circ \cdot \right)
      \right) &  = \left(Y_j u\right) \left(\z\circ \cdot\right), \qquad j=1, \dots , n, \\
      Z \left(u \left(\z\circ \cdot \right)
      \right) & = \left(Z u\right) \left(\z\circ \cdot\right),
\end{split}
\end{equation*}
for every function $u \in C^\infty(\rnn)$, and for any $\z\in \rnn$;
   \item[{\it ii)}] $Y_{1},\dots,Y_{n}$ are $\d_{\l}$-homogeneous of degree one and $Z$ is $\d_{\l}$-homogeneous of degree two:
\begin{equation*}
\begin{split}
 Y_j \left(u \left( \d_{\ll}z \right) \right) & = \ll \left(Y_j u\right) \left(\d_{\ll}z\right), \qquad j=1, \dots , n, \\
 Z \left(u \left( \d_{\ll}z \right) \right) & = \ll^2 \left(Z u\right) \left(\d_{\ll}z\right),
\end{split}
\end{equation*}
for every function $u \in C^\infty(\rnn)$, and for any $z\in\rnn, \ll >0$.
\end{description}
\end{description}
To illustrate Property \Ass{L} we recall the Lie group structure of the Kolmogorov operator corresponding to $k=1$ in \eqref{PROC_2}.
\begin{example}{\sc (Kolmogorov operators)}\label{ex.Kolmo.op} \ $\L :=\frac 12 \Delta_{x_{1,n}} +\sum_{i=1}^n x_i \p_{x_2} - \p_t$. The Kolmogorov group is $\mathbb{K} = \left( \R^{n+2}, \circ, \d_\l \right)$, where
$$(x,t) \circ (\x, \t) = \big(
x_{1,n}+\x_{1,n}, x_{n+1}+\x_{n+1}  - \sum_{i=1}^n x_i \t, t+\t \big), \qquad \d_\l(x,t) = \left(\l
x_{1,n}, \l^3 x_{n+1}, \l^2 t \right).$$
Clearly, $\L$ can be written as in \eqref{OPERATOR} with $Y_i = \p_{x_i}, \ i\in \leftB 1,n \rightB$, and $Z = \sum_{i=1}^nx_i\p_{x_{n+1}} - \p_t$, and satisfies \emph{\Ass{L}}.
\end{example}

It is known that the composition law $\circ$ is always a sum with respect to the $t$ variable (see Propostion 10.2 in \cite{KogojLanconelli2}). Moreover, the family $\left(\d_{\ll} \right)_{\ll > 0}$ acts on $\rnn$ as follows:
$$
    \d_{\ll}(x_1, x_2, \dots , x_N,t) =  \left( \ll^{\sigma_1} x_1, \ll^{\sigma_2} x_2, \dots ,
    \ll^{\sigma_N} x_N, \ll^2 t \right), \quad \text{for every} \ (x,t) \in \rnn,
$$
where $\sigma = \left( {\sigma_1} , {\sigma_2}, \dots , {\sigma_N} \right) \in \N^N$ is a multi-index. The natural number $Q=\sum_{k=1}^{N}\sigma_k +2$ is called the {\it homogeneous dimension} of $\mathbb{G}$ with respect to $\d_{\lambda}$. We shall assume that $Q\geq 3$.
Observe that the diagonal decay of the heat kernel on the homogeneous Lie group is given by the characteristic time scale $t^{-(Q-2)/2} $. For the above example we have $Q=n+3+2$, matching the diagonal exponent in \eqref{DENS_K} $(Q-2)/2=(n+3)/2$.

Write the operator $\L$ as follows
\begin{equation*}
    \L = \sum_{i,j=1}^N a_{i,j}(x) \p_{x_i,x_j} + \sum_{j=1}^N b_{j}(x) \p_{x_j} - \p_t,
\end{equation*}
for suitable smooth coefficients $a_{i,j}$'s and $b_{j}$'s only depending on the vector fields $Y_0, \dots, Y_n$. As $n < N$, $\L$ is strictly degenerate, since the rank$\left( A(x) \right) \le n$ at every $x$ (here $A(x) := \left( a_{i,j}(x) \right)_{i,j \in \leftB 1,n \rightB}$). In Example \ref{ex.Kolmo.op} we see that rank$(A)$ never vanishes. We say that $\L$ is \emph{not totally degenerate} if
\begin{description}
  \item[{\rm \Ass{B}}] $\qquad\;$ for every $x\in \R^N$ there exists $\nu\in\R^N\setminus \{0 \}$ such that $\langle A(x) \nu, \nu \rangle > 0$.
\end{description}
This property holds for a more general class of operators. Indeed, if $\L$ satisfies \Ass{H} and \Ass{L}, then there exists a $\nu \in \R^N \setminus \{0 \}$ such that
\begin{equation} \label{e-B}
   \langle A(x) \nu, \nu \rangle > 0, \quad \text{for every} \ x \in \R^N.
\end{equation}
We refer to Section 1.3 in the monograph \cite{LibroBLU} for the proof of this statement.

Fix now $T>0$ and define $I:=[0,T]$. We call {\it diffusion trajectory} any absolutely continuous curve on $I$ such that
\begin{equation}\label{e-gdot}
   \g'(s)=\sum_{k=1}^{n} \omega_k(s) Y_k(\g(s)), \quad \text{for every}  \ s \in I,
\end{equation}
where $\omega_{1},\dots,\omega_{n}$ are piecewise constant real functions. A {\it drift trajectory} is any positively oriented integral curve of $Z$. We say that a curve $\g:[0,T] \to \rnn$ is {\it $\L$-admissible} if it is absolutely continuous and is a sum of a finite number of diffusion and drift trajectories.

Let $\O$ be any open subset of $\rnn$, and let $z_0\in\O$. We define the {\it attainable set} $\A_{z_0} := \overline {A_{z_0}}$ as the closure in $\O$ of the following set
\begin{equation}\label{e-Anew}
\begin{split}
 A_{z_0} = \big\{ & z\in\O : \text{there exists an $\L$-admissible path} \\
 & \qquad \g: [0,T] \to \O \
 \text{such that} \  \g(0)= z_0, \gamma(T)=z  \big\}.
 \end{split}
\end{equation}

The main result of the section is the following
\begin{theorem}\label{tt-1}
Let $\L$ be an operator in the form \eqref{e1} satisfying \emph{\Ass{H}} and \emph{\Ass{L}}, let $\O\subseteq\rnn$ be an open set, and let $z_0\in\O$. Then,
\begin{equation}\label{e-tt}
   \text{\it for every compact set ${\K} \subset$ {\rm Int}}\left(\A_{z_0}\right), \quad \sup_{\K} u \le C_{\K} \, u(z_0),
\end{equation}
for any non-negative solutions $u$ to $\L u = 0$ in $\O$. Here  $C_{\K}$ is a positive constant depending on $\O, {\K}, z_0$ and on $\L$.
\end{theorem}

We recall that a Harnack inequality for operators satisfying \Ass{H} and \Ass{B} is due to Bony (see \cite{Bony}). Another result analogous to Theorem \ref{tt-1} is given in \cite[Theorem 1.1]{CintiNystromPolidoro} by Cinti, Nystrom and Polidoro, assuming \Ass{L} and the following {\it controllability} condition:
\begin{description}
  \item[{\rm \Ass{C}}] for every $(x,t), (\x,\t) \in \rnn$ with $t>\t$, there exists an $\L$-admissible path $\g:[0,T] \to \rnn$
  such that $\g(0)=(x,t)$, $\g(T)=(\x,\t)$.
\end{description}
Our Theorem \ref{tt-1} improves Bony's one in that it gives an explicit geometric description of the set ${\K}$ in \eqref{e-tt}. Also, it is more general than the one in \cite{CintiNystromPolidoro}, since \Ass{L} and \Ass{C} imply \Ass{H} (see Proposition 10.1 in  \cite{KogojLanconelli2}).

\medskip
The proof of Theorem \ref{tt-1} is based on a general result from Potential Theory. In Section \ref{PotTh} we recall the basic results of Potential Theory needed in our work, then we apply them to operators $\L$ satifying \Ass{H} and \Ass{L}. We explicitly remark that condition \Ass{L} 
is not satisfied by the Kolmogorov operators \eqref{K-PROC} and \eqref{K-PROC_2}. 

\subsection{Potential Theory}
\label{PotTh}
For the first part of the section, we assume $\L $ to be a general abstract parabolic differential operator satisfying \Ass{B} and \Ass{L}. 

Let $\O$ be any open subset of $\rnn$. If $u:\O\to\mathbb{R}$ is a smoothfunction such that $\L u=0$ in $\O$, we say that $u$ is  $\L$-\emph{harmonic} in $\O$. We denote by $\mathcal{H}(\O)$ the linear space of functions which are $\L$-harmonic in $\O$.

Let $V$ be a bounded open subset of $\mathbb{R}^{N+1}$  with Lipschitz-continuous boundary. We say that $V$ is $\L$-\emph{regular} if, for every $z_0\in \partial V$, there exists a neighborhood $U$ of $z_0$ and a smooth function $w: U\to \R$ satisfying
\begin{equation*}
    w(z_0)=0,\quad \L w (z_0)<0, \quad w>0\;\text{in}\;\overline{V}\cap U \setminus \{z_0 \}.
\end{equation*}
Note that the function $\psi(x,t)=\frac{1}{2}+\frac{1}{\pi}\arctan t$ verifies
\begin{equation}\label{e-I}
    0\leq \psi \leq 1, \quad \L \psi < 0 \ \text{ in } \ \rnn.
\end{equation}
As a first consequence of \eqref{e-I}, the classical Picone's maximum principle holds on any bounded open set $\O\subset\rnn$. Precisely, if $u\in C^2(\O)$ satisfies
\begin{equation*}
 \L u \ge 0 \quad \text{in} \ \O, \qquad \limsup_{z\to\z} u(z) \le 0 \quad \text{for every} \ \z\in\partial\O,
\end{equation*}
then $u\le 0$ in $\O$ (see e.g. Bonfiglioli and Uguzzoni \cite{BonfiglioliUguzzoni3}). Then, for every $\L$-regular open set $V\subset\mathbb{R}^{N+1}$, and for any $\varphi\in C(\partial V)$ there exists a unique function $H_\varphi^V$ satisfying
\begin{equation}\label{H.varphi.V}
 H_\varphi^V\in\mathcal{H}(V), \qquad \lim_{z\to \z}H_\varphi^V(z)=\varphi(\z) \quad
 \text{for every }\, \z\in\partial V.
\end{equation}
Moreover, $H_\varphi^V\geq 0$ whenever $\varphi\geq 0$ (see Bauer \cite{Bauer} and Constantinescu and Cornea \cite{CC}). Hence, if $V$ is $\L$-regular, for every fixed $z\in V$ the map $\varphi\mapsto H_\varphi^V(z)$ defines a linear positive functional on $C(\partial V,\mathbb{R})$. Thus, the Riesz representation theorem implies that there exists a Radon measure $\mu_z^V$, supported in $\partial V$, such that
\begin{equation}\label{harmonic.measure.L.regular}
 H_\varphi^V(z)=\int_{\partial V}\varphi(\zeta)\,d \mu_z^V(\zeta),\quad\text{ for every }\,\varphi\in
 C(\partial V,\mathbb{R}).
\end{equation}
We will refer to $\mu_z^V$ as the $\L$-\emph{harmonic measure} defined with respect to $V$ and $z$.

\medskip

A lower semi-continuous function $u:\O\to\,]-\infty,\infty]$ is said
to be $\L$-\emph{superharmonic} in $\O$ if $u<\infty$ in a dense subset of $\O$ and if
$$u(z)\geq\int_{\partial V}u(\zeta)\,d\mu_z^V(\zeta),$$
for every open $\L$-regular set $V\subset\overline{V}\subset\O$ and for every $z\in V$.
We denote by $\overline{\mathcal{S}}(\O)$ the set of $\L$-superharmonic functions in $\O$, and by $\overline{\mathcal{S}}^+(\O)$ the set of the functions in $\overline{\mathcal{S}}(\O)$
which are non-negative.
A function $v:\O\rightarrow [-\infty,\infty[$ is said to be $\L$-\emph{subharmonic} in $\O$ if
$-v\in\overline{\mathcal{S}}(\O)$ and we write $\underline{\mathcal{S}}(\O):=-\overline{\mathcal{S}}(\O)$.
Since the collection of $\mathcal{L}$-regular sets is a basis for the Euclidean topology (as we will see in a moment), we have $\overline{\mathcal{S}}(\O)\cap\underline{\mathcal{S}}(\O)=\mathcal{H}(\O)$.

 This last property and Picone's maximum principle are the main tools in order to show the following criterion of $\L$-superharmonicity for functions of class $C^2$ (a proof can be found in the monograph \cite[Proposition 7.2.5]{LibroBLU}).
\begin{remark}\label{smooth.subarm.funct}
 Let $u\in C^2(\O)$. Then $u$ is $\L$-superharmonic if and only if $\L u \leq 0$ in $\O$.
\end{remark}

\medskip

With the terminology of Potential Theory (we refer to the monographs \cite{Bauer, CC}), the map
 $\rnn\supseteq\O\mapsto\mathcal{H}(\O)$ is said \emph{harmonic sheaf} and $(\rnn,\mathcal{H})$ is said
 \emph{harmonic space}. Since the constant functions are $\L$-harmonic, the last statement is a consequence of the following properties:
 \begin{itemize}
  \item[-] the $\L$-regular sets form a basis for the Euclidean topology (by \eqref{e-B}, $\L$ is a not totally degenerate operator%
  , so that this statement is a consequence of \cite[Corollaire 5.2]{Bony});
  \item[-] $\mathcal{H}$ satisfies \emph{the Doob convergence
   property}, \emph{i.e.}, the pointwise limit $u$ of any increasing
   sequence $\{u_n\}_n$ of $\L$-harmonic functions, on any open set $V$, is
   $\L$-harmonic whenever $u$ is finite in a dense set $T\subseteq V$ (as in
   \cite[Proposition 7.4]{KogojLanconelli2}, we can rely on the weak Harnack inequality due to Bony stated in \cite[Theoreme 7.1]{Bony});
  \item[-] the family $\overline{\mathcal{S}}(\rnn)$
   \emph{separates the points of} $\rnn$, \emph{i.e.}, for every $z,\zeta\in\rnn$, $z\neq\zeta$, there   exists $u\in\overline{\mathcal{S}}(\rnn)$ such that $u(z)\neq u(\zeta)$.
 \end{itemize}

This last \emph{separation property} is proved in Lemma \ref{separation.property}. 
We will in fact show a stronger result: actually, the family $\overline{\mathcal{S}}^+(\mathbb{R}^{N+1})\cap C(\mathbb{R}^{N+1})$ separates the points of $\mathbb{R}^{N+1}$. A harmonic space $(\mathbb{R}^{N+1}, \mathcal{H})$ satisfying this property is said to be a $\mathfrak{B}$-\emph{harmonic space}.

In order to prove the separation property we use a fundamental solution $\G$ of $\L$. To prove the existence of a fundamental solution we now rely on condition \Ass{H} that we assumed to be in force through the paper.
We recall that a fundamental solution is a function $\Gamma$ with the following properties:
\begin{description}
\item[{\it i)}] the map $(z,\z)\mapsto\Gamma(z,\z)$ is defined, non-negative and smooth away from the set $\{(z,\z)\in \rnn\times \rnn: z\neq \zeta \}$;
\item[{\it ii)}] for any $z \in \rnn, \Gamma(\cdot,z)$ and $\Gamma(z, \cdot)$
    are locally integrable;
\item[{\it iii)}] for every $\phi \in C_0^\infty(\rnn)$ and $z \in \rnn$ we have
  \begin{equation*}
    \L \!\int_{\rnn} \Gamma(z,\z)\phi(\z)\, d \z =
    \int_{\rnn} \Gamma(z,\z)\L \phi(\z)\, d \z = - \phi(z);
  \end{equation*}
\item[{\it iv)}] $\L \Gamma( \cdot, \z) = - \d_\z$ (Dirac measure supported at $\z$);
\item[{\it v)}] if we define $\Gamma^*(z,\z):= \Gamma(\z,z)$, then $\Gamma^*$ is
    the fundamental solution for the formal adjoint $\L^*$ of $\L$, satisfying the dual statements of $\textit{iii)}$, $\textit{iv)}$;
\item[{\it vi)}] $\Gamma(x,t,\x, \t) = 0$ \ if  \ $t < \t$.
\end{description}

\begin{remark} \label{REM-GAMMA}
Assumption \Ass{H} implies the existence of a smooth density $p(t,\x,x)dx:=\P_{\x}[X_t\in dx],\ t>0,$ for the process $(X_t)_{t\ge 0} $ associated to $\L $ see e.g. Stroock \cite{stro:83} or Nualart \cite{nual:95}. Actually,
$$\Gamma(x,t,\x,\tau):=p(t-\tau,\xi,x)$$
is a fundamental solution for $\L$ in the above sense. Indeed $p$ satisfies the Kolmogorov equation $\L p=0 $, in $ \R^{N+1}\backslash \{(\xi,\tau)\}  $. We refer to Bonfiglioli and Lanconelli  \cite{BonfiglioliLanconelli1} for a purely analytic proof of existence of fundamental solutions for operators satisfying \Ass{H}, \Ass{L}.
\end{remark}

If condition \Ass{L} holds, then we also have:
\begin{description}
\item[{\it vii)}] $\Gamma(z,\z)=\Gamma(\alpha\circ z,\alpha \circ \z)$ \ for every $\alpha,z,\z\in\rnn$, $z\neq \z$;
\item[{\it viii)}]  $\Gamma(\d_\lambda(z),\d_\lambda(\z))=\lambda^{-Q+2}\Gamma(z,\z), \qquad z,\z\in\rnn,\,z\neq \z, \;\lambda>0$.
\end{description}

We next prove the separation property for $\L$ by adapting the argument in \cite[Proposition 7.1]{CintiLanconelli}.

\begin{lemma}\label{separation.property}
 For every $z_1,z_2\in\mathbb{R}^{N+1}$, $z_1\neq z_2$, there exists a function $u\in\overline{\mathcal{S}}^+(\mathbb{R}^{N+1})\cap C(\mathbb{R}^{N+1})$ such that
 $u(z_1)\neq u(z_2)$.
\end{lemma}
\noindent {\it Proof.}\ Let us denote $z_i=(x_i,t_i)$ for $i=1,2$.
First we suppose that $t_1<t_2$. The properties of $\Gamma$ yield that there exists $z_0=(x_0,t_0)$ with $t_0> 0$ such that $\Gamma(z_0,0)>0$. On the other hand,
since \Ass{H} and \Ass{L} yield \Ass{B}, there exists a $\L$-regular open set $V_0$ containing the origin, a small $r_0>0$ and a large $\l_0 >1$ such that
\begin{equation}\label{eq.prop.V1}
 U_{r_0} \subseteq V_0 \subseteq \delta_{\lambda_0}(U_{r_0}),\qquad U_{r_0}=\{(x_1,\ldots,x_N,t)\in\rnn : |x_i|<r_0, |t|<r_0\}.
\end{equation}
By the smoothness of $\Gamma$, there exists $\e>0$ such that $\Gamma>0$ in the set $z_0\circ U_{\e}$. 
For a fixed $\lambda\in \left]0,\sqrt{\frac{t_2-t_1}{2(t_0+\e)}}\right[$ and a non-negative function $\varphi\in C_0^{\infty}(z_2\circ \bigl( \delta_\lambda(z_0\circ U_\e)\bigr)^{-1}\cap\{t<t_2\})$, we set
\begin{equation}\label{def.pot.varphi}
 u_\varphi(z)=\int_{\mathbb{R}^{N+1}}\Gamma(z,\z)\,\varphi(\zeta)\,d \zeta,\qquad z\in \rnn.
\end{equation}
Hence, we obtain $u_\varphi\in C^\infty(\mathbb{R}^{N+1})$, $u_\varphi\geq 0$ and $\L u_\varphi=-\varphi\leq 0$, so that, by Remark \ref{smooth.subarm.funct}, $u_\varphi\in\overline{\mathcal{S}}(\mathbb{R}^{N+1})$. Moreover the choice of $\varphi$ implies that $u_\varphi(z_1)=0$ and $u_\varphi(z_2)>0$.

In the case $t_1=t_2$, $x_1\neq x_2$, we consider the sequence
\begin{equation}\label{eq-seq.L.balls}
 \O_n(z_2)=\left\{ \zeta\in\mathbb{R}^{N+1}:\Gamma(z_2, \z)>n^{Q-2}\right\}\!, \qquad n\in\N.
\end{equation}
We note that $\O_n(z_2)$ shrinks to $\{z_2\}$ as $n\to\infty$, by property \textit{viii)} of the fundamental solution. For any $\varphi_n\in C_0^{\infty}(\O_n(z_2))$ such that $\int\varphi_n = 1$ and $\varphi_n\geq 0$, we define $u_{\varphi_n}$ as in \eqref{def.pot.varphi}. Then, $u_{\varphi_n}$ is a smooth non-negative function in $\mathbb{R}^{N+1}$ satisfying $\L u_{\varphi_n}\leq 0$, and so $u_{\varphi_n}$ is $\L$-superharmonic. It holds
\begin{gather*}
   u_{\varphi_n}(z_2)
   =\int_{\rnn}\Gamma(z_2,\z)\,\varphi_n(\zeta)\,d\zeta
    \geq n^{Q-2}\quad\text{ for every }\,n\in\mathbb{N};\\
    u_{\varphi_n}(z_1)
   \leq\max_{\zeta\in\overline{\O_1(z_2)}}\Gamma(z_1,\z)=C,
  \end{gather*}
  where $C$ is a real positive constant independent of $n$. This ends the proof.
  \hfill $\square$

\medskip

We summarize the above facts in the following

\begin{proposition} \label{r-harm}
Let $\L$ be an operator in the form \eqref{e1} and assume that \emph{\Ass{H}} and \emph{\Ass{L}} are satisfied. The map $\mathcal{H}$ which associates any open set $\O \subseteq \rnn$ with the linear space of the $\L$-harmonic functions in $\O$ is a \emph{harmonic sheaf}, and $(\rnn,\mathcal{H})$ is a  $\mathfrak{B}$-\emph{harmonic space}.
\end{proposition}

\medskip

A remarkable feature of a $\mathfrak{B}$-harmonic space is that the \emph{Wiener resolutivity theorem} holds (see \cite{Bauer, CC}). In order to state it, we introduce some additional notations. We recall that if $\O\subset\mathbb{R}^{N+1}$ is a bounded open set, then an extended real function $f:\partial\O\to[-\infty,\infty]$ is called \emph{resolutive} if
 \begin{equation*}
  \inf\overline{\mathcal{U}}_f^\O=
 \sup\underline{\mathcal{U}}_f^\O=:H_f^\O\in\mathcal{H}(\O),
 \end{equation*}
 where
 \begin{gather*}
  \begin{split}
   \overline{\mathcal{U}}_f^\O
  &=\big\{u\in\overline{\mathcal{S}}(\O)
   :\inf_{\O}u>-\infty\,\text{ and }\,\liminf_{z\to\zeta}u(z)\geq f(\zeta),
   \,\forall\,\zeta\in\partial\O\big\},\\
    \underline{\mathcal{U}}_f^\O
  &=\big\{u\in\underline{\mathcal{S}}(\O)
   :\sup_{\O}u<\infty\,\text{ and }\,\limsup_{z\to\zeta}u(z)\leq f(\zeta),
   \,\forall\,\zeta\in\partial\O\big\}.
  \end{split}
 \end{gather*}
 We say that $H_f^\O$ is the \emph{generalized solution in the sense of Perron-Wiener-Brelot} to the problem
 \begin{equation*}
 u\in\mathcal{H}(\O), \qquad  u=f\quad\text{on }\,\partial\O.
\end{equation*}
 The Wiener resolutivity theorem yields that any $f\in C(\partial\O,\mathbb{R})$ is resolutive. The map $C(\partial\O,\mathbb{R})\ni f\mapsto H_f^\O(z)$ defines a linear positive functional for every $z\in\O$. Again, there exists a Radon measure  $\mu_z^\O$ on $\partial\O$ such that
 \begin{equation}\label{harmonic.measure}
  H_f^\O(z)=\int_{\partial\O}f(\zeta)\,\mathrm{d}\mu_z^\O(\zeta).
 \end{equation}
We call $\mu_z^\O$ the $\L$-\emph{harmonic measure} relative to $\O$ and $z$, and when $\O$ is $\L$-regular this definition coincides with the one in \eqref{harmonic.measure.L.regular}.
Finally, a point $\zeta\in\partial\O$ is called $\L$-\emph{regular} for $\O$ if
\begin{equation}\label{e-boundary}
    \lim_{\O\ni z\to\zeta}H_f^\O(z)=f(\zeta),
  \quad\text{ for every }\,f\in C(\partial\O,\mathbb{R}).
 \end{equation}
Obviously, $\O$ is $\L$-regular if and only if every $\zeta\in\partial\O$ is $\L$-regular.

 \subsection{Harnack inequalities}
Let $\O \subset\rnn$ be an open set. A closed subset $F$ of $\O$ is called an
\emph{absorbent set} if, for any $z\in F$ and any $\L$-regular neighborhood $V\subset \overline{V}\subset\O$ of $z$, it holds $\mu_z^V(\partial V \setminus F)=0$. For any given
$z_0\in\O$ we set
$$\mathscr{F}_{z_0}=\{ F \subset\O : F \ni z_0, F\text{ is an absorbent set}\}.$$
Then,
\begin{equation}\label{e-absorbent}
   \O_{z_0}=\bigcap_{F\in\mathscr{F}_{z_0}}F
\end{equation}
is the \emph{smallest} absorbent set containing $z_0$. The Potential Theory provides us with the following Harnack inequality.
{\it Let $(\mathbb{R}^{N+1},\mathcal{H})$ be a $\mathfrak{B}$-harmonic space, let $\O$ be an open subset of $\rnn$ and let $z_0\in\O$. Then,}
\begin{equation}\label{e-Hp}
	\text{\it for every compact set $K \subset$ {\rm Int}}\left(\O_{z_0}\right), \quad \sup_K u \le C_K \, u(z_0),
\end{equation}
{\it for any non-negative function $u \in \mathcal{H}(\O)$. Here  $C_K$ is a positive constant depending on $\O, K, z_0$.}
We refer to Theorem 1.4.4 in \cite{Bauer} and Proposition 6.1.5 in \cite{CC}. Proposition \ref{r-harm} implies that \eqref{e-Hp} applies to our operator $\L$. We summarize the above argument in the following

\begin{proposition} \label{t-Harnack.ineq.absorbent.set}
Let $\L$ be an operator in the form \eqref{e1} satisfying \emph{[H]} and \emph{[L]}, let $\O\subseteq\rnn$ be an open set, and let $z_0\in\O$. Then,
\begin{equation*}
    \text{\it for every compact set $K \subset$ {\rm Int}}\left(\O_{z_0}\right), \quad \sup_K u \le C_K \, u(z_0),
\end{equation*}
for any non-negative solutions $u$ to $\L u = 0$ in $\O$. Here $C_K$ is a positive constant depending on $\O, K, z_0$ and on $\L$.
\end{proposition}


In order to prove Theorem \ref{tt-1} we give the following

\begin{lemma}\label{th.A.subset.Omega}
Let $\L$ be an operator as in \eqref{e1} satisfying \emph{[H]} and \emph{[L]}, and let $\O$ be an open subset of $\rnn$. For any given $z_0\in\O$, we have $\A_{z_0}\subseteq \O_{z_0}$ with $\A_{z_0} $ defined in \eqref{e-Anew}.
\end{lemma}

\noindent {\it Proof.}\ Since $\O_{z_0}$ is a closed set, and $\A_{z_0}$ is the closure of the set $A_{z_0}$ defined in \eqref{e-Anew}, it is sufficient to show that $A_{z_0}\subseteq \O_{z_0}$. By contradiction, assume that $\overline{z}\in A_{z_0}\!\setminus \O_{z_0}$. Then, there exists an $\L$-admissible path $\g: [0,T] \to \O$ such that $\g(0)=z_0, \g(T)=\overline{z}$.

We set
\begin{equation*}
    t_1:=\inf\{t>0 : \g(]t,T])\cap\O_{z_0}=\emptyset\}.
\end{equation*}
Note that, since $\O\setminus\O_{z_0}$ is an open set containing $\overline{z}$ and $\g$ is a continuous curve, there exists an open neighborhood $U\subseteq\O$ of $\overline{z}$ such that $U\cap \O_{z_0}=\emptyset$, and a positive $\s$ satisfying $\g(]T-\s,T])\subseteq U$. Hence, $t_1 \in [0,T[$ is well defined and we have $\g(t)\notin \O_{z_0}$ for every $t\in ]t_1,T]$. Again, by the continuity of $\g$, we have
$$z_1=\g(t_1) \in \O_{z_0}.$$

\noindent Let $V\subset \overline{V}\subset\O$ be a $\L$-regular neighborhood of $z_1$ with $\overline{z}\notin \overline{V}$. Arguing as above, we can find $t_2\in ]t_1,T[$ such that $\gamma([t_1, t_2[)\subset {V}$ and $z_2=\g( t_2)\in\partial V$. Consider any neighborhood $W$ of $z_2$, such that $W \subset \O\setminus\O_{z_0}$. Let $\varphi\in C(\partial V)$ be any non-negative function, supported in $W \cap \partial V$, and such that $\varphi(z_2)>0$. Recalling that the harmonic function $H_\varphi^V$ is non-negative, we aim to show that
\begin{equation}\label{eq1-Omega.equiv.A}
    H_\varphi^V(z_1)> 0.
\end{equation}
By contradiction, we suppose that $H_\varphi^V$ vanishes at $z_1$. In other terms, $H_\varphi^V$ attains its minimum value at $z_1$, then Bony's minimum principle implies $H_\varphi^V\equiv 0$ in $\gamma([t_1, t_2[)$. As a consequence, since $H_\varphi^V$ satisfies \eqref{H.varphi.V},
\begin{equation}\label{eq2-Omega.equiv.A}
    \lim_{t \to t_2^-}H_\varphi^V(\g(t))=0.
\end{equation}
On the other hand, by the choice of $\varphi$
\begin{equation*}
    \lim_{V \ni z \to z_2}H_\varphi^V(z)=\varphi( z_2)>0.
\end{equation*}
This contradicts \eqref{eq2-Omega.equiv.A} and proves \eqref{eq1-Omega.equiv.A}. By using representation \eqref{harmonic.measure.L.regular} of $H_\varphi^V$ in terms of the $\L$-harmonic measure, \eqref{eq1-Omega.equiv.A} reads as follows
\begin{equation}\label{eq3-Omega.equiv.A}
    H_\varphi^V(z_1) = \int_{\partial V \cap W}\varphi(\zeta)\,d \mu_{z_1}^V(\zeta)>0, \qquad \text{ then }\quad \mu_{z_1}^V(\partial V \cap W) >0.
\end{equation}
On the other hand, $z_1$ belongs to the absorbent set $\O_{z_0}$, so that $\mu_{z_1}^V(\partial V \setminus \O_{z_0})=0$. But this clashes with \eqref{eq3-Omega.equiv.A}, being $W \subseteq \O \setminus \O_{z_0}$. This accomplishes the proof. \hfill $\square$

\medskip

\noindent {\it Proof of Theorem \ref{tt-1}.}\ It is a plain consequence of Proposition \ref{t-Harnack.ineq.absorbent.set} and Lemma \ref{th.A.subset.Omega}
\hfill $\square$


\medskip

As the following proposition shows, we are able to give a complete characterization of the set $\O_{z_0}$ if $\A_{z_0}$ is an absorbent set as well.

\begin{proposition}\label{cor.Aass.equiv.Omega}
Let $\L$ be an operator as in \eqref{e1} satisfying \emph{[H]} and \emph{[L]}, let $\O\subseteq \rnn$ be an open set, and let $z_0\in\O$. If $\A_{z_0}$ is an absorbent set, then $\A_{z_0}\equiv \O_{z_0}$.
\end{proposition}
\noindent {\it Proof.}\ The claim directly follows from Lemma \ref{th.A.subset.Omega}, recalling the definition of $\O_{z_0}$. \hfill $\square$

\medskip

The first statement in next proposition is a classical result in abstract potential theory (see e.g.\ \cite[Theorem 1.4.1]{Bauer} and \cite[Proposition 6.1.1]{CC}). For the convenience of the reader, we explicitly give here its simple proof.

\begin{proposition}\label{pr.cond.suff.Aass}
Let $\L$ be an operator as in \eqref{e1} satisfying \emph{[H]} and \emph{[L]}, let $\O\subseteq \rnn$ be an open set, and let $z_0\in\O$. Assume that there exists a solution $u \ge 0$ to $\L u= 0$ in $\O$ such that $u \equiv 0$ in $\A_{z_0}$ and $u>0$ in $\O\setminus \A_{z_0}$. Then $\A_{z_0}$ is an absorbent set, and $\A_{z_0}\equiv \O_{z_0}$.
\end{proposition}
\noindent {\it Proof.}\ Since $u$ is continuous and non-negative,
\begin{equation*}
 \A_{z_0}=\{z\in\O: u(z)\le 0\}
\end{equation*}
is a closed subset of $\O$. Let $z\in \A_{z_0}$, and let $V\subset \overline{V}\subset\O$ be a $\L$-regular neighborhood of $z$. As $u\in\mathcal{H}(\O)$, we have
\begin{equation*}
    0\geq u(z)=\int_{\partial V}u(\zeta)\,d\mu_z^V(\zeta) \geq 0, \qquad \text{ so that}\quad
    \mu_z^V\!\left(\partial V \setminus \A_{z_0}\right)=0.
\end{equation*}
Hence $\A_{z_0}$ is an absorbent set. The last statement plainly follows from Proposition \ref{cor.Aass.equiv.Omega}. \hfill $\square$



\label{PROOFS}
\subsection{Lifting and Harnack inequalities}
\label{CHAINS}
We first consider the PDE \eqref{K-PROC} for $k=2$. Note that, in this case, it is equivalent to \eqref{K-PROC_2}, and reads as follows
\begin{equation}\label{K-PROC-k=2}
	\L = \frac12 \Delta_{x_{1,n}} + |x_{1,n}|^2 \p_{x_{n+1}} - \p_t.
\end{equation}
It is homogeneous with respect to the following dilation
\begin{equation}\label{e-dil-K}
  \d_\lambda(x,t)=\big(\lambda x_{1,n}, \lambda^{4} x_{n+1}, \lambda^2 t\big).
\end{equation}

Even if $\L$ does not satisfy \Ass{L}--\emph{i)}, it has a fundamental solution $\Gamma$ which shares several properties of the usual heat kernels. We remark that, since $\L$ does not satisfy the controllability condition \Ass{C}, the support of $\G$ is strictly contained in the half space $\big\{ t < \t \big\}$. 

We next show that $\L$ can be \emph{lifted} to a suitable operator $\tilde \L$ in the form \eqref{e1} satisfying both \Ass{H} and \Ass{L}. By adding a new variable $y= y_{1,n} \in \R^n$, we define the following vector fields on $\R^{2n+2}$
\begin{equation} \label{e-LL}
\tilde Y_i = Y_i = \partial_{x_i}, i \in \leftB 1,n \rightB, \qquad \widetilde{Z}  = |x_{1,n}|^2 \p_{x_{n+1}} + \sum_{i=1}^n x_i \partial_{y_i} - \p_t.
\end{equation}
Clearly, if we denote $v(x,y,t) = u(x,t)$ for any $u \in C^\infty(\R^{n+2})$, we have
\begin{equation*}
  \tilde Y_i v (x,y,t)  = Y_i u(x,t), \quad \forall i \in \leftB 1,n \rightB, \qquad \tilde Z v (x,y,t) = Z u(x,t),
\end{equation*}
then, if we consider the \emph{lifted} operator $\tilde \L  = \tfrac12 \sum_{i=1}^n \tilde Y_i^2 +\tilde Z$, we find $\tilde \L v (x,y,t) = \L u(x,t)$.
By a standard procedure (see e.g., \cite[Chapter 1]{LibroBLU}), we explicitly write the group law $\circ$ of the homogeneous Lie group $\GG = \big( \R^{2n+2}, \circ, (\dd_\l)_{\l >0} \big)$ such that $\tilde \L$ is $\GG$-Lie-invariant:
\begin{equation}\label{eq-group-2}
 (x,y,t)\circ (\x,\y,\t)= (x_{1,n}+ \x_{1,n},\, x_{n+1}+\x_{n+1} +2 \langle x_{1,n}, \y_{1,n} \rangle - \t |x_{1,n}|^2,\, y_{1,n}+\y_{1,n}-\t x_{1,n},\, t+\t),
\end{equation}
and the dilation $\dd_\l$:
\begin{equation}\label{def.dil.k.lift}
 \tilde \d_\lambda (x,y,t)
   = \big(\lambda x_{1,n}, \lambda^{4}  x_{n+1}, \lambda^3 y_{1,n}, \lambda^2 t \big).
\end{equation}
Therefore, the lifted operator $\tilde \L$ satisfies \Ass{H} and \Ass{L}.In the sequel we will consider admissible paths in the following form
\begin{equation*}
    \tilde\g'(s)= \sum_{j=1}^n \omega_j \widetilde Y_j(\tilde\g(s)) + \widetilde{Z}(\tilde\g(s)), \quad \ s \in [0,\tilde \t],
\end{equation*}
for some constant vector $\omega = (\omega_1, \dots, \omega_n), \tilde\g (0) = (x,y,t)$. Its explicit expression is
\begin{equation} \label{e-lemma43}
	\tilde \g(s) = \left( x_{1,n} + s \omega, x_{n+1} + s |x_{1,n}|^2 + s^2 \langle x_{1,n}, \omega \rangle + \frac{s^3}{3}|\omega|^2, y + s  x_{1,n} + \frac{s^2}{2} \omega, t-s \right).
\end{equation}

\bigskip

In order to prove an invariant Harnack inequality for the non-negative solutions to $\tilde \L v=0$, we describe the sets $\O_{z_0}$ and $\A_{z_0}$ in the case when $z_0$ is the origin and 
\begin{equation} \label{eq-O-K}
 \O = \Big\{ (x,y,t) \in \R^{2n+2} \mid |x_{1,n}| < 1, -1 < x_{n+1} < 1, |y| < 1, - 1 < t < 1 \Big\}.
\end{equation}

\begin{lemma} \label{pr-Oz} Let $\O$ be the open set defined in \eqref{eq-O-K}, and let $z_0=(0,0,0)$. Then
\begin{equation} \label{e-4prop}
    \A_{z_0} = \left\{ (x,y,t) \in \O \mid  0 \le x_{n+1} \le -t, |y|^2  \le -t \, x_{n+1} \right\},
\end{equation}
and $\O_{z_0}=\A_{z_0}$.
\end{lemma}
\noindent{\it Proof.} \ In order to prove \eqref{e-4prop}, we consider any $\tilde \L$-admissible curve $\gamma$ in $\O$. In our setting, the components $x_{n+1}, y_{1,n}$ and $t$ of every diffusion trajectory are constant functions. Moreover, any drift trajectory $\g: [0, T] \to \O$ starting from $(\bar{x}, \bar{y}, \bar{t})$ is given by
\begin{equation} \label{eq-pp}
 \g(s) = (\bar{x}_{1,n}, \bar{x}_{n+1} + s |\bar{x}_{1,n}|^2, \bar{y} + s \bar{x}_{1,n}, \bar{t} - s).
\end{equation}
Hence, any $\tilde \L$-admissible curve $\gamma:[0,T] \to \O$ with $\gamma(0)=(0,0,0)$ is given by
\begin{equation*}
 \g(s)=\left(x_{1,n}(s), \int_{0}^{s} \sum_{k=1}^m |c_k|^2 \I_{I_k}(r)\, d r, \int_{0}^{s} \sum_{k=1}^m c_k \I_{I_k}(r) \,d r, \, - \!\int_{0}^{s} \sum_{k=1}^m \I_{I_k}(r)\, d r\right)\!,\quad s \in [0,T].
\end{equation*}
Here $I_1, \dots, I_m$ are disjoint intervals contained in $[0,T]$ and $\I_{I_k}$ denotes the characteristic function of $I_k$. The function $x_{1,n}$ is constant on every $I_k$, and any $c_k$ is a constant vector such that $|c_k| \le 1$ for $k=1, \dots, m$. As a consequence of the H\"older inequality we find 
\begin{equation*}
   \A_{z_0} \subseteq \left\{ (x,y,t) \in \O \mid 0 \le x_{n+1} \le -t, |y|^2 \le -t \, x_{n+1} \right\}.
\end{equation*}

In order to prove the opposite inclusion, we consider any point
\begin{equation*}
   (\bar{x}, \bar{y}, \bar{t}) \in \left\{ (x,y,t) \in \O \mid 0 < x_{n+1} < -t, |y|^2 < -t \, x_{n+1} \right\}, \qquad \bar y \ne 0,
\end{equation*}
and we show that there exists a $\tilde \L$-admissible curve $\gamma = \g_1+ \g_2 + \dots + \g_5$ contained in $\O$, which steers $(0,0,0)$ to $(\bar{x}, \bar{y}, \bar{t})$. To this aim, we fix a small positive $\e$, that will be specified in the sequel, and we set
\begin{equation*}
	s_\e = \frac{- \bar{t} \bar{x}_{n+1}- |\bar{y}|^2}{\bar{x}_{n+1} - 2 |\bar{y}|(1- \e) - \bar{t}(1-\e)^2}.
\end{equation*}
Note that $- \bar{t} \bar{x}_{n+1} + 2 |\bar{y}| \bar{t}(1- \e) + \bar{t}^2 (1-\e)^2 \ge \big( |\bar{y}|  + \bar{t}(1- \e) \big)^2$, so that $0 < s_\e < -\bar{t}$. We set $\tilde x_{1,n} = \frac{1-\e}{|\bar y|}\bar y$ and we choose $\g_1$ as a diffusion trajectory connecting $(0,0,0)$ to $(\tilde x_{1,n},0,0,0)$, and $\g_2:[0, s_\e] \to \R^{2n+2}$ as a drift trajectory starting from $(\tilde x_{1,n},0,0,0)$. Hence, according to \eqref{eq-pp}, we find $\g_2(s_\e) = \big( \tilde x_{1,n}, s_\e(1-\e)^2, s_\e\frac{1-\e}{|\bar y|} \bar y, - s_\e \big)$.
Then, by a diffusion trajectory $\g_3$, we connect $\g_2(s_\e)$ to the point $\left( \frac{|\bar y| - s_\e (1- \e)}{(-\bar t - s_\e)| \bar y|} \bar y, s_\e(1-\e)^2, s_\e\frac{1-\e}{|\bar y|} \bar y, - s_\e \right)$. We next consider a drift path $\g_4:[0, - \bar{t} - s_\e] \to \R^{2n+2}$ which, by \eqref{eq-pp}, and by our choice of $s_\e$, steers the end point of $\g_3$ to $\left(\frac{|\bar y| - s_\e (1- \e)}{(-\bar t - s_\e)| \bar y|} \bar y, \bar{x}_{n+1}, \bar{y}, \bar{t} \right)$. Finally, we can find a diffusion path $\g_5$ connecting $\g_4(- \bar{t} - s_\e)$ to $(\bar{x}, \bar{y}, \bar{t})$.

Clearly, $\gamma = \g_1+ \g_2 + \dots + \g_5$ is a $\tilde \L$-admissible curve of $\R^{2n+2}$ connecting $(0,0,0)$ to $(\bar{x}, \bar{y}, \bar{t})$. Next we prove that, for sufficiently small $\e$, the trajectory $\g$ is contained in $\O$. To this aim, as the set $\O$ is convex and the paths $\g_1, \g_2, \dots, \g_5$ are segments, we only need to show that the end-points of $\g_1, \g_2, \g_3, \g_4$ belong to $\O$. The inequalities $-1 < \frac{|\bar{y}| - s_\e(1-\e)}{-\bar{t} - s_\e} < 1$ directly follow from the definition of $s_\e$, for sufficiently small positive $\e$. The other inequalities are a plain consequence of the fact that $0 < s_\e < -\bar{t} < 1$, as previously noticed. Since $\A_{z_0}$ is the closure of the set of the points that can be reached by a $\tilde \L$-admissible path, we get
\begin{equation*}
    \left\{ (x,y,t) \in \O \mid 0 \le x_{n+1} \le -t, |y|^2 \le -t \, x_{n+1} \right\} \subseteq \A_{z_0}.
\end{equation*}
This concludes the proof of \eqref{e-4prop}.

To complete the proof, by Proposition \ref{pr.cond.suff.Aass} it is sufficient to find a non-negative solution $v$ of $\tilde \L v = 0$, such that $v \equiv 0$ in $\A_{z_0}$, and $v>0$ in $\O \setminus \A_{z_0}$. Let $\varphi$ be any function in $C(\partial \O)$, such that
$\varphi \equiv 0$ in $\partial \O \cap \A_{z_0}$ and $\varphi > 0$ in
$\partial \O \setminus \A_{z_0}$. Then the Perron-Wiener-Brelot solution $v := H_\varphi^\O$ of the following Cauchy-Dirichlet problem
\begin{equation*}
 \begin{cases}
  \tilde \L v = 0 & \text{in } \O\\
  v = \varphi & \text{in } \partial \O
 \end{cases}
\end{equation*}
is non-negative. Next we prove that $v>0$ in $\O\setminus \A_{z_0}$. By contradiction, let $(x,y,t)\in \O \setminus \A_{z_0}$ be such that $v(x,y,t)=0$. Then $(x,y,t)$ is a minimum for $v$, so that from Bony's minimum principle \cite[Th\'eor\`eme 3.2]{Bony} it follows that
$v(\tilde x_{1,n},x_{n+1},y,t) = \varphi(\tilde x_{1,n},x_{n+1},y,t)=0$, for every $\tilde x_{1,n} \in \partial (]-1,1[^n)$. Since every point $(\tilde x_{1,n},x_{n+1},y,t)$ is regular for the Dirichlet problem, and belongs to $\partial \O \setminus \A_{z_0}$, we find a contradiction with our assumption on $\varphi$. Suppose now that there exists $(x,y,t)\in \A_{z_0}$ such that $v(x,y,t) > 0$. Since every point of the set $\partial \O \cap \A_{z_0}$ is $\tilde \L$-regular, $v$ is continuous in ${\A_{z_0}}$. Hence there exists a $(\bar x, \bar y, \bar t) \in {\A_{z_0}}$ such that $v(\bar x, \bar y, \bar t) = \max_{{\A_{z_0}}} v > 0$. By Bony's minimum principle we have $v(\tilde x_{1,n}, \bar x_{n+1}, \bar y, \bar t) = \varphi(\tilde x_{1,n}, \bar x_{n+1}, \bar y, \bar t) >0$, for any $\tilde x_{1,n} \in \partial (]-1,1[^n)$, and this fact contradicts our assumption on $\varphi$. \hfill $\square$

\medskip

Next we introduce some notations to state a Harnack inequality which is invariant with respect to the group law $\circ$ defined in \eqref{eq-group-2} and the dilation $\tilde \delta_r$ introduced in \eqref{def.dil.k.lift}. Consider the box $Q_r = ]-r,r[^n \times ]-r^4,r^4[ \times ]-r^3,r^3[^n \times ]-r^2,0]$, and note that $Q_r = \tilde \delta_r Q_1$. For every compact set ${\K} \subseteq Q_1$, for any positive $r$ and for any $z_0 \in \R^{2n+2}$ we denote by
\begin{equation} \label{QKr}
    Q_r(z_0) = z_0 \circ \tilde\delta_r Q_1 = \big\{ z_0 \circ \tilde\delta_r \z \mid \z \in Q_1 \big\}, \qquad {\K}_r(z_0) = z_0 \circ \tilde\delta_r {\K}.
\end{equation}

\begin{corollary} \label{th-Harnack2} For every compact set ${\K} \subseteq \left\{ (x,y,t) \in Q_1 \mid 0 < x_{n+1} < -t, |y|^2 < -t x_{n+1} \right\} $, $r>0$ and $z_0\in\R^{2n+2}$ there exists a positive constant $C_{\K}$, depending only on $\tilde \L$ and ${\K}$, such that
\begin{equation*}
    \sup_{{\K}_r(z_0)}v\le C_{\K} \,v(z_0),
\end{equation*}
for every non-negative solution $v$ of $\tilde \L v=0$ on any open set containing $\overline Q_r(z_0)$.
\end{corollary}

\medskip \noindent{\it Proof.} \
Consider the function $w(z)=v\big(z_0\circ \tilde \delta_r z \big)$. By the invariance with respect to $\tilde \delta_r$ and $\circ$, we have $\tilde \L w =0$ in $Q_1$. Aiming to apply Theorem \ref{tt-1}, we consider the open set $\O$ defined in \eqref{eq-O-K}, and we note that $\O \cap \big\{ t< 0 \big\} \subset Q_1$. Then $w$ is defined as a continuous function on $\partial \O \cap \big\{ t< 0 \big\}$. We extend $w$ to a continuous function on $\partial \O$, and we solve the boundary value problem $\tilde \L \widetilde w =0$ in $\widetilde Q_1$, with $\widetilde w = w$ in $\partial \O$. Then we apply Theorem \ref{tt-1} and Lemma \ref{pr-Oz}, and we get $\sup_{{\K}} \widetilde w\le C_{\K} \, \widetilde w(0,0,0)$. By the comparison principle we have $\widetilde w = w$ in $\O \cap \big\{ t \le 0 \big\}$, then the claim plainly follows from the inclusion ${\K} \subset \O \cap \big\{ t < 0 \big\}$. 
\hfill $\square$

\medskip

We are now ready to build a Harnack chain for \eqref{K-PROC-k=2} by using the following set
\begin{equation} \label{eq-kappa}
    {\K} = \Big\{ (x,y,t) \in \R^{2n+2} \mid |x_{1,n}| \le \tfrac{1}{2}, \tfrac{1}{32} \le x_{n+1} \le \tfrac{1}{4}, |y| \le \tfrac{1}{8},  t = -\tfrac{1}{2} \Big\}
\end{equation}
which is a compact subset of $\left\{ (x,y,t) \in Q_1 \mid 0 < x_{n+1} < -t, |y|^2 < -t x_{n+1} \right\}$. Before doing that for $k=2$ only, we extend the above procedure to equations \eqref{K-PROC} and \eqref{K-PROC_2} for $k>2$.


\medskip

We next show that, in both cases \eqref{K-PROC} and \eqref{K-PROC_2}, $\L$ can be lifted to a suitable operator $\tilde \L$ in the form \eqref{e1} satisfying \Ass{H} and \Ass{L}. We introduce a new variable $y \in \R^{(k-1)n}$, that will be denoted as follows $y= (y_{1}, y_{2}, \dots ,y_{(k-1)})$, with $y_j = (y_{j1}, \dots, y_{jn}) \in \R^n$ for $j \in \leftB 1, k-1 \rightB$. We then define the  lifted vector fields on $\R^{kn+2}$:
\begin{equation} \label{e-LL.k.1}
 \tilde Y_i = Y_i = \partial_{x_i},\quad i \in \leftB 1,n \rightB, \qquad \widetilde{Z}  = Z + \sum_{i=1}^{k-1}\sum_{j=1}^n x_j^i \partial_{y_{i j}};
\end{equation}
where $Z = |x_{1,n}|^k \p_{x_{n+1}} - \p_t$ for \eqref{K-PROC}, and $Z = \sum_{j=1}^n x_{j}^k \p_{x_{n+1}} - \p_t$ for \eqref{K-PROC_2}. If we denote $v(x,y,t) = u(x,t)$ for any $u \in C^\infty(\R^{n+2})$, we have
\begin{equation*}
  \tilde Y_i v (x,y,t)  = Y_i u(x,t), \quad \forall i \in \leftB 1,n \rightB, \qquad \tilde Z v (x,y,t) = Z u(x,t).
\end{equation*}
Then, setting $\tilde \L  = \tfrac12 \sum_{i=1}^n \tilde Y_i^2\! +\tilde Z$, we plainly find $\tilde \L v (x,y,t) = \L u(x,t)$.

Since $\text{dim}\big( \text{Lie} \{\tilde Y_1, \dots, \tilde Y_n,  \tilde Z \} \big) = kn+2$ and $\text{rank}\big(\text{Lie}\{\tilde Y_1, \dots, \tilde Y_n,  \tilde Z \}(x,y,t)\big) = kn+2$ at every point $(x,y,t) \in \R^{kn+2}$, 
Theorem 1.1 in \cite{BonfiglioliLanconelli1} yields the existence of a homogeneous Lie group $\GG = \big( \R^{kn+2}, \circ, ({\dd_\l})_{\l >0} \big)$ such that $\tilde \L$ is Lie-invariant on $\GG$. Therefore, the lifted operators $\tilde \L$ satisfy \Ass{H} and \Ass{L}. The dilation $\dd_\l$ acts as follows:
\begin{equation} \label{eq-ddl}
		\dd_\l(x,y,t) = \left( \l x_{1,n}, \l^{k+2} x_{n+1}, \l^3 y_1, \dots, \l^{k+1} y_{k-1}, \l^2 t\right),
\end{equation}
for every $(x,y,t) \in \R^{kn + 2}$, and $\l >0$.
We next aim to apply Theorem \ref{tt-1} in order to prove a Harnack inequality on the lifted space $\R^{kn+2}$. 
For any $\omega \in L^2([-T,T], \R^n)$ for every $(x, y,t) \in \R^{kn+2}$ and $T>0$, we denote by $\tilde \gamma: [-T,T] \to \R^{kn+2}$ the solution of the Cauchy problem
\begin{equation} \label{e-gamma-2+}
\begin{cases}
   \tilde \g'(s) & \!\!\!\!\! = \sum_{j=1}^n \omega_j(s) \widetilde{Y}_j(\tilde\g(s)) + \widetilde{Z}(\tilde\g(s)),
   \quad \ s \in [-T,T], \\
   \tilde \g(0) & \!\!\!\!\!\! = (x,y,t).
\end{cases}
\end{equation}
In order to simplify the notation, in the sequel we will denote the solution of \eqref{e-gamma-2+} as
\begin{equation} \label{harnack-loc-66}
   \g(s) = \left( x_{1,n}(s), x_{n+1}(s), y(s), t(s) \right), \qquad s \in [-T,T].
\end{equation}
Note that $t(s) = t-s$ for every $s \in [-T,T]$, so that $t(T) = t - T$.

The composition law ``$\circ$'' of $G$ is related to \eqref{e-gamma-2+} as follows: if $\left(\bar x, \bar y, \bar t\right) = \tilde \g(T)$ is the end point of the path $\tilde \g$ defined by \eqref{e-gamma-2+} with $\tilde \g (0) = (0,0,0)$ and $\left(\tilde x, \tilde y, \tilde t\right) = \tilde \g(T)$ is the end point of the path $\tilde \g$ defined by \eqref{e-gamma-2+} with $\tilde \g (0) = \left(\xi, \eta, \tau\right)$, then
\begin{equation} \label{eq-ccirc}
		\left(\tilde x, \tilde y, \tilde t\right) = \left(\xi, \eta, \tau\right) \circ \left(\bar x, \bar y, \bar t \right),
\end{equation}
with $\bar t = - T$ (see for instance Corollary 1.2.24 in \cite{LibroBLU}). The above identity also holds when computing $\tilde \gamma$ at $s= -T$. In particular, if we choose any $\bar \omega \in \R^n$, and $\bar t >0$, we let $T = \bar t, \omega(s) = \bar \omega$ for any $s \in [-\bar t,\bar t]$, we find
\begin{equation} \label{eq-ccirc-2}
\begin{split}
		& \bar x = \left( - \bar t \bar \omega, - \tfrac{\bar t^{k+1}}{k+1} | \bar \omega|^k \right), \qquad \tilde t = \tau + \bar t, \\
		& \tilde x = \bigg( \xi_{1,n} - \bar t \bar \omega, \xi_{n+1} - \int_0^{\bar t} \left|\x_{1,n} -s \bar \omega \right|^k d s\bigg).
\end{split}
\end{equation}
According with Remark \ref{REM-GAMMA}, the fundamental solution $\tilde \Gamma$ of $\tilde \L$ exists and is invariant with respect to the group operations \eqref{eq-ccirc} and \eqref{eq-ddl}. Then function
\begin{equation} \label{eq-FS-unlift}
		\Gamma(x,t, \x, \tau) = \int_{\R^{(k-1)n}}\tilde \Gamma(x,y,t, \x, 0, \tau) d y
\end{equation}
is a fundamental solution to $\L$ and gets from $\tilde \Gamma$ the following invariance properties:
\begin{equation} \label{eq-FS-prop}
\begin{split}
        & \Gamma(\tilde x, \tilde t, \xi, \tau) = \Gamma(\bar x, \bar t, 0, 0), \\
        & \Gamma\big(\lambda \bar x_{1,n}, \lambda^{k+2} \bar x_{n+1},\lambda^2 \bar t,0_{1,n+1},0 \big) = \frac{1}{\lambda^{n+k+2}}
        \Gamma(\bar x,\bar t, 0_{1,n+1} ,0), 
\end{split}
\end{equation}
for every $(\xi, \tau) (\bar  x, \bar  t) \in \R^{n+2}, \lambda >0$, where $(\tilde x, \tilde t)$ is defined in \eqref{eq-ccirc}. In the following remark we summarize the above properties when $(\bar  x, \bar  t)$ has the form \eqref{eq-ccirc-2}.

\begin{remark} \label{rem-Gamma-omo}
For every $(\xi, \tau) \in \R^{n+2}, t >0$ and for any constant vector $\omega \in \R^n,$ we have
\begin{equation*}
    \Gamma\bigg(\xi_{1,n} - \sqrt{t} \omega, \xi_{n+1} - \int_0^{t} \left|\x_{1,n} - \frac{s}{\sqrt{t}} \omega \right|^k d s, \tau + t, \xi, \tau \bigg) = \frac{1}{t^{\frac{n+k}{2} +1}}\Gamma\Big(- \omega, - \tfrac{|\omega|^k}{k+1}, 0_{1,n+1}, 0\Big).
\end{equation*}
\end{remark}

We next focus on the attainable set $\A_{z_0}$ of the unit cylinder
\begin{equation} \label{eq-O-K+}
 \O = \Big\{ (x,y,t) \in \R^{kn+2} \mid |x_{1,n}| < 1, -1 < x_{n+1} < 1, |y| < 1, - 1 < t < 1 \Big\},
\end{equation}
with respect to the point $z_0 = (0,0,0)$. Here $|x_{1,n}|$ and $|y|$ denote, respectively, the Euclidean norm of the vectors $x_{1,n} \in \R^n$ and $y \in \R^{(k-1)n}$.

Unlike the case $k=2$, as $k>2$ we are not able to give a complete characterization of the sets $\A_{z_0}$ and $\O_{z_0}$ as we did in Lemma \ref{pr-Oz}. We will consider instead the differential of the \emph{end point map} related to \eqref{e-gamma-2+} to find some interior points of $\A_{z_0}$.
With obvious meaning of the notations, we set $\big(x(T), y(T), t(T)\big) = \tilde \g(T)$, we note that $t(T) = t-T$, and we define
\begin{equation} \label{epm}
 E : L^2([0, T]) \to \R^{kn+1}, \qquad E(\omega) = E(\omega, x,y,t, T) := \big(x(T), y(T)\big).
\end{equation}
We refer to the classical literature (see \emph{e.g.} \cite[Theorem 3.2.6]{BressanPiccoli}) for the differentiability properties of $E$.
We next show that the differential $D E(\omega)$ of $E$, computed at some given $\omega \in L^2([0,T])$ is surjective. Hence $E(\omega)$ is an interior point of $\A_{z_0}$, so that we can apply Theorem \ref{tt-1}.

\begin{lemma} \label{lemma+}
Let $\bar w$ be any given vector of $\R^n$ such that $\bar w_j \ne 0$ for every $j \in \leftB 1, n \rightB$. Consider the solution $\tilde \gamma$ to the problem \eqref{e-gamma-2+}, with $\omega \equiv \bar w$. Then $D E(\omega)$ is surjective.
\end{lemma}

\medskip \noindent {\it Proof.}\ By the invariance of the vector fields $\tilde Y_i, i \in \leftB 1,n \rightB$, and $\tilde Z$ with respect to the homogeneous Lie group $\GG$, is not restrictive to assume $(x,y,t) = (0,0,0)$ and $T=1$. To prove our claim, we compute
$$
  D E(\omega)\tilde \omega = \lim_{h\to 0} \frac{1}{h} \big( E(\omega + h \tilde \omega) - E(\omega) \big),
$$
where 
\begin{eqnarray} \label{eq-epmap}
			\tilde \omega (s) &=& \frac{1}{b-a} v \quad \text{for} \quad s \in [a,b], \ a, b \in [0,1],\ a<b,\ v \text{ is any vector of }\R^n,\nonumber\\
			\tilde \omega (s) &=&0 \quad \text{for} \quad s \not \in [a,b].
\end{eqnarray}
In the sequel, we denote by $\tilde \g^h(s) = \big(x^h(s), y^h(s), t^h(s) \big)$ the solution of \eqref{e-gamma-2+} relevant to $\omega + h \tilde \omega$.
Clearly, $t^h(s) = -s$, and $x^h(1)= \bar w + h v$, so that
\begin{equation} \label{exh}
 \lim_{h\to 0} \frac{1}{h} \left( x_{1,n}^h(1) - x_{1,n}(1) \right) = v.
\end{equation}
We next show that, for every $j \in \leftB 1, n \rightB$ and $i \in \leftB 1, k-1 \rightB$, we have
\begin{equation} \label{eyh}
 \lim_{h\to 0} \frac{y_{ij}^h(1) - y_{ij}(1)}{h} = \left( \frac{i}{i+1} \frac{b^{i+1}- a^{i+1}}{b-a} - a \frac{b^{i}- a^{i}}{b-a} + 1-b^i \right) \bar w_j^{i-1} v_j.
\end{equation}
Indeed, we have
\begin{gather*}
 \begin{split}
  y_{ij}^h(1) & = \int_0^{a} (t \bar w_j)^i d t + \int_{a}^{b} \left( t \bar w_j +h \frac{t-a}{b-a} v_j \right)^{i} d t + \int_{b}^{1} \left( t \bar w_j +h v_j \right)^{i} d t \\
  & = \int_0^{a} (t \bar w_j)^i d t + \int_{a}^{b} \left( t \bar w_j \right)^{i} d t + \int_{b}^{1} \left( t \bar w_j \right)^{i} d t \\
  & + i h \bar w_j^{i-1}  v_j \left( \int_{a}^{b} t^{i-1} \frac{t-a}{b-a} d t + \int_{b}^{1} t^{i-1} d t \right) + o(h) , \quad \text{as} \ h\to 0, \\
  & = \,y_{ij} (1) + \left( \frac{i}{i+1} \frac{b^{i+1}- a^{i+1}}{b-a} - a \frac{b^{i}- a^{i}}{b-a} + 1-b^i \right) \bar w_j^{i-1} v_j h+  o(h), \quad \text{as $h\to 0$,}
 \end{split}
\end{gather*}
where $o(h)$ vanishes as $h$ goes to zero. This proves \eqref{eyh}. Analogously,
\begin{equation} \label{exn1}
 \lim_{h\to 0} \frac{x_{n+1}^h(1) - x_{n+1}(1)}{h} = \left( \frac{k}{k+1} \frac{b^{k+1}- a^{k+1}}{b-a} - a \frac{b^{k}- a^{k}}{b-a} + 1-b^{k} \right) |\bar w |^{k-2} \langle \bar w, v \rangle,
\end{equation}
when considering system \eqref{PROC}, and
\begin{equation} \label{exn1_2}
 \lim_{h\to 0} \frac{x_{n+1}^h(1) - x_{n+1}(1)}{h} = \left( \frac{k}{k+1} \frac{b^{k+1}- a^{k+1}}{b-a} - a \frac{b^{k}- a^{k}}{b-a} + 1-b^{k} \right) \sum_{j=1}^n \bar w_j^{k-1} v_j,
\end{equation}
in the case of \eqref{PROC_2}. Note that for all $i\in \leftB 1,k\rightB $ one has
\begin{equation} \label{ec_iab}
\frac{i}{i+1} \frac{b^{i+1}- a^{i+1}}{b-a} - a \frac{b^{i}- a^{i}}{b-a} = O(b-a)
, \qquad \text{as} \quad b-a \to 0,
\end{equation}
for any $i \in \leftB 1, k \rightB$. Then, from \eqref{exh}, \eqref{eyh}, \eqref{exn1}, in the case \eqref{PROC}, it follows that
\begin{equation} \label{DEomegatilde}
\begin{split}
			D E(\omega)\tilde \omega = & \bigg( v, \big(1-b^{k}\big) |\bar w |^{k-2} \langle \bar w, v \rangle, (1-b) v, \big(1-b^2\big)\bar w_1 v_1, \dots,  \big(1-b^2\big) \bar w_n v_n, \\
			& \dots, \big(1-b^{k-1}\big) \bar w_1^{k-2} v_1, \dots, \big(1-b^{k-1}\big) \bar w_n^{k-2} v_n \bigg) + O(b-a),
\end{split}
\end{equation}
as $b-a \to 0$. We next choose $b_0, \dots, b_k \in ]0,1]$ such that $b_i \ne b_m$ if $i \ne m$ and we let $v$ be any unit vector $e_j$ of the canonical basis of $\R^n$.
Then the $j-th$, the $n+j+1-th$ $\dots$, the $(k-1)n+j+1-th$ components of $D E(\omega)\tilde \omega$ are
$$
  \left(1, 1-b_i, (1-b_i^2) w_j, \dots,  \big(1-b_i^{k-1}\big) \bar w_j^{k-2} \right),
$$
while the $n+1-th$ component is $\big(1-b_i^{k}\big) |\bar w |^{k-2} \bar w_j$. By our assumption, $\bar w_j \neq 0$, and the
following $(k+1)\times(k+1)$ matrix
\begin{equation*}
M(b_0,b_1,\dots, b_k)=
\left(
  \begin{array}{cccccc}
    1 & 1-b_0 & 1-b_0^2 & \dots & 1-b_0^k \\
    1 & 1-b_1 & 1-b_1^2 & \dots & 1-b_1^k \\
    \vdots & \vdots & \vdots & \ddots & \vdots \\
    1 & 1-b_k & 1-b_k^2 & \dots &  1-b_k^k\\
  \end{array}
\right).
\end{equation*}
is non singular, since
\begin{equation*}
    \det \,M(b_0,b_1,\dots, b_k)= (-1)^k \prod_{i\neq m}\big(b_i-b_m\big)\neq 0,
\end{equation*}
because of our choice of the $b_i$'s. Thus, if we choose $v = e_j$ and each $a_i$ sufficiently close to $b_i$, then \eqref{DEomegatilde} restores $k+1$ linearly independent vectors. In conclusion, it is possible to find $v_1, \dots, v_n, b_0, \dots, b_k, a_0, \dots, a_k$, such that the vectors $D E(\omega)\tilde \omega$ defined by using $v_j, a_i, b_i$ in \eqref{eq-epmap}, span $\R^{kn + 1}$.
This proves our claim for system \eqref{PROC}. The proof in the case \eqref{PROC_2} is analogous, we only need to replace \eqref {exn1} by \eqref {exn1_2}. We omit the details. \hfill $\square$

\medskip

We next obtain, as a corollary, a Harnack inequality which is invariant with respect to the Lie group $\GG = \big( \R^{kn+2}, \circ, ({\dd_\l})_{\l >0} \big)$.
For every compact subset ${\K}$ of the unit cylinder $\O$ defined in \eqref{eq-O-K+}, any positive $r$ and any $z_0 = (x_0, y_0, t_0) \in \R^{kn+2}$ we set
\begin{equation} \label{QKrtilde}
    \O_r(z_0) = z_0 \circ \tilde\delta_r \O = \big\{ z_0 \circ \tilde\delta_r \z \mid \z \in \O \big\}, \qquad {\K}_r(z_0) = z_0 \circ \tilde\delta_r {\K}.
\end{equation}

We recall the definition of the end-point map $E$ introduced in \eqref{epm}, and we define
the following sets
\begin{equation} \label{Ktilde}
\begin{split}
    \tilde I & := \Big\{ \omega :\left[0,  \tfrac{1}{2}\right] \to  \R^n \mid \omega(s) \equiv \tilde \omega \ \text{with} \ \tfrac{1}{4} \le | \tilde \omega_j | \le 1, j \in \leftB 1, n \rightB \Big\}, \\
    \widehat I & := \Big\{ \omega :\left[0,  \tfrac{1}{2}\right] \to  \R^n \mid \omega_j(s) = \hat c \ \text{for any} \ s \in \left[0, \tfrac{1}{8}\right] \cup \left[\tfrac{3}{8}, \tfrac{1}{2} \right], \\
    & \qquad \omega_j(s) = - \hat c \ \text{for any} \ s \in \left[\tfrac{1}{8}, \tfrac{3}{8} \right], j \in \leftB 1, n \rightB \
    \text{with} \ 1 \le \hat c \le \sqrt 2,  \Big\}, \\
    \tilde {\K} & := \Big\{ \left(E(\omega, 0,0,0,  \tfrac{1}{2} ),  \tfrac{1}{2} \right) \mid \omega \in \tilde I \Big\}\qquad \widehat {\K}
     := \Big\{ \left(E(\omega, 0,0,0,  \tfrac{1}{2} ),  \tfrac{1}{2} \right) \mid \omega \in \widehat I \Big\}.
\end{split}
\end{equation}
We remark that $\tilde {\K}$ and $\widehat {\K}$ are compact subset of the unit cylinder defined in \eqref{eq-O-K+}, being $E$ a continuous map.
We also note that, if we denote by $\bar z = (\bar x, \bar y, \bar t)$ the point $\left(E(\omega, 0,0,0,  \tfrac{1}{2} ),  \tfrac{1}{2} \right)$ with $\omega$ as described in $\widehat I$, we have $\bar x_{1,n} = 0_{1,n}, \bar x_{n+1} = \hat c^k a_{n+1,k}$, with:
\begin{equation}
\label{Qbarz}
   \begin{array}{ll}
     a_{n+1,k} = \tfrac{4 n^{k/2}}{k+1}8^{-(k+1)}  & \text{for \eqref{K-PROC}},\\
     a_{n+1,k} = \tfrac{4 n}{k+1}8^{-(k+1)}  & \text{for \eqref{K-PROC_2} and} \ k \ \text{even},  \\
     a_{n+1,k} = 0  & \text{for \eqref{K-PROC_2} and} \ k \ \text{odd}. 
    \end{array}
 \end{equation}

\begin{proposition} \label{t-Harnack.ineq+++}
Let $\O\subseteq \R^{kn+2}$ be the unit cylinder defined in \eqref{eq-O-K+}, let $r>0$ and let $z_0=(x_0,y_0,t_0) \in \rnn$.
If $\tilde \L$ is the lifted operator of $\L$ in \eqref{K-PROC} or \eqref{K-PROC_2}, then there exists a positive constant $\tilde c$ such that the sets $\tilde {\K}$  and $\widehat {\K}$ defined in \eqref{Ktilde} are compact subsets of \, {\rm Int}$\left(\A_{(0,0,0)}\right)$. Moreover there exist two positive constants $ C_{\tilde \K},  C_{\widehat \K}$, only depending on $\O$ and on $\tilde \L$, such that
\begin{equation*}
    \sup_{\tilde {\K}_r(z_0)} \tilde u \le  C_{\tilde \K} \, \tilde u(z_0), \qquad \sup_{\widehat {\K}_r(z_0)} \tilde u \le  C_{\widehat \K} \, \tilde u(z_0),
\end{equation*}
for every positive solution $\tilde u$ of $\tilde L \tilde u = 0$ in $\O_r(z_0)$.

\end{proposition}

\medskip \noindent {\it Proof.}\ By the invariance of $\tilde \L$ with respect to the homogeneous Lie group $\GG$, it is not restrictive to assume $(x_0,y_0,t_0) = (0,0,0)$ and $r=1$.
By Lemma \ref{lemma+} $E(\omega)$ is an interior point of $\A_{(0,0,0)}$ for any $\omega \in \tilde I$, and for any $\omega \in \widehat I$. The conclusion follows from Theorem \ref{tt-1}. \hfill $\square$

 \setcounter{equation}{0}
\setcounter{theorem}{0}

\section{Harnack chains and lower bounds }
\label{CHAINS_SECTION}

In this section we build Harnack chains and we prove asymptotic lower bounds for positive solutions to $\L u = 0 $. In the first Lemma \ref{th-Harnack2-path} we capture paths that give Gaussian lower bounds as $|x_{1,n} - \x_{1,n}|^2 \ge K (t-\t)$, for suitably big $K$ and asymptotic bounds for points $x, \x$ with $|x_{n+1} - \x_{n+1}|$ suitably big  with respect to $(t - \t)^{1+k/2}$.

Lemma \ref{th-Harnack2-eps} applies when $|x_{n+1} - \x_{n+1}|$ is small with respect to $(t - \t)^{1+k/2}$.
Moreover, in this lemma, we consider points with non degenerate components set to zero.
In such case, if we denote $(\bar x,\bar y, \bar  t) = (\tilde x, \tilde y, \tilde t) \circ (\widehat x,\widehat y,\widehat t)$, then
\begin{equation} \label{eq-ccl}
\tilde x_{1, n}= \widehat x_{1, n}=0_{1,n} \quad \Longrightarrow	\quad	\bar  x_{ n+1} = \tilde x_{n+1} +  \widehat x_{n+1},
\end{equation}
that is the group law $ \circ$ is additive w.r.t. the $(n+1)$\emph{th} component. In some sense this allows to move in the direction of the vector field
$\underbrace{[\partial_{x_1},[\partial_{x_1},\cdots,[\partial_{x_1},Z]\cdots]]}_{k\ {\rm times}}=k!\partial_{x_{n+1}}. $

The proof of the two lemmas is based on the lifting procedure introduced in Section 5 and on the construction of a finite sequence of cylinders contained in the lifted domain of the solution. Specifically, we find points along the trajectory of the integral path introduced in \eqref{e-gamma-2+}. The bounds depend on the length of the Harnack chain that in turns depends on the slope of the trajectory. Asymptotic lower bounds are proved in Propositions \ref{th-Harnack2-6} and \ref{th-Harnack2-7}.

\begin{lemma} \label{th-Harnack2-path} Let $\L$ be the operator defined in \eqref{K-PROC} or in \eqref{K-PROC_2},
let $T_1, \t,t, T_2$ be such that $T_1 < \t < t < T_2$ and $t-\t \le \t - T_1$. Let $\g: [0,t-\t] \to \R^{n+1}\times \,]T_1,T_2[$ be a path satisfying
\begin{equation}\label{eq-admiss-e6}
    \g'(s)= \sum_{j=1}^n \omega_j  Y_j(\g(s)) + {Z}(\g(s)), \qquad \gamma(0) = (x,t), \quad \gamma(t-\t) = (\x,\t),
\end{equation}
for some constant vector of $\R^n$ such that
$$\max_{j \in \leftB 1, n \rightB} |\omega_j| \le 2 \min_{i \in \leftB 1, n \rightB} |\omega_i|,\quad
\text{and that}\quad (t-\t)\max_{j\in \leftB 1,n\rightB}\omega_j^2 \ge 4 .$$
Then there exists a positive constant $C$, only depending on $\L$, such that:
\begin{equation*}
    u(\x,\t) \le \exp \left( C \big((t-\t) \max_{j\in \leftB 1,n\rightB}\omega_j^2 + 1 \big) \right) u(x,t),
\end{equation*}
for every non-negative solution $u$ to $\L u = 0$ in $\R^{n+1} \times \,]T_1,T_2[$.
\end{lemma}

\medskip \noindent{\it Proof.} \
Define the function $\widetilde u$ by setting $\widetilde u (x, y, t)=u(x, t)$ for every $(x,y,t) \in \R^{k n + 1}\times \,]T_1,T_2[$. Clearly, $\widetilde u$ is a non-negative solution to $\tilde \L \widetilde u =0$. Let $\tilde \gamma: [0,t-\t] \to \R^{k n+1}\times \,]T_1,T_2[$ be the solution of the Cauchy problem
\begin{equation*}
\begin{cases}
   \tilde \g'(s) & \!\!\!\!\! = \sum_{j=1}^n \omega_j \widetilde{Y}_j(\tilde\g(s)) + \widetilde{Z}(\tilde\g(s)),
   \quad \ s \in [0,t-\t], \\
   \tilde \g(0) & \!\!\!\!\!\! = (x,0,t),
\end{cases}
\end{equation*}
where $\omega$ is the constant vector in \eqref{eq-admiss-e6}. Note that, if $\tilde x_{1,n+1}$ and $x_{1,n+1}$ are the first $n+1$ components of $\tilde \g$ and $\g$, respectively, then we have $\tilde x_{1,n+1}(s) = x_{1,n+1}(s)$, for every $s \in [0, t - \t]$.

We next apply the Harnack inequalities stated in Proposition \ref{t-Harnack.ineq+++} to a suitable set of points $z_1, \dots, z_m$ lying on $\tilde \gamma([0, t-\t])$. We suppose $\omega_1^2 = \max_{j \in \leftB 1, n \rightB} \omega_j^2$, as it is not restrictive, and we let $m$ be the unique positive integer such that $ m-1 < \frac{(t - \t) \o_1^2}{2} \le m$. We set $\tilde s=\frac{t - \t}{m}$ and we define $z_j=\tilde \gamma(j \tilde s)$ for $j\in \leftB 1,m \rightB$. In order to apply Proposition \ref{t-Harnack.ineq+++}, we put $\tilde r  = \sqrt{2\tilde s}$, and $\tilde \omega = \frac{\tilde r}{2} \omega$. Note that $m \ge \frac{(t - \t) \o_1^2}{2} \ge 2$ and, as a consequence, we have $\tilde \omega \in \tilde I$. Thus, if we denote by $\tilde z$ the point $\left(E(\tilde \omega, 0,0,0, \tfrac{1}{2} ), \tfrac{1}{2} \right)$ defined in \eqref{Ktilde}, we have $\tilde z \in \tilde {\K}$.
Moreover, $z_j = z_{j-1} \circ \delta_{\tilde r} \, \tilde z$, thus
\begin{equation*}
    z_j \in \tilde {\K}_{\tilde r} (z_{j-1}), \quad \text{for any} \quad j\in \leftB 1, m \rightB.
\end{equation*}
Note that by our assumption, $\frac{\tilde r^2}{2} \le \frac{2}{\omega_1^2} \le t - \t \le \t - T_1$, hence, ${\cal O}_{\tilde r} (z_{j}) \subset \R^{k n + 1}\times \,]T_1,T_2[$ for every $j \in \leftB 0, m-1 \rightB$. Thus, by Proposition \ref{t-Harnack.ineq+++}, there exists a constant $C_{\tilde \K} >1$ such that $\tilde u(z_j) \le C_{\tilde \K} \, \tilde u(z_{j-1})$ for every $j\in \leftB 1, m \rightB$. In particular, being $m <  \frac{(t - \t) \o_1^2}{2} + 1, z_0 = (x,0,t), z_{m} = \tilde \gamma( t - \t)$, we find
\begin{equation*}
    \tilde u \big(\tilde \gamma(t - \t)\big) \le C_{\tilde \K}^{\frac{(t - \t) \o_1^2}{2} +1} \tilde u(x,0,t).
\end{equation*}
 Hence,
\begin{equation*}
    u \big(\x,\t\big) = u \big(\gamma( t-\t)\big) = \tilde u \big(\tilde \gamma( t-\t)\big) \le C_{\tilde \K}^{\frac{(t - \t) \o_1^2}{2} +1} \tilde u (x,0,t) = C_{\tilde \K}^{\frac{(t - \t) \o_1^2}{2} +1} u(x,t),
\end{equation*}
and our claim follows by choosing $C := \log(C_{\tilde \K})$.
\hfill $\square$

\medskip
Note that, whenever Lemma \ref{th-Harnack2-path} applies, we have $x_{1,n}(t - \t) = x_{1,n} +  (t - \t) \o \ne x_{1,n}$. Next result gives a bound  along a trajectory $\tilde \gamma$ such that $x_{1,n}(s) = 0$ for any $s \in [0, t - \t]$.

\begin{lemma} \label{th-Harnack2-eps} Let $\L$ be the operator defined in \eqref{K-PROC} or in \eqref{K-PROC_2}, with k \emph{even}, and let $a_{n+1,k}$ be as in \eqref{Qbarz}. Let $T_1, \t,t, T_2$ be such that $T_1 < \t < t < T_2$ and $t-\t \le \t - T_1$. Then there exists a positive constant $C$, only depending on $\L$, such that:
\begin{equation*}
    u \left(0_{1,n} ,x_{n+1} + \tilde  \x_{n+1}, \t \right) \le \exp\bigg( C \bigg( \frac{(t - \t)^{1 + 2/k}}{\tilde \x_{n+1}^{2/k}} + 1\bigg) \bigg) u(0_{1,n}, x_{n+1},t).
\end{equation*}
for every $(x,t) \in \R^{n+1} \times \,]T_1,T_2[, \t \in ]T_1,t[$ with $\tilde \x_{n+1} \in \left]0, {2}\, {a_{n+1,k}} (t- \t)^{1+k/2} \right[$,
for every non-negative solution $u$ to $\L u = 0$ in $\R^{n+1} \times \,]T_1,T_2[$.
\end{lemma}

\medskip \noindent{\it Proof.} \ As in the proof of Lemma \ref{th-Harnack2-path}, we consider the function $\widetilde u$ 
defined as $\widetilde u (x, y, t)=u(x, t)$ for every $(x,y,t) \in \R^{k n + 1}\times \,]T_1,T_2[$. Let $m$ be the unique positive integer such that
\begin{equation} \label{harnack-loc-6bbb}
	m- 1 < \big( 2 (t - \t) \big)^{1 + 2/k}  \left(\frac{a_{n+1,k}}{\tilde \x_{n+1}}\right)^{2/k} \le m.
\end{equation}
Let $\tilde r = \sqrt{\frac{2 (t - \t)}{m}}$ and let $\hat c = \left(\frac{\tilde \x_{n+1}}{a_{n+1,k}}\right)^{1/k}\frac{m^{1/2}}{\big( 2 (t - \t) \big)^{1/2+1/k}}$.
Note that from our assumption $0 < \tilde \x_{n+1} < {2}\, {a_{n+1,k}} (t- \t)^{1+k/2}$ it follows that $m\ge 2 $, hence $1 \le \hat c \le \sqrt 2$. Then, if we denote by $\bar z$ the point $\left(E(\omega, 0,0,0, \tfrac{1}{2} ), \tfrac{1}{2} \right)$ defined in \eqref{Ktilde}, with $\hat c$ as above, we have $\bar z \in \widehat {\K}$.

Let $z_0 = (0_{1,n}, x_{n+1},0,t)$, and let $z_{j} = z_{j-1} \circ \delta_{\tilde r} \, \bar z$, for $j \in \leftB 1, m \rightB$. By our assumption we also have ${\tilde r^2} \le \t - T_1$, then ${\cal O}_{\tilde r} (z_{j}) \subset \R^{kn + 1} \times ]T_1, T_2[$, for every $j \in \leftB 0, m \rightB$. Thus Proposition \ref{t-Harnack.ineq+++} yields $\tilde u(z_m) \le C_{\widehat \K} ^m \tilde u (z_0).$
According with \eqref{eq-ccl}, and with our choice of $m, \hat c$, and $\tilde r$, the first $n+1$ components of $z_m$ are $\left( 0_{1,n}, x_{n+1} + \tilde \x_{n+1}\right)$. Then
\begin{equation*} 
    u( 0_{1,n}, x_{n+1} + \tilde x_{n+1}, \t) = \tilde u (z_m) \le C_{\widehat \K} ^m \, \tilde u (z_0) = C_{\widehat \K} ^m \, u ( 0_{1,n}, x_{n+1} , t),
\end{equation*}
and the conclusion follows from \eqref{harnack-loc-6bbb}. 
 \hfill $\square$

\medskip

\begin{proposition}\label{th-Harnack2-6}
Let $\L$ be the operator defined in \eqref{K-PROC} and let $k$ be a positive even integer. 
Let $u: \R^{n+1} \times ]T_1, T_2[ \to \R$ be a non-negative solution to $\L u = 0$, and let $t, \t \in \R$ be such that $T_1 < \t < t < T_2$, and $t - \t \le 2(\t - T_1)$. Then there exists a positive constant $C_1$, only depending on $\L$, such that
\begin{description}
	\item[{\it i)}] for any $x, \x \in \R^{n+1}$ such that $\x_{n+1} - x_{n+1} \ge 15^k \big(  (t - \t)(|x_{1,n}|^k + |\x_{1,n}|^k) + n^{k/2} (t - \t)^{1+k/2} \big)$ we have
\begin{equation*}
\!\!\!\!    u(\x,\t) \le \exp \bigg( C_1 \bigg( \frac{|x_{1,n}- \x_{1,n}|^2}{t-\t} + \frac{\Big(\x_{n+1}-x_{n+1}-\frac{t-\t}{2^{k+4}}\big(|x_{1,n}|^k+|\x_{1,n}|^k\big)\Big)^{2/k}\!\!\!\!}{(t-\t)^{1 + 2/k}} + 1 \bigg)\!\! \bigg) u(x,t);
\end{equation*}
	\item[{\it ii)}] for any $x, \x \in \R^{n+1}$ such that $0 < \x_{n+1} - x_{n+1} \le \frac{t - \t}{2(k+1)} \left(|x_{1,n}|^k + |\x_{1,n}|^k\right)+
	\frac{n^{k/2}}{8^{k+1}(k+1)}(t - \t)^{1+k/2}$ 
	we have
\begin{equation*}
    u(\x,\t) \le \exp \bigg( C_1 \bigg( \frac{|x_{1,n}|^{k+2} + |\x_{1,n}|^{k+2}}{\x_{n+1} - x_{n+1}} + \frac{(t - \t)^{1+2/k}}{(\x_{n+1} - x_{n+1})^{2/k}} + 1 \bigg) \bigg) u(x,t).
\end{equation*}
\end{description}
\end{proposition}

\medskip \noindent{\it Proof of (i).}
We divide the proof into two steps. In the first one we find a path $\g: \left[0, \frac{t-\t}{2}\right] \to \R^{kn+ 2}$ that steers $x_{1,n}$ to $\x_{1,n}$. In the second step another path $\g_1 + \g_2$ steers the $(n+1)$\emph{th} component $x_{n+1}\left( \frac{t-\t}{2}\right)$ of $\g\left( \frac{t-\t}{2}\right)$ to $\x_{n+1}$.

\medskip \noindent
\emph{Step 1:} 
If  $\o=\frac{2}{t-\t} (\X1n- \xx1n)$ satisfies the assumptions of Lemma   \ref{th-Harnack2-path}, then one readily gets $x_{1,n}\left(\frac{t-\t}{2}\right)=\xx1n $ and  $u\left(\gamma\left(\frac{t-\t}{2}\right) \right) \le \exp(C[2\frac{|\X1n- \xx1n|^2}{t-\t}+1])u(x,t)$ and the first step is achieved.
If this is not the case, we rely on the following construction. Set $K=\max_{j\in \leftB 1,n \rightB}\frac{|\x_j-x_j|}{\sqrt {t-\t}}$ and $M:=\max\{ 3 K, \frac{3}{2}\} $. For every $j\in \leftB 1,n\rightB $, we set
$$\tilde \o_j:=\frac{4M}{\sqrt{t-\t}},\quad \widehat \o_j:=4\frac{\x_j-x_j}{t-\t}-\frac{4M}{\sqrt{t-\t}}, $$
so that
\begin{eqnarray*}
\frac {8} 3 \frac{M}{\sqrt{t-\t}} \le |\widehat \o_j|\le \frac {16} 3 \frac{M}{\sqrt{t-\t}},\quad
\frac{t-\t}{4}\tilde \o_j^2\ge 9,\quad \frac{t-\t}{4}\widehat \o_j^2\ge 4.
\end{eqnarray*}
Consider now the path $\gamma $ associated to \eqref{eq-admiss-e6} with $\o(s)=\tilde \o \I_{s\in [0,\frac{t-\t}4]}+\widehat \o\I_{s\in [\frac{t-\t}4,\frac{t-\t}2]} $ for which $x_{1,n}\left( \frac{t-\t}{2}\right)=\xx1n$. The assumptions of Lemma \ref{th-Harnack2-path} are clearly satisfied. Hence:
\begin{equation} \label{e-g/2-6}
\begin{split}
	u\left(\g\left(\tfrac{t - \t}{2}\right) \right) & \le \exp \left( C' \left( \frac{|x_{1,n}- \x_{1,n}|^2}{t-\t} + 1 \right) \right) u(x,t), \\
	\g\left(\tfrac{t - \t}{2}\right) & = \left(\x_{1,n},x_{n+1} + \int_0^{\frac{t - \t}{2}} \left| x_{1,n} + \int_0^s\o(u)du \right|^k d s, \frac{t+\t}{2} \right),
\end{split}
\end{equation}
for some positive constant $C'$. By a plain change of variable in the above integral we find
\begin{equation} \label{e-g/2-6666}
 x_{n+1}\left(\tfrac{t - \t}{2}\right) - x_{n+1} = \int_0^{\frac{t - \t}{4}} \left| x_{1,n} +s \widetilde \o\right|^k d s + \int_0^{\frac{t - \t}{4}} \left| \x_{1,n} - s \widehat \o\right|^k d s.
\end{equation}
Note that, for $s \in \left[ 0, \frac{t - \t}{4} \right]$ and $j \in \leftB 1, n \rightB$, we have
\begin{equation} \label{e-g/2-677}
\begin{split}
 | x_{j}+s \widetilde \o_{j}| & \le | x_{j} |+ \tfrac{4 s}{\sqrt{t - \t}} M \le | x_{j} | + 3 \left( {|x_{j} - \x_{j}|} + \tfrac{{\sqrt{t - \t}}}{2} \right) \le 4 | x_{j} | + 3 |\x_{j}| + 3 \tfrac{{\sqrt{t - \t}}}{2} \\ | \x_{j} - s \widehat \o_{j}| & \le | \x_{j} |+ \tfrac{4 s}{\sqrt{t - \t}} \left( M + \tfrac{|x_{j} - \x_{j}|}{\sqrt{t - \t}}\right) \le 4 | x_{j} | + 5 |\x_{j}| + 3 \tfrac{{\sqrt{t - \t}}}{2},
\end{split}
\end{equation}
thus, from the elementary inequality $(a^2 + b^2 + c^2)^{k/2} \le 3^{k/2-1}(a^k + b^k + c^k)$, we get
\begin{equation*} 
 x_{n+1}\left(\tfrac{t - \t}{2}\right) - x_{n+1} \le \tfrac{t - \t}{4}3^{k-1} \left( 2^{2k+1} |x_{1,n}|^k + \big(3^k+5^k \big) |\x_{1,n}|^k + 2 \left(\tfrac{3^{k}n^{k/2}}{2^k}\right) (t - \t)^{k/2} \right).
\end{equation*}
Recalling that $\x_{n+1} - x_{n+1} \ge 15^k \big( (t - \t) (|x_{1,n}|^k + |\x_{1,n}|^k) + n^{k/2} (t - \t)^{1+k/2} \big)$, we find
\begin{equation} \label{e-corollary43-6}
 \x_{n+1} - x_{n+1}\left(\tfrac{t - \t}{2}\right) \ge \frac{t - \t}{2}15^k \big(  |x_{1,n}|^k + |\x_{1,n}|^k + n^{k/2} (t - \t)^{k/2} \big).
\end{equation}
We next prove a similar lower bound for $x_{n+1}\left(\tfrac{t - \t}{2}\right) - x_{n+1}$. To this aim, we note that, if $|\X1n| \ge \frac{t - \t}{8}|\tilde \o|$, then $|\X1n + s \tilde \o| \ge \frac{1}{2}|\X1n|$ for every $s \in \left[0, \frac{t-\t}{16} \right]$. Analogously, if $|\X1n| \le \frac{t - \t}{8}|\tilde \o|$, then $|\X1n + s \tilde \o| \ge \frac{t - \t}{16}|\tilde \o|$ for every $s \in \left[\frac{3}{16}(t-\t), \frac{t-\t}{4} \right]$. The same remark holds if we replace  $\X1n$ and $\tilde \o$ by $\xi_{1,n}$ and $\widehat \o$, respectively.
As a consequence we find,
\begin{equation*}
 x_{n+1} - x_{n+1}\left(\tfrac{t - \t}{2}\right) \ge \frac{t - \t}{16} \left(\max\left\{ \frac{|\X1n|^k}{2^k}, \left(\frac {t-\t}{16}|\tilde \o| \right)^k
  \right\} + \max\left\{ \frac{|\x_{1,n} |^k}{2^k},  \left(\frac {t-\t}{16}|\widehat \o| \right)^k \right\} \right),
\end{equation*}
so that, in particular,
\begin{equation} \label{e-corollary43-6bis}
 x_{n+1}\left(\tfrac{t - \t}{2}\right) - x_{n+1} \ge \frac{t - \t}{2^{k+4}} \left( |\X1n|^k + |\x_{1,n}|^k \right).
\end{equation}

\medskip \noindent
\emph{Step 2:}
We next denote by $\bar \o$ the vector in $\R^n$ such that $\bar \o_j = 1$ if $\xi_j \ge 0$, $\bar \o_j = - 1$ if $\xi_j < 0$, for $ j \in \leftB 1, n \rightB$, and we fix a real parameter $b \ge \frac{4}{\sqrt{t-\t}}$, that will be specified later. 
We consider the path  $\g_1 : \left[ \frac{t-\t}{2}, \frac{3}{4} (t-\t) \right] \to \R^{n+1} \times ]T_1, T_2[$, starting from $\g\left( \frac{t-\t}{2} \right)$ and defined as in \eqref{eq-admiss-e6} with $\o = b\, \bar \o$, then the path  $\g_2 : \left[\frac{3}{4} (t-\t) , t-\t \right] \to \R^{n+1} \times ]T_1, T_2[$, starting from $\g_1 \left( \frac{3}{4} (t-\t) \right)$ and defined by setting $\o = - b\, \bar \o$. From Lemma \ref{th-Harnack2-path} it then follows
\begin{equation} \label{e-g/2-6b}
	u\left(\x_{1,n}, \varphi(b), \t \right) \le \exp \left( 2 C \left( \frac{(t-\t)b^2}{4} + 1 \right) \right) u\left(\g\left(\tfrac{t - \t}{2}\right) \right), 
\end{equation}
where
\begin{equation*} 
	\varphi(b) = x_{n+1}\left(\tfrac{t - \t}{2}\right) + 2 \int_{0}^{\frac{t-\t}{4}}\left| \x_{1,n} + s  b \, \bar \o \right|^k ds
\end{equation*}
is an increasing continuous function of $b \in \left[ \frac{4}{\sqrt{t-\t}} , + \infty \right[$. An elementary computation shows that
\begin{equation*} 
	2 \int_{0}^{\frac{t-\t}{4}}\left| \x_{1,n} + s  b\, \bar \o \right|^k ds \le 2^k \frac{t-\t}{4}|\x_{1,n}|^k + \frac{b^kn^{k/2}}{k+1}\frac{(t-\t)^{k+1}}{2^{k+2}},
\end{equation*}
then
\begin{equation*} 
	\varphi\left( \frac{4}{\sqrt{t-\t}} \right) \le x_{n+1}\left(\tfrac{t - \t}{2}\right) +2^k \frac{t-\t}{4}|\x_{1,n}|^k + \frac{2^{k-2}n^{k/2}}{{k+1}} (t-\t)^{1+k/2} < \x_{n+1},
\end{equation*}
by \eqref{e-corollary43-6}.
On the other hand we have
\begin{equation} \label{e-g/2-6c}
	\varphi(b) \ge x_{n+1}\left(\tfrac{t - \t}{2}\right) + \frac{b^k}{k+1}\frac{(t-\t)^{k+1}n^{k/2}}{2^{2k+1}} \to + \infty, \quad \text{as} \ b \to + \infty.
\end{equation}
Hence, there exists a unique $\tilde b \ge \frac{4}{\sqrt{t-\t}}$ such that $\varphi(\tilde b) = \x_{n+1}$. Moreover from \eqref{e-g/2-6c} it also follows that
\begin{equation*} 
    \tilde b \le 2(k+1)^{1/k}\left(\tfrac{2}{t-\t} \right)^{1 + 1/k} \left(\x_{n+1} - x_{n+1}\left(\tfrac{t - \t}{2}\right) \right)^{1/k}n^{-1/2},\end{equation*}
which, together with \eqref{e-corollary43-6bis}, gives
\begin{equation*} 
    \tilde b \le (k+1)^{1/k}\left(\tfrac{2}{t-\t} \right)^{1 + 1/k} \left(\x_{n+1} - x_{n+1} - \frac{t-\t}{2^{k+4}}\left( |x_{1,n}|^k + |\x_{1,n}|^k \right) \right)^{1/k}.
\end{equation*}
Eventually, equation \eqref{e-g/2-6b} yields
\begin{equation*}
	u (\xi, \tau ) \le \exp \bigg( C_1 \bigg( \frac{(\x_{n+1}-x_{n+1}-\frac{t-\t}{2^{k+4}} (|x_{1,n}|^k+|\x_{1,n}|^k))^{2/k}}{(t-\t)^{1 + 2/k}} + 1 \bigg) \bigg) u\left(\g\left(\tfrac{t - \t}{2}\right) \right).
\end{equation*}
with $C_1 = 4^{1+ 1/k}(k+1)^{2/k} C$. The above inequality, with \eqref{e-g/2-6}, proves the claim \emph{(i)}.

\medskip \noindent
\emph{Proof of (ii)}. \ We prove our claim by applying Lemma \ref{th-Harnack2-eps} in a suitable interval $[\t + t_2, t - t_1] \subsetneq [\t, t]$, and Lemma \ref{th-Harnack2-path} in the remaining intervals $[t - t_1,t]$ and $[\t, \t + t_2]$. We first suppose that $x_{1,n} \ne 0$, $\x_{1,n} \ne 0$, we set
\begin{equation*}
\begin{split}
    t_1 & = \min \left\{\frac{\max_{j\in \leftB 1, n \rightB} x_{j}^2}{4}, \frac{\x_{n+1} - x_{n+1}}{|x_{1,n}|^k}, \frac{t-\t}{3} \right\}, \\
    t_2 & = \min \left\{\frac{\max_{j\in \leftB 1, n \rightB} \x_{j}^2}{4}, \frac{\x_{n+1} - x_{n+1}}{|\x_{1,n}|^k}, \frac{t-\t}{3} \right\},
\end{split}
\end{equation*}
and we consider the paths
\begin{equation*}
\begin{split}
	\g_1(s) & = \left( \Big(1- \frac{s}{t_1} \Big) x_{1,n}, x_{n+1} + \frac{(s- t_1)^{k+1} + t_1^{k+1}}{({k+1}) t_1^k} |x_{1,n}|^k, t - s \right), \quad s \in [0,t_1] \\
	\g_2(s)&  = \left( \frac{s}{t_2} \x_{1,n}, \x_{n+1} + \frac{s^{k+1}-t_2^{k+1}}{({k+1}) t_2^k} |\x_{1,n}|^k, \t + t_2 - s \right), \qquad s \in [0, t_2].
\end{split}
\end{equation*}

We next proceed assuming that $\max_{j\in \leftB 1, n \rightB} |x_j| \le 2 \min_{i \in \leftB 1, n \rightB} |x_i|$ and $\max_{j\in \leftB 1, n \rightB} |\x_j| \le 2 \min_{i \in \leftB 1, n \rightB} |\x_i|$. In this case we apply Lemma \ref{th-Harnack2-path} in the interval $[0, t_1]$ with $\overline \omega= - \frac{1}{t_1} x_{1,n}$, so that we have $t_1 \max_{j\in \leftB 1, n \rightB} \overline \o_j^2 \ge 4$ and $\max_{j\in \leftB 1, n \rightB} |\overline \o_j| \le 2 \min_{i \in \leftB 1, n \rightB} |\overline \o_i|$. We find
\begin{equation}\label{e-g1g2-6}
\begin{split}
		u \big(\g_1(t_1) \big) & \le \exp \left( C \left( \max \left\{ \frac{|x_{1,n}|^{k+2}}{\x_{n+1} - x_{n+1}}, 3 \frac{|x_{1,n}|^2}{t-\t}, 4 \right\} + 1 \right) \right) u (x,t), \\
		u (\x,\t ) & \le \exp \left( C \left( \max \left\{ \frac{|\x_{1,n}|^{k+2}}{\x_{n+1} - x_{n+1}}, 3 \frac{|\x_{1,n}|^2}{t-\t} , 4 \right\} + 1 \right) \right) u \big(\g_2(0)\big),
\end{split}
\end{equation}
with
$$
\g_1(t_1) = \left(0, x_{n+1} + \tfrac{t_1}{k+1} |x_{1,n}|^k, t-t_1\right), \quad	\g_2(0) = \left(0, \x_{n+1} - \tfrac{t_2}{k+1} |\x_{1,n}|^k, \t + t_2\right).
$$
If $\max_{j\in \leftB 1, n \rightB} |x_j| > 2 \min_{i \in \leftB 1, n \rightB} |x_i|$, we rely on the argument used at the beginning of the proof of \emph{(i)}. Specifically, we set $M: = \max\{3\max_{j\in \leftB 1,n \rightB} \frac{|x_j|}{\sqrt {t_1}}, \frac{3}{2}\}$, and $\tilde \o_j:= 2 \frac{M}{\sqrt{t_1}}, \widehat \o_j:= \frac{2}{t_1}x_j+ 2 \frac{M}{\sqrt{t_1}}$ for every $j\in \leftB 1,n\rightB $, then we consider the path $\gamma $ associated to \eqref{eq-admiss-e6}
with $\o(s)=\tilde \o \I_{s\in [0,\frac{t_1}2]} - \widehat \o\I_{s\in [\frac{t_1}2,t_1]}$. Also in this case we get \eqref{e-g1g2-6}, with some bigger constant $C$. Aiming to simplify our exposition we omit the details of the proof.

We next conclude the proof by using Lemma \ref{th-Harnack2-eps}. We set $\tilde \x_{n+1} = \x_{n+1} - x_{n+1}- \frac{t_1}{k+1} |x_{1,n}|^k- \frac{t_2}{k+1} |\x_{1,n}|^k$, and we recall that $\x_{n+1} - x_{n+1} \le \frac{t - \t}{2(k+1)} \left(|x_{1,n}|^k + |\x_{1,n}|^k\right)+ \frac{n^{k/2}}{8^{k+1}(k+1)}(t - \t)^{1+k/2}$.
Thus
\begin{equation} \label{e-ttxx-6}
		\frac{t-\t}{3} < (t -t_1) - (\t + t_2) < t-\t, \quad
		\frac{k-1}{k+1}
        \left(\x_{n+1} - x_{n+1}\right) \le \tilde \x_{n+1} \le \frac{n^{k/2}}{8^{k+1}(k+1)}(t - \t)^{1+k/2}.
\end{equation}
Then, from Lemma \ref{th-Harnack2-eps} we get
\begin{equation} \label{e-g1-to-g2-6}
\begin{split}
    u \left(\g_2(0) \right) & \le \exp\left( C \left( \frac{(t- t_1 - t_2 - \t)^{1+2/k}}{(\tilde \x_{n+1})^{2/k}} + 1\right) \right) u \big(\g(t_1) \big)\\
    & \le \exp\left( C \left( \left(\frac{k+1}{k-1}\right)^{2/k} \frac{(t - \t)^{1+2/k}}{(\x_{n+1} - x_{n+1})^{2/k}} + 1\right) \right) u \big(\g(t_1) \big).
\end{split}
\end{equation}
From inequalities \eqref{e-g1g2-6} and \eqref{e-g1-to-g2-6} it follows that
\begin{equation*}
    u(\x,\t) \le \exp \bigg( C' \bigg( \frac{|x_{1,n}|^{k+2} + |\x_{1,n}|^{k+2}}{\x_{n+1} - x_{n+1}} + \frac{(t - \t)^{1+2/k}}{(\x_{n+1} - x_{n+1})^{2/k}} + \frac{|x_{1,n}|^2 + |\x_{1,n}|^2}{t- \t} +1 \bigg) \bigg) u(x,t),
\end{equation*}
for some positive constant $C'$ only depending on $C$ and on $k$. Note that the last term in the above expression is bounded by the first one. Indeed, the inequality
$$
  \frac{|x_{1,n}|^2 + |\x_{1,n}|^2}{t- \t} \le \frac{2}{k+2} \frac{|x_{1,n}|^{k+2} + |\x_{1,n}|^{k+2}}{(t- \t)^{1+k/2}} + \frac{k}{k+2}
$$
combined with \eqref{e-ttxx-6}, gives
$$
 \frac{|x_{1,n}|^2 + |\x_{1,n}|^2}{t- \t} \le \frac{2}{4^{k+1}(k+2)(k-1)} \frac{|x_{1,n}|^{k+2} + |\x_{1,n}|^{k+2}}{\x_{n+1} - x_{n+1}} + \frac{k}{k+2}.
$$
This concludes the proof of \emph{(ii)} when $x_{1,n} \ne 0$, and $\x_{1,n} \ne 0$.

If $x_{1,n}=0$, we simply omit the construction of $\g_1$, and we rely on $\g_2$ and on the application of Lemma \ref{th-Harnack2-eps} in the interval $[\t + t_2, t]$. Analogously, if $\x_{1,n}=0$, we avoid the construction of $\g_2$. This concludes the proof. \hfill $\square$

\medskip

\begin{proposition}\label{th-Harnack2-7}
Let $\L$ be the operator defined in \eqref{K-PROC_2} and let $k$ be a positive integer. 
Let $u: \R^{n+1} \times ]T_1, T_2[ \to \R$ be a non-negative solution to $\L u = 0$, and let $t, \t \in \R$ be such that $T_1 < \t < t < T_2$, and $t - \t \le 2(\t - T_1)$. Then there exists a positive constant $C_1$, only depending on $\L$, such that
\begin{description}
	\item[{\it i)}] if $k$ is even, then for any $x, \x \in \R^{n+1}$ such that $\x_{n+1} - x_{n+1} \ge 15^k (t - \t) \sum_{j=1}^n \left( x_j^k + \x_j^k \right) + n \left(\frac32\right)^k (t - \t)^{1+k/2}$ we have
\begin{eqnarray*}
    u(\x,\t) \!\!\!\! & \le & \!\!\!\! \exp \bigg( C_1 \bigg( \frac{|x_{1,n}- \x_{1,n}|^2}{t-\t} + \\
   & & \qquad + \frac{(\x_{n+1}-x_{n+1}-
\frac{1}{2^{k+4}}   \sum_{j=1}^n(x_j^k+\x_j^k)(t-\t))^{2/k}}{(t-\t)^{1 + 2/k}} + 1 \bigg) \bigg) u(x,t);
\end{eqnarray*}
	\item[{\it ii)}] if $k$ is even, then for any $x, \x \in \R^{n+1}$ such that $0 < \x_{n+1} - x_{n+1} \le \frac{t - \t}{2(k+1)} \sum_{j=1}^n \left( x_j^k + \x_j^k \right) +
	\frac{(t - \t)^{1+k/2}}{4^{k+1}(k+1)}$ we have
\begin{equation*}
    u(\x,\t) \le \exp \bigg( C_1 \bigg( \frac{|x_{1,n}|^{k+2} + |\x_{1,n}|^{k+2}}{\x_{n+1} - x_{n+1}} + \frac{(t - \t)^{1+2/k}}{(\x_{n+1} - x_{n+1})^{2/k}} + 1 \bigg) \bigg) u(x,t);
\end{equation*}
	\item[{\it iii)}] if $k$ is odd, we have
\begin{eqnarray*}
    u(\x,\t) \!\!\!\! & \le & \!\!\!\! \exp \bigg( C_1 \bigg( \frac{|x_{1,n}- \x_{1,n}|^2}{t-\t} + \\
   & & \qquad + \frac{|\x_{n+1}-x_{n+1}- 
\frac{1}{2^{k+4}}   \sum_{j=1}^n(x_j^k+\x_j^k)(t-\t)|^{2/k}}{(t-\t)^{1 + 2/k}} + 1 \bigg) \bigg) u(x,t).
\end{eqnarray*}
\end{description}
\end{proposition}

\medskip \noindent{\it Proof.} The proof of  \emph{(i)} and \emph{(ii)} is analogous to that of  Proposition \ref{th-Harnack2-6}. The unique modification is due to the fact that here we have
\begin{equation} \label{e-g/2-666}
 x_{n+1}\left(\tfrac{t - \t}{2}\right) - x_{n+1} = \int_0^{\frac{t - \t}{4}} \sum_{j=1}^n \left( x_{j} +s \widetilde \o_j \right)^k d s + \int_0^{\frac{t - \t}{4}} \sum_{j=1}^n \left( \x_{j} - s \widehat \o_j \right)^k d s.
\end{equation}
instead of  \eqref{e-g/2-6666}. From \eqref{e-g/2-677}, by the same argument used in the proof of \eqref{e-corollary43-6} and \eqref{e-corollary43-6bis}, we obtain
\begin{equation*} 
  \frac{t - \t}{2^{k+4}} \sum_{j=1}^n \left( x_j^k + \x_j^k \right) \le x_{n+1}\left(\tfrac{t - \t}{2}\right) - x_{n+1} \le \frac{t-\t}{2} 15^k \sum_{j=1}^n \left( |x_j|^k + |\x_j|^k \right) + 2 n \left( \tfrac{3}{2} \right)^k (t - \t)^{k/2 + 1}.
\end{equation*}
We omit the other  details of the proof of \emph{(i)} and \emph{(ii)}.

\medskip

The proof of \emph{(iii)} follows from the same argument used in the proof of \emph{(i)}. Note that, as $k$ is odd, the function
\begin{equation*} 
	\varphi(b) = x_{n+1}\left(\tfrac{t - \t}{2}\right) + 2 \int_{0}^{\frac{t-\t}{4}}\sum_{j=1}^n \left( x_{j} +  b \, \bar \o_j \right)^k ds
\end{equation*}
defined for any $b \in \R$ is surjective, then in this case Lemma \ref{th-Harnack2-path} is sufficient to conclude the proof.  \hfill $\square$

\subsection*{Final derivation of the estimates}

\medskip \noindent{\it Proof of Theorem \ref{MTHM}.} It follows from the bounds proved in Sections \ref{upp-bd} and \ref{CHAINS_SECTION}. Consider first equation \eqref{PROC}.

The upper bound of \emph{(i)} is given in \eqref{EST_OFF_DIAG_DENS}. In order to prove the lower bound, we fix a constant vector $\omega \in \R^n, \alpha \in ]0,1[$, and we set
\begin{equation} \label{eq-xtilde}
    \tilde x_{1,n} = \xi_{1,n} - \alpha \sqrt{t} \, \omega, \qquad \tilde x_{n+1} = \xi_{n+1} - \int_0^{\alpha^2t} \Big| \x_{1,n} - \tfrac{s}{\alpha \sqrt{t}} \omega \Big|^k d s.
\end{equation}
According with Remark \ref{rem-Gamma-omo}, we have
\begin{equation} \label{eq-p-alpha-t}
	p({\alpha^2t}, \tilde x, \x) = \Gamma(\tilde x, {\alpha^2t}, \x, 0) = 
\frac{C}{{(\alpha^2t})^{1 + \frac{n+k}{2}}},  \qquad C = \Gamma\Big( \omega, - \tfrac{ |\omega|^k}{k+1}, 1 \Big) 
\end{equation}
For our purpose, we choose here $\omega = (1, \dots, 1) \in \R^n$, and $\alpha = 1/\sqrt{2}$.
Note that the constant $C$ is strictly positive, being $\omega \ne 0$. Also note that
\begin{equation} \label{eq-xn+1}
    0 < \xi_{n+1} - \tilde x_{n+1} \le 2^k \alpha^2 \left( t |\x_{1,n}|^k + \alpha^k t^{1 + k/2}|\omega|^k \right).
\end{equation}
We then apply Proposition \ref{th-Harnack2-6} \textit{(i)}, and we find
\begin{equation*}
\begin{split}
    p(t, x,\x) \ge \exp \bigg( - C_1 \bigg( & \frac{|x_{1,n}- \tilde x_{1,n}|^2}{t} + \\
    & \frac{\Big(\tilde x_{n+1}-x_{n+1}-\frac{t}{2^{k+5}}\big(|x_{1,n}|^k+|\tilde x_{1,n}|^k\big)\Big)^{2/k}\!\!\!\!}{t^{1 + 2/k}} + 1 \bigg)\!\! \bigg) p(t/2, \tilde x,\x).
\end{split}
\end{equation*}
The lower bound \textit{(i)} thus follows from the inequalities $|\tilde x_{1,n} - \xi_{1,n}| \le \sqrt{t}/2\, |\omega|$ and $\xi_{n+1} > \tilde x_{n+1}$.


The upper bound of \emph{(ii)} is equation \eqref{GAUSSIAN_BOUND}, the lower bound is given in Lemma \ref{MIN_CPT_MET}.

The upper bound of \emph{(iii)} is given in \eqref{MOMODERATE_BOUND}. In order to prove the lower bound we rely on Proposition \ref{th-Harnack2-6} \emph{(ii)}. In order to satisfy its hypothesis, we choose $\alpha \in ]0,1/\sqrt{2}]$ such that $\xi_{n+1} - \tilde x_{n+1} \le \frac12 (\xi_{n+1} - x_{n+1})$. Note that
\begin{equation} \label{eq-xn+2}
    0 < \xi_{n+1} - \tilde x_{n+1} \le 2^k \alpha^2 \left( t |\tilde x_{1,n}|^k + \alpha^k t^{1 + k/2}|\omega|^k \right),
\end{equation}
then it is sufficient to choose $\alpha$ such that
\begin{equation*}
    \alpha^2 t \le \frac{\xi_{n+1} - x_{n+1}}{2^{k+2} |\tilde x_{1,n}|^k}, \qquad
    ( \alpha \sqrt{t})^{k+2} \le \frac{\xi_{n+1} - x_{n+1}}{2^{k+2}|\omega|^k}.
\end{equation*}
Using our assumption $|\xi_{n+1} - x_{n+1}| \le K t^{1+k/2}$ we find
\begin{equation*}
    \frac{1}{\alpha^2 t} \ge
    \frac{C_{K,k}}{t} \max \left(\frac{|\tilde x_{1,n}|^k}{|\xi_{n+1} - x_{n+1}|^{\frac{k}{k+2}}}, 1 \right)
\end{equation*}
for some positive constant $C_{K,k}$ depending on $K$ and $k$. Then \eqref{eq-p-alpha-t} gives
$$
p({\alpha^2t}, \tilde x, \x) \ge \frac{\tilde C_{K,k}}{t^{\frac{n+k+2}{2}}} \left( 1 + \frac{|\tilde x_{1,n}|^{k+2}}{|\xi_{n+1} - x_{n+1}|} \right)^{{\frac{(k+2)(n+k+2)}{2k}}}
$$
With this choice of $\alpha$ Proposition \ref{th-Harnack2-6} \emph{(ii)} then gives
\begin{equation*}
\begin{split}
    p(t, x,\x) \ge & \frac{\tilde C_{K,k}}{t^{\frac{n+k+2}{2}}} \left( 1 + \frac{|\tilde x_{1,n}|^{k+2}}{|\xi_{n+1} - x_{n+1}|} \right)^{{\frac{(k+2)(n+k+2)}{2k}}} \cdot \\
    & \exp \bigg( - C_1 \bigg( \frac{|x_{1,n}|^{k+2} + |\tilde x_{1,n}|^{k+2}}{\tilde x_{n+1} - x_{n+1}} +
    \frac{((1- \alpha^2)t)^{1+2/k}}{(\tilde x_{n+1} - x_{n+1})^{2/k}} + 1 \bigg) \bigg)
\end{split}
\end{equation*}
and our lower bound follows from the inequalities $t/2 \le (1- \alpha^2) t \le t, 1/2(\x_{n+1} - x_{n+1}) \le \tilde x_{n+1} - x_{n+1} \le \x_{n+1} - x_{n+1}$ and $|\tilde x_{1,n} - \xi_{1,n}| \le \sqrt{t}/2\, |\omega|$.

When considering equation \eqref{PROC_2}, the lower bounds for points \textit{(i)} and \textit{(iii)} directly follow from \ref{th-Harnack2-7}. The proof of the remaining bounds can be done by the same arguments used for equation \eqref{PROC}. \hfill $\square$

\begin{remark}
We can assume by symmetry  that w.l.o.g. $|\xi_{1,n}|\ge |x_{1,n}| $.
In this case, observe from equation \eqref{eq-xn+1} that if $|\xi_{1,n}|\ge K t^{1/2} $ for $K$ large enough, then we can derive from
Lemma \ref{MIN_CPT_MET} that the Gaussian diagonal regime holds for the lower bound. From the proof leading to \eqref{EST_OFF_DIAG_DENS}, this means that the 
non-exponential estimates in case \textit{i)} could be alternatively rewritten changing the $t^{-\{\frac{n+k}2+1\}} $ term into
\begin{equation}
\label{LAST_EQ}
t^{-n/2}t^{-3/2}(\{ |x_{1,n}|^{k-1}+|\xi_{1,n}|^{k-1}\}\vee t^{\frac{k-1}2} )^{-1},
\end{equation}
that emphasizes the regime transition depending on the magnitude of the non-degenerate components w.r.t. to their characteristic time-scale.
From the above computations, the expression in \eqref{LAST_EQ} can also substitute the non-exponential term in case \textit{iii)}.
\end{remark}

\bigskip \noindent {\sc Acknowledgments.} \ We express our gratitude to Denis Talay, because this study has started in occasion of a workshop he organized, with the specific aim at fostering new research. We thank E. Lanconelli for his interest in our work. Also, the last author thanks University Paris 7 and the Laboratoire de Probabilit\'es et Mod\`eles Al\'eatoires for the one month invitation during which this work developed.

\bibliography{bib}

\end{document}

%% file: macro.tex
\def\I{{\rm{I\!\!I}}}
\def\leftB{[\![}
\def\rightB{]\!]}

\def\F{{\mathscr{ F}}}

\def\R{{\mathbb{R}} }
\def\N{{\mathbb{N}} }
\def\E{{\mathbb{E}}  }

\def\P{{\mathbb{P}}  }

\def\I{{\mathbb{I}}}
\def\det{{\rm{det}}}

\def\bint#1^#2{\displaystyle{\int_{#1}^{#2}}}
\def\bsum#1^#2{\displaystyle{\sum_{#1}^{#2}}}

\def\xdt_#1{X_#1(\Delta t)}

\newtheorem{PROP}{Proposition}[section]